# A Friendly
# Introduction to
# Differential Equations

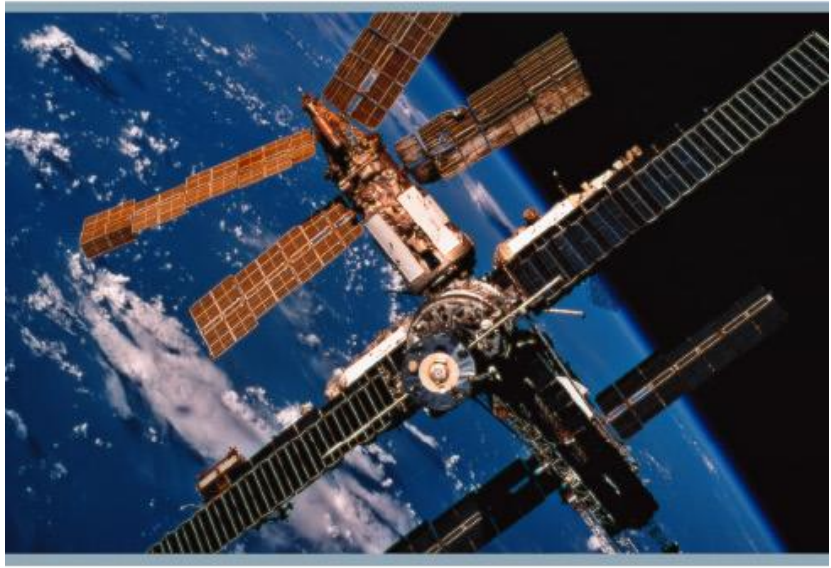

# First Edition

# Mohammed Kaabar

# A Friendly Introduction to Differential Equations

# A Friendly Introduction to Differential Equations

First Edition (Updated on December 28, 2015)

## By

# Mohammed K A Kaabar

# About the Author

**Mohammed Kaabar** has a Bachelor of Science in Theoretical Mathematics from Washington State University, Pullman, WA. He is a graduate student in Applied Mathematics at Washington State University, Pullman, WA, and he is a math tutor at the Math Learning Center (MLC) at Washington State University, Pullman. He is the author of *A First Course in Linear Algebra* Book, and his research interests are applied optimization, numerical analysis, differential equations, linear algebra, and real analysis. He was invited to serve as a Technical Program Committee (TPC) member in many conferences such as ICECCS 14, ENCINS 15, eQeSS 15, SSCC 15, ICSoEB 15, CCA 14, WSMEAP 14, EECSI 14, JIEEEC 13 and WCEEENG 12. He is an online instructor of two free online courses in numerical analysis: Introduction to Numerical Analysis and Advanced Numerical Analysis at Udemy Inc, San Francisco, CA. He is a former member of Institute of Electrical and Electronics Engineers (IEEE), IEEE Antennas and Propagation Society, IEEE Consultants Network, IEEE Smart Grid Community, IEEE Technical Committee on RFID, IEEE Life Sciences Community, IEEE Green ICT Community, IEEE Cloud Computing Community, IEEE Internet of Things Community, IEEE Committee on Earth Observations, IEEE Electric Vehicles Community, IEEE Electron Devices Society, IEEE Communications Society, and IEEE Computer Society. He also received several educational awards and certificates from accredited institutions. For more information about the author and his free online courses, please visit his personal website: http://www.mohammed-kaabar.net.

This book is available under a Creative Commons license.

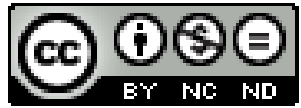

Note: If you find this book helpful and useful, and you see an opportunity to cite it in any of your publications, I would be very happy and appreciative of the citation. Please cite this book as:

Kaabar, M.K.: A Friendly Introduction to Differential Equations. Printed by CreateSpace, San Bernardino, CA (2015)

# Table of Contents

# Introduction

In this book, I wrote five chapters: The Laplace Transform, Systems of Homogenous Linear Differential Equations (HLDE), Methods of First and Higher Orders Differential Equations, Extended Methods of First and Higher Orders Differential Equations, and Applications of Differential Equations. I also added exercises at the end of each chapter above to let students practice additional sets of problems other than examples, and they can also check their solutions to some of these exercises by looking at "Answers to Odd-Numbered Exercises" section at the end of this book. This book is a very useful for college students who studied Calculus II, and other students who want to review some concepts of differential equations before studying courses such as partial differential equations, applied mathematics, and electric circuits II. According to my experience as a math tutor, I follow the steps of my Professor Ayman Badawi who taught me Differential Equations and Linear Algebra at the American University of Sharjah, and I will start with laplace transforms as the first chapter of this book [3]. I used what I learned from Dr. Badawi to write the content of this book depending on his class notes as my main reference [3]. If you have any comments related to the contents of this book, please email your comments to mkaabar@math.wsu.edu.

I wish to express my gratitude and appreciation to my father, my mother, and my only lovely 13-year old brother who is sick, and I want to spend every dollar in his heath care. I would also like to give a special thanks to my Professor Dr. Ayman Badawi who supported me in every successful achievement I made, and I would like to thank all administrators and professors of mathematics at WSU for their educational support. In conclusion, I would appreciate to consider this book as a milestone for developing more math books that can serve our mathematical society in the area of differential equations.

Mohammed K A Kaabar

# Table of Laplace Transform

| | | |
|---|---|---|
| 1 | $\mathcal{L}\{f(t)\} = \int_{0}^{\infty} e^{-st} f(t)\,dt$ | |
| 2 | $\mathcal{L}\{1\} = \dfrac{1}{s}$ | $\mathcal{L}\{t^m\} = \dfrac{m!}{s^{m+1}}$ where m is a positive integer (whole number) |
| 3 | $\mathcal{L}\{\sin ct\} = \dfrac{c}{s^2 + c^2}$ | $\mathcal{L}\{\cos ct\} = \dfrac{s}{s^2 + c^2}$ |
| 4 | $\mathcal{L}\{e^{bt}\} = \dfrac{1}{s - b}$ | |
| 5 | $\mathcal{L}\{e^{bt} f(t)\} = F(s)\,|s \to s - b$ | $\mathcal{L}^{-1}\{F(s)\,|s \to s - b\,\} = e^{bt} f(t)$ |
| 6 | $\mathcal{L}\{h(t)U(t - b)\} = e^{-bs}\mathcal{L}\{h(t + b)\}$ $\mathcal{L}^{-1}\{e^{-bs} F(s)\} = f(t - b)U(t - b)$ | |
| 7 | $\mathcal{L}\{U(t - b)\} = \dfrac{e^{-bs}}{s}$ | $\mathcal{L}^{-1}\left\{\dfrac{e^{-bs}}{s}\right\} = U(t - b)$ |
| 8 | $\mathcal{L}\{f^{(m)}(t)\} = s^m F(s) - s^{m-1} f(0) - s^{m-2} f'(s) - \cdots - f^{(m-1)}(s)$ where m is a positive integer | |
| 9 | $\mathcal{L}\{t^m f(t)\}(s) = (-1)^m \dfrac{d^m F(s)}{ds^m}$ where m is a positive integer | $\mathcal{L}^{-1}\left\{\dfrac{d^m F(s)}{ds^m}\right\} = (-1)^m t^m f(t)$ where m is a positive integer |
| 10 | $\mathcal{L}\{f(t) * h(t)\} = F(s) \cdot H(s)$ | $f(t) * h(t) = \int_{0}^{t} f(\psi) h(t - \psi)\,d\psi$ |
| 11 | $\mathcal{L}\left\{\int_{0}^{t} f(\psi)\,d\psi\right\} = \dfrac{F(s)}{s}$ | $\mathcal{L}^{-1}\{F(s) \cdot H(s)\} = f(t) * h(t)$ |
| 12 | $\mathcal{L}\{\delta(t)\} = 1$ | $\mathcal{L}\{\delta(t - b)\} = e^{-bs}$ |
| 13 | Assume that $f(t)$ is periodic with period $P$, then: $\mathcal{L}\{f(t)\} = \dfrac{1}{1 - e^{-Ps}} \int_{0}^{P} e^{-st} f(t)\,dt$ | |

Table 1.1.1: Laplace Transform

*This page intentionally left blank*

# Chapter 1

# The Laplace Transform

In this chapter, we start with an introduction to Differential Equations (DEs) including linear DEs, nonlinear DEs, independent variables, dependent variables, and the order of DEs. Then, we define the laplace transforms, and we give some examples of Initial Value Problems (IVPs). In addition, we discuss the inverse laplace transforms. We cover in the remaining sections an important concept known as the laplace transforms of derivatives, and we mention some properties of laplace transforms. Finally, we learn how to solve systems of linear equations (LEs) using Cramer's Rule.

## 1.1 Introductions to Differential Equations

In this section, we are going to discuss how to determine whether the differential equation is linear or nonlinear, and we will find the order of differential equations. At the end of this section, we will show the purpose of differential equations.

let's start with the definition of differential equation and with a simple example about differential equation.

**Definition 1.1.1** A mathematical equation is called differential equation if it has two types of variables: dependent and independent variables where the dependent variable can be written in terms of independent variable.

**Example 1.1.1** Given that $y' = 15x$.

   a) Find $y$. (Hint: Find the general solution of $y'$)
   b) Determine whether $y' = 15x$ is a linear differential equation or nonlinear differential equation. Why?
   c) What is the order of this differential equation?

**Solution: Part a:** To find $y$, we need to find the general solution of $y'$ by taking the integral of both sides as follows:

$$\int y' \, dy = \int 15x \, dx$$

Since $\int y' \, dy = y$ because the integral of derivative function is the original function itself (In general, $\int f'(x) \, dx = f(x)$), then $y = \int 15x \, dx = \frac{15}{2}x^2 + c = 7.5x^2 + c$. Thus, the general solution is the following: $y(x) = 7.5x^2 + c$ where $c$ is constant. This means that $y$ is called dependent variable because it depends on $x$, and $x$ is called independent variable because it is independent from $y$.

**Part b:** To determine whether $y' = 15x$ is a linear differential equation or nonlinear differential equation, we need to introduce the following definition:

**Definition 1.1.2** The differential equation is called linear if the dependent variable and all its derivatives are to the power 1. Otherwise, the differential equation is nonlinear.

According to the above question, we have the following: $y(x) = 7.5x^2 + c$ where $c$ is constant. Since the dependent variable $y$ and all its derivatives are to the power 1, then using definition 1.1.2, this differential equation is linear.

**Part c:** To find the order of this differential equation, we need to introduce the following definition:

**Definition 1.1.3** The order of differential equation is the highest derivative in the equation (i.e. The order of $y''' + 3y'' + 2y' = 12x^2 + 22$ is 3).

Using definition 1.1.3, the order of $y' = 15x$ is 1.

**Example 1.1.2** Given that $z''' + 2z'' + z' = 2x^3 + 22$.

   a) Determine whether $z''' + 2z'' + z' = 2x^3 + 22$ is a linear differential equation or nonlinear differential equation. Why?
   b) What is the order of this differential equation?

**Solution: Part a:** Since $z$ is called dependent variable because it depends on $x$, and $x$ is called independent variable because it is independent from $z$, then to determine whether $z''' + 2z'' + z' = 2x^3 + 22$ is a linear differential equation or nonlinear differential equation, we need to use definition 1.1.2 as follows: Since the dependent variable $z$ and all its derivatives are to the power 1, then this differential equation is linear.

**Part b:** To find the order of $z''' + 2z'' + z' = 2x^3 + 22$, we use definition 1.1.3 which implies that the order is 3 because the highest derivative is 3.

**Example 1.1.3** Given that $m^{(4)} + (3m'')^3 - m = \sqrt{x+1}$.
(Hint: Do not confuse between $m^{(4)}$ and $m^4$ because $m^{(4)}$ means the fourth derivative of $m$, while $m^4$ means the fourth power of m).

  a) Determine whether $m^{(4)} + (3m'')^3 - m = \sqrt{x+1}$ is a linear differential equation or nonlinear differential equation. Why?
  b) What is the order of this differential equation?

**Solution: Part a:** Since $m$ is called dependent variable because it depends on $x$, and $x$ is called independent variable because it is independent from $m$, then to determine whether $m^{(4)} + (3m'')^3 - m = \sqrt{x+1}$ is a linear differential equation or nonlinear differential equation, we need to use definition 1.1.2 as follows: Since the dependent variable $m$ and all its derivatives are not to the power 1, then this differential equation is nonlinear.

**Part b:** To find the order of $m^{(4)} + (3m'')^3 - m = \sqrt{x+1}$, we use definition 1.1.3 which implies that the order is 4 because the highest derivative is 4.

The following are some useful notations about differential equations:
$z^{(m)}$ is the $m^{\text{th}}$ derivative of $z$.
$z^m$ is the $m^{\text{th}}$ power of $z$.

The following two examples are a summary of this section:

**Example 1.1.4** Given that
$2xb^{(3)} + (x+1)b^{(2)} + 3b = x^2 e^x$. Determine whether it is a linear differential equation or nonlinear differential equation.

**Solution:** To determine whether it is a linear differential equation or nonlinear differential equation, We need to apply what we have learned from the previous examples in the following five steps:

**Step 1:** $b$ is a dependent variable, and $x$ is an independent variable.

**Step 2:** Since $b$ and all its derivatives are to the power 1, then the above differential equation is linear.

**Step 3:** Coefficients of $b$ and all its derivatives are in terms of the independent variable $x$.

**Step 4:** Assume that $C(x) = x^2 e^x$. Then, $C(x)$ must be in terms of $x$.

**Step 5:** Our purpose from the above differential equation is to find a solution where $b$ can be written in term of $x$.

Thus, the above differential equation is a linear differential equation of order 3.

**Example 1.1.5** Given that

$(w^2 + 1)h^{(4)} - 3wh' = w^2 + 1$. Determine whether it is a linear differential equation or nonlinear differential equation.

**Solution:** To determine whether it is a linear differential equation or nonlinear differential equation, We need to apply what we have learned from the previous examples in the following five steps:

**Step 1:** $h$ is a dependent variable, and $w$ is an independent variable.

**Step 2:** Since $h$ and all its derivatives are to the power 1, then the above differential equation is linear.

**Step 3:** Coefficients of $h$ and all its derivatives are in terms of the independent variable $w$.

**Step 4:** Assume that $C(w) = w^2 + 1$. Then, $C(w)$ must be in terms of $w$.

**Step 5:** Our purpose from the above differential equation is to find a solution where $h$ can be written in term of $w$.

Thus, the above differential equation is a linear differential equation of order 4.

# 1.2 Introductions to the Laplace Transforms

In this section, we are going to introduce the definition of the laplace transforms in general, and how can we use this definition to find the laplace transform of any function. Then, we will give several examples about the laplace transforms, and we will show how the table 1.1.1 is helpful to find laplace transforms.

**Definition 1.2.1** the laplace transform, denoted by $\mathcal{L}$, is defined in general as follows:

$$\mathcal{L}\{f(x)\} = \int_0^\infty f(x)e^{-sx}dx$$

**Example 1.2.1** Using definition 1.2.1, find $\mathcal{L}\{1\}$.

**Solution:** To find $\mathcal{L}\{1\}$ using definition 1.2.1, we need to do the following steps:

**Step 1:** We write the general definition of laplace transform as follows:

$$\mathcal{L}\{f(x)\} = \int_0^\infty f(x)e^{-sx}dx$$

**Step 2:** Here in this example, $f(x) = 1$ because $\mathcal{L}\{f(x)\} = \mathcal{L}\{1\}$.

**Step 3:**

$$\mathcal{L}\{1\} = \int\limits_0^\infty (1)e^{-sx}dx$$

By the definition of integral, we substitute $\int_0^\infty (1)e^{-sx}dx$ with $\lim_{b \to \infty} \int_0^b (1)e^{-sx}dx$.

**Step 4:** We need to find $\lim_{b \to \infty} \int_0^b \boxed{e^{-sx}dx}$ as follows:

It is easier to find what it is inside the above box $\left(\int_0^b e^{-sx}dx\right)$, and after that we can find the limit of $\int_0^b e^{-sx}dx$.

Thus, $\int_0^b e^{-sx}dx = -\frac{1}{s}e^{-sx}\Big|_{x=0}^{x=b} = -\frac{1}{s}e^{-sb} + \frac{1}{s}e^{-s(0)} = -\frac{1}{s}e^{-sb} + \frac{1}{s}$.

**Step 5:** We need find the limit of $-\frac{1}{s}e^{-sb} + \frac{1}{s}$ as follows:

$\lim_{b \to \infty}\left(-\frac{1}{s}e^{-sb} + \frac{1}{s}\right) = \frac{1}{s}$ where $s > 0$.

To check if our answer is right, we need to look at table 1.1.1 at the beginning of this book. According to table 1.1.1 section 2, we found $\mathcal{L}\{1\} = \frac{1}{s}$ which is the same answer we got. Thus, we can conclude our example with the following fact:

**Fact 1.2.1** $\mathcal{L}\{any\ constant, say\ m\} = \frac{m}{s}$.

**Example 1.2.2** Using definition 1.2.1, find $\mathcal{L}\{e^{7x}\}$.

**Solution:** To find $\mathcal{L}\{e^{7x}\}$ using definition 1.2.1, we need to do the following steps:

**Step 1:** We write the general definition of laplace transform as follows:

$$\mathcal{L}\{f(x)\} = \int\limits_0^\infty f(x)e^{-sx}dx$$

**Step 2:** Here in this example, $f(x) = e^{7x}$ because $\mathcal{L}\{f(x)\} = \mathcal{L}\{e^{7x}\}$.

**Step 3:**

$$\mathcal{L}\{e^{7x}\} = \int\limits_0^\infty (e^{7x})e^{-sx}dx$$

By the definition of integral, we substitute $\int_0^\infty (e^{7x})e^{-sx}dx$ with $\lim_{b\to\infty} \int_0^b (e^{7x})e^{-sx}dx$.

**Step 4:** We need to find $\lim_{b\to\infty} \int_0^b \boxed{(e^{7x})e^{-sx}dx \text{ as}}$ follows:
It is easier to find what it is inside the above box $\left(\int_0^b (e^{7x})e^{-sx}dx\right)$, and after that we can find the limit of $\int_0^b (e^{7x})e^{-sx}dx$.

Thus, $\int_0^b (e^{7x})e^{-sx}dx = \int_0^b e^{(7-s)x}dx = \frac{1}{(7-s)}e^{(7-s)x}\Big|_{x=0}^{x=b} = \frac{1}{(7-s)}e^{(7-s)b} - \frac{1}{(7-s)}e^{(7-s)(0)} = \frac{1}{(7-s)}e^{(7-s)b} - \frac{1}{(7-s)}$.

**Step 5:** We need find the limit of $\frac{1}{(7-s)}e^{(7-s)b} - \frac{1}{(7-s)}$ as follows:

$\lim_{b\to\infty} \left(\frac{1}{(7-s)}e^{(7-s)b} - \frac{1}{(7-s)}\right) = -\frac{1}{(7-s)} = \frac{1}{(s-7)}$ where $s > 7$.
To check if our answer is right, we need to look at table 1.1.1 at the beginning of this book. According to table 1.1.1 section 4, we found $\mathcal{L}\{e^{7x}\} = \frac{1}{s-7}$ which is the same answer we got.

The following examples are two examples about finding the integrals to review some concepts that will help us finding the laplace transforms.

**Example 1.2.3** Find $\int x^3 e^{2x}dx$.

**Solution:** To find $\int x^3 e^{2x} dx$, it is easier to use a method known as the table method than using integration by parts. In the table method, we need to create two columns: one for derivatives of $x^3$, and the other one for integrations of $e^{2x}$. Then, we need to keep deriving $x^3$ till we get zero, and we stop integrating when the corresponding row is zero. The following table shows the table method to find $\int x^3 e^{2x} dx$:

| Derivatives Part | Integration Part |
|:---:|:---:|
| $x^3$ | $e^{2x}$ |
| $3x^2$ | $\frac{1}{2}e^{2x}$ |
| $6x$ | $\frac{1}{4}e^{2x}$ |
| $6$ | $\frac{1}{8}e^{2x}$ |
| $0$ | $\frac{1}{16}e^{2x}$ |

Table 1.2.1: Table Method for $\int x^3 e^{2x} dx$

We always start with positive sign, followed by negative sign, and so on as we can see in the above table 1.2.1. Now, from the above table 1.2.1, we can find $\int x^3 e^{2x} dx$ as follows:

$$\int x^3 e^{2x} dx = \frac{1}{2}x^3 e^{2x} - \frac{1}{4}(3)x^2 e^{2x} + \frac{1}{8}(6)x e^{2x} - \frac{1}{16}(6)e^{2x} + C$$

Thus, $\int x^3 e^{2x} dx = \frac{1}{2}x^3 e^{2x} - \frac{3}{4}x^2 e^{2x} + \frac{3}{4}x e^{2x} - \frac{3}{8}e^{2x} + C$.

In conclusion, we can always use the table method to find integrals like $\int (polynomial)e^{ax} dx$ and $\int (polynomial)\sin(ax) dx$.

**Example 1.2.4** Find $\int 3x^2 \sin(4x) dx$.

**Solution:** To find $\int 3x^2 \sin(4x) dx$, it is easier to use the table method than using integration by parts. In the table method, we need to create two columns: one for derivatives of $3x^2$, and the other one for integrations of $\sin(4x)$. Then, we need to keep deriving $3x^2$ till we get zero, and we stop integrating when the corresponding row is zero. The following table shows the table method to find $\int 3x^2 \sin(4x)dx$:

| Derivatives Part | Integration Part |
|---|---|
| $3x^2$ | $\sin(4x)$ |
| $6x$ | $-\frac{1}{4}\cos(4x)$ ▬ |
| $6$ | $-\frac{1}{16}\sin(4x)$ ▬ |
| $0$ | $\frac{1}{64}\cos(4x)$ ▬ |

Table 1.2.2: Table Method for $\int 3x^2 \sin(4x)dx$

We always start with positive sign, followed by negative sign, and so on as we can see in the above table 1.2.2. Now, from the above table 1.2.2, we can find $\int 3x^2 \sin(4x)dx$ as follows:

$$\int 3x^2 \sin(4x)dx$$
$$= -\frac{1}{4}(3)x^2\cos(4x) - \left(-\frac{1}{16}\right)6x\sin(4x)$$
$$+ \left(\frac{1}{64}\right)6\cos(4x) + C$$

Thus, $\int 3x^2 \sin(4x)dx = -\frac{3}{4}x^2\cos(4x) + \frac{3}{8}x\sin(4x) + \frac{3}{32}\cos(4x) + C$.

**Example 1.2.5** Using definition 1.2.1, find $\mathcal{L}\{x^2\}$.

**Solution:** To find $\mathcal{L}\{x^2\}$ using definition 1.2.1, we need to do the following steps:

**Step 1:** We write the general definition of laplace transform as follows:

$$\mathcal{L}\{f(x)\} = \int\limits_0^\infty f(x)e^{-sx}dx$$

**Step 2:** Here in this example, $f(x) = x^2$ because $\mathcal{L}\{f(x)\} = \mathcal{L}\{x^2\}$.

**Step 3:** $\mathcal{L}\{x^2\} = \int_0^\infty (x^2)e^{-sx}dx$

By the definition of integral, we substitute $\int_0^\infty (x^2)e^{-sx}dx$ with $\lim_{b\to\infty}\int_0^b (x^2)e^{-sx}dx$.

**Step 4:** We need to find $\lim_{b\to\infty}\int_0^b \boxed{(x^2)e^{-sx}dx}$ as follows:

It is easier to find what it is inside the above box $\left(\int_0^b (x^2)e^{-sx}dx\right)$, and after that we can find the limit of $\int_0^b (x^2)e^{-sx}dx$.

Now, we need to find $\int_0^b (x^2)e^{-sx}dx$ using the table method. In the table method, we need to create two columns: one for derivatives of $x^2$, and the other one for integrations of $e^{-sx}$. Then, we need to keep deriving $x^2$ till we get zero, and we stop integrating when the corresponding row is zero. The following table shows the table method to find $\int (x^2)e^{-sx}dx$:

| Derivatives Part | Integration Part |
|:---:|:---:|
| $x^2$ | $e^{-sx}$ |
| $2x$ | $-\frac{1}{s}e^{-sx}$ |
| $2$ | $\frac{1}{s^2}e^{-sx}$ |
| $0$ | $-\frac{1}{s^3}e^{-sx}$ |

Table 1.2.3: Table Method for $\int (x^2)e^{-sx}dx$

We always start with positive sign, followed by negative sign, and so on as we can see in the above table 1.2.3. Now, from the above table 1.2.3, we can find $\int 3x^2\sin(4x)dx$ as follows:

Thus, $\int (x^2)e^{-sx}dx = -\frac{1}{s}x^2 e^{-sx} - \frac{2}{s^2}xe^{-sx} - \frac{2}{s^3}e^{-sx} + C$.

Now, we need to evaluate the above integral from 0 to $b$ as follows:

$\int_0^b (x^2) e^{-sx} dx = -\frac{1}{s} x^2 e^{-sx} - \frac{2}{s^2} x e^{-sx} - \frac{2}{s^3} e^{-sx} \Big|_{x=0}^{x=b} = \left( -\frac{1}{s} b^2 e^{-sb} - \frac{2}{s^2} b e^{-sb} - \frac{2}{s^3} e^{-sb} + \frac{2}{s^3} \right).$

**Step 5:** We need find the limit of $\left( -\frac{1}{s} b^2 e^{-sb} - \frac{2}{s^2} b e^{-sb} - \frac{2}{s^3} e^{-sb} + \frac{2}{s^3} \right)$ as follows:

$\lim_{b \to \infty} \left( -\frac{1}{s} b^2 e^{-sb} - \frac{2}{s^2} b e^{-sb} - \frac{2}{s^3} e^{-sb} + \frac{2}{s^3} \right) = \lim_{b \to \infty} \left( 0 + \frac{2}{s^3} \right) = \frac{2}{s^3}$ where $s > 0$.

To check if our answer is right, we need to look at table 1.1.1 at the beginning of this book. According to table 1.1.1 section 2 at the right side, we found $\mathcal{L}\{x^2\} = \frac{2!}{s^{2+1}} = \frac{2}{s^3}$ which is the same answer we got.

**Example 1.2.6** Using table 1.1.1, find $\mathcal{L}\{\sin(5x)\}$.

**Solution:** To find $\mathcal{L}\{\sin(5x)\}$ using table 1.1.1, we need to do the following steps:
**Step 1:** We look at the transform table (table 1.1.1).
**Step 2:** We look at which section in table 1.1.1 contains *sin* function.
**Step 3:** We write down what we get from table 1.1.1 (Section 3 at the left side) as follows:

$\mathcal{L}\{\sin ct\} = \frac{c}{s^2 + c^2}$

**Step 4:** We change what we got from step 3 to make it look like $\mathcal{L}\{\sin(5x)\}$ as follows:

$\mathcal{L}\{\sin 5x\} = \frac{5}{s^2 + 5^2} = \frac{5}{s^2 + 25}.$

Thus, $\mathcal{L}\{\sin 5x\} = \frac{5}{s^2 + 5^2} = \frac{5}{s^2 + 25}.$

**Example 1.2.7** Using table 1.1.1, find $\mathcal{L}\{\cos(-8x)\}$.

**Solution:** To find $\mathcal{L}\{\cos(-8x)\}$ using table 1.1.1, we need to do the following steps:
**Step 1:** We look at the transform table (table 1.1.1).

**Step 2:** We look at which section in table 1.1.1 contains *cos* function.

**Step 3:** We write down what we get from table 1.1.1 (Section 3 at the right side) as follows:

$\mathcal{L}\{\cos ct\} = \dfrac{s}{s^2 + c^2}$

**Step 4:** We change what we got from step 3 to make it look like $\mathcal{L}\{\cos(-8x)\}$ as follows:

$\mathcal{L}\{\cos(-8x)\} = \dfrac{s}{s^2 + (-8)^2} = \dfrac{s}{s^2 + 64}$.

Thus, $\mathcal{L}\{\cos(-8x)\} = \dfrac{s}{s^2 + (-8)^2} = \dfrac{s}{s^2 + 64}$.

We will give some important mathematical results about laplace transforms.

**Result 1.2.1** Assume that $c$ is a constant, and $f(x)$, $g(x)$ are functions. Then, we have the following: (Hint: $F(s)$ and $G(s)$ are the laplace transforms of $f(x)$ and $g(x)$, respectively).

  a) $\mathcal{L}\{g(x)\} = G(s)$ (i.e. $\mathcal{L}\{8\} = \dfrac{8}{s} = G(s)$ where $g(x) = 8$).

  b) $\mathcal{L}\{c \cdot g(x)\} = c \cdot \mathcal{L}\{g(x)\} = c \cdot G(s)$.

  c) $\mathcal{L}\{f(x) \mp g(x)\} = F(s) \mp G(s)$.

  d) $\mathcal{L}\{f(x) \cdot g(x)\}$ is not necessary equal to $F(s) \cdot G(s)$.

**Example 1.2.8** Using definition 1.2.1, find $\mathcal{L}\{y'\}$.

**Solution:** To find $\mathcal{L}\{y'\}$ using definition 1.2.1, we need to do the following steps:

**Step 1:** We write the general definition of laplace transform as follows:

$$\mathcal{L}\{f(x)\} = \int\limits_{0}^{\infty} f(x)e^{-sx}dx$$

**Step 2:** Here in this example, $f(x) = y'$ because $\mathcal{L}\{f(x)\} = \mathcal{L}\{y'\}$.

**Step 3:** $\mathcal{L}\{y'\} = \int_{0}^{\infty}(y')e^{-sx}dx$

By the definition of integral, we substitute $\int_0^\infty (y')e^{-sx}dx$ with $\lim_{b \to \infty} \int_0^b (y')e^{-sx}dx$.

**Step 4:** We need to find $\lim_{b \to \infty} \int_0^b \boxed{y'e^{-sx}dx}$ as follows:

It is easier to find what it is inside the above box $\left( \int_0^b y'e^{-sx}dx \right)$, and after that we can find the limit of $\int_0^b y'e^{-sx}dx$.

To find $\int_0^b y'e^{-sx}dx$, we need to use integration by parts as follows:

$u = e^{-sx}$ $\qquad\qquad\qquad\qquad dv = y'dx$

$du = -se^{-sx}dx$ $\qquad\qquad\quad v = \int dv = \int y' \, dx = y$

$\int u \, dv = uv - \int v \, du$

$\int y'e^{-sx}dx = ye^{-sx} - \int y(-se^{-sx}) \, dx$

Now, we can find the limit as follows:

$\lim_{b \to \infty} \int_0^b y'e^{-sx}dx = \lim_{b \to \infty} ye^{-sx} - \lim_{b \to \infty} \int_0^b y(-se^{-sx}) \, dx$

$\lim_{b \to \infty} \int_0^b y'e^{-sx}dx = \lim_{b \to \infty} ye^{-sx}\Big|_{x=0}^{x=b} + \int_0^\infty y(se^{-sx}) \, dx$

Because $\int_0^\infty y(se^{-sx}) \, dx = \lim_{b \to \infty} \int_0^b y(-se^{-sx}) \, dx$.

$\lim_{b \to \infty} \int_0^b y'e^{-sx}dx = \lim_{b \to \infty} (ye^{-sb} - ye^{-s(0)}) + \int_0^\infty y(se^{-sx}) \, dx$

$\lim_{b \to \infty} \int_0^b y'e^{-sx}dx = \lim_{b \to \infty} \left( y(b)e^{-sb} - y(0)e^{-s(0)} \right) + s\int_0^\infty ye^{-sx} \, dx$

Since we have $\int_0^\infty ye^{-sx} \, dx$, then using result 1.2.1 $\left( \int_0^\infty ye^{-sx} \, dx = Y(s) \right)$, we obtain the following:

$$\lim_{b \to \infty} \int_0^b y' e^{-sx} dx = \lim_{b \to \infty} \left( y(b) e^{-sb} - y(0) \right) + s \int_0^\infty y e^{-sx} dx$$

$$\int_0^b y' e^{-sx} dx = \lim_{b \to \infty} \left( y(b) e^{-sb} - y(0) \right) + sY(s)$$

$$\int_0^b y' e^{-sx} dx = \left( y(b) e^{-s(\infty)} - y(0) \right) + sY(s)$$

$$\int_0^b y' e^{-sx} dx = (0 - y(0)) + sY(s)$$

$$\int_0^b y' e^{-sx} dx = sY(s) - y(0)$$

Thus, $\mathcal{L}\{y'\} = sY(s) - y(0)$.

We conclude this example with the following results:

**Result 1.2.2** Assume that $f(x)$ is a function, and $F(s)$ is the laplace transform of $f(x)$. Then, we have the following:

   a) $\mathcal{L}\{f'(x)\} = \mathcal{L}\{f^{(1)}(x)\} = sF(s) - f(0)$.
   b) $\mathcal{L}\{f''(x)\} = \mathcal{L}\{f^{(2)}(x)\} = s^2 F(s) - sf(0) - f'(0)$.
   c) $\mathcal{L}\{f'''(x)\} = \mathcal{L}\{f^{(3)}(x)\} = s^3 F(s) - s^2 f(0) - sf'(0) - f''(0)$.
   d) $\mathcal{L}\{f^{(4)}(x)\} = s^4 F(s) - s^3 f(0) - s^2 f'(0) - sf''(0) - f'''(0)$.

For more information about this result, it is recommended to look at section 8 in table 1.1.1. If you look at it, you will find the following:
$\mathcal{L}\{f^{(m)}(t)\} = s^m F(s) - s^{m-1} f(0) - s^{m-2} f'(s) - \cdots - f^{(m-1)}(s)$.

**Result 1.2.3** Assume that $c$ is a constant, and $f(x)$ is a function where $F(s)$ is the laplace transform of $f(x)$.

Then, we have the following: $\mathcal{L}\{cf(x)\} = c\mathcal{L}\{f(x)\} = cF(s)$.

# 1.3 Inverse Laplace Transforms

In this section, we will discuss how to find the inverse laplace transforms of different types of mathematical functions, and we will use table 1.1.1 to refer to the laplace transforms.

**Definition 1.3.1** The inverse laplace transform, denoted by $\mathcal{L}^{-1}$, is defined as a reverse laplace transform, and to find the inverse laplace transform, we need to think about which function has a laplace transform that equals to the function in the inverse laplace transform. For example, suppose that $f(x)$ is a function where $F(s)$ is the laplace transform of $f(x)$. Then, the inverse laplace transform is $\mathcal{L}^{-1}\{F(s)\} = f(x)$. (i.e. $\mathcal{L}^{-1}\left\{\frac{1}{s}\right\}$ we need to think which function has a laplace transform that equals to $\frac{1}{s}$, in this case the answer is 1).

**Example 1.3.1** Find $\mathcal{L}^{-1}\left\{\frac{34}{s}\right\}$.

**Solution:** First of all, $\mathcal{L}^{-1}\left\{\frac{34}{s}\right\} = \mathcal{L}^{-1}\left\{34\left(\frac{1}{s}\right)\right\}$.
Using definition 1.3.1 and table 1.1.1, the answer is 34.

**Example 1.3.2** Find $\mathcal{L}^{-1}\left\{\frac{s}{s^2+2}\right\}$.

**Solution:** Using definition 1.3.1 and the right side of section 3 in table 1.1.1, we find the following:

$\mathcal{L}\{\cos ct\} = \dfrac{s}{s^2 + c^2}$

$\dfrac{s}{s^2+2}$ is an equivalent to $\dfrac{s}{s^2+(\sqrt{2})^2}$

Since $\mathcal{L}\{\cos\sqrt{2}t\} = \dfrac{s}{s^2+(\sqrt{2})^2}$ , then this means that

$\mathcal{L}^{-1}\left\{\dfrac{s}{s^2+2}\right\} = \cos\sqrt{2}t.$

**Example 1.3.3** Find $\mathcal{L}^{-1}\left\{\dfrac{-3}{s^2+9}\right\}$.

**Solution:** Using definition 1.3.1 and the left side of section 3 in table 1.1.1, we find the following:

$\mathcal{L}\{\sin ct\} = \dfrac{c}{s^2 + c^2}$

$\dfrac{-3}{s^2+9}$ is an equivalent to $\dfrac{-3}{s^2+(-3)^2}$

Since $\mathcal{L}\{\sin(-3)t\} = \mathcal{L}\{\sin(-3t)\} = \dfrac{-3}{s^2+(-3)^2}$ , then this

means that $\mathcal{L}^{-1}\left\{\dfrac{-3}{s^2+9}\right\} = \sin(-3t).$

**Example 1.3.4** Find $\mathcal{L}^{-1}\left\{\dfrac{1}{s+8}\right\}$.

**Solution:** Using definition 1.3.1 and section 4 in table 1.1.1, we find the following:

$\mathcal{L}\{e^{bt}\} = \dfrac{1}{s-b}$ .

$\dfrac{1}{s+8}$ is an equivalent to $\dfrac{1}{s-(-8)}$

Since $\mathcal{L}\{e^{-8t}\} = \dfrac{1}{s+8}$ , then this means that $\mathcal{L}^{-1}\left\{\dfrac{1}{s+8}\right\} = e^{-8t}.$

**Result 1.3.1** Assume that $c$ is a constant, and $f(x)$ is a function where $F(s)$ and $G(s)$ are the laplace transforms of $f(x),$ and $g(x),$ respectively.

    a) $\mathcal{L}^{-1}\{cF(s)\} = c\mathcal{L}^{-1}\{F(s)\} = cf(x).$

    b) $\mathcal{L}^{-1}\{F(s) \mp G(s)\} = \mathcal{L}^{-1}\{F(s)\} \mp \mathcal{L}^{-1}\{G(s)\} = f(x) \mp g(x).$

**Example 1.3.5** Find $\mathcal{L}^{-1}\left\{\dfrac{5}{s^2+4}\right\}$.

**Solution:** Using result 1.3.1, $\mathcal{L}^{-1}\left\{\frac{5}{s^2+4}\right\} = 5\mathcal{L}^{-1}\left\{\frac{1}{s^2+4}\right\}$. Now, by using the left side of section 3 in table 1.1.1, we need to make it look like $\frac{c}{s^2+c^2}$ because $\mathcal{L}\{\sin ct\} = \frac{c}{s^2+c^2}$.

Therefore, we do the following: $\frac{5}{2}\mathcal{L}^{-1}\left\{(2)\frac{1}{s^2+4}\right\} = \frac{5}{2}\mathcal{L}^{-1}\left\{\frac{2}{s^2+4}\right\}$, and $\frac{2}{s^2+4}$ is an equivalent to $\frac{2}{s^2+(2)^2}$.

Since $\frac{5}{2}\mathcal{L}\{\sin(2)t\} = \mathcal{L}\left\{\frac{5}{2}\sin(2t)\right\} = \frac{5}{2}\left(\frac{2}{s^2+(2)^2}\right) = \frac{5}{s^2+4}$ , then by using definition 1.3.1, this means that $\mathcal{L}^{-1}\left\{\frac{5}{s^2+4}\right\} = \frac{5}{2}\sin(2t)$.

**Example 1.3.6** Find $\mathcal{L}^{-1}\left\{\frac{7}{2s-3}\right\}$.

**Solution:** Using result 1.3.1, $\mathcal{L}^{-1}\left\{\frac{7}{2s-3}\right\} = 7\mathcal{L}^{-1}\left\{\frac{1}{2s-3}\right\}$. Now, by using section 4 in table 1.1.1, we need to make it look like $\frac{1}{s-b}$ because $\mathcal{L}\{e^{bt}\} = \frac{1}{s-b}$.

Therefore, we do the following: $7\mathcal{L}^{-1}\left\{\frac{1}{2s-3}\right\} = 7\mathcal{L}^{-1}\left\{\frac{1}{2\left(s-\frac{3}{2}\right)}\right\} = \frac{7}{2}\mathcal{L}^{-1}\left\{\frac{1}{\left(s-\frac{3}{2}\right)}\right\}$.

Since $\frac{7}{2}\mathcal{L}\left\{e^{\frac{3}{2}t}\right\} = \mathcal{L}\left\{\frac{7}{2}e^{\frac{3}{2}t}\right\} = \frac{7}{2}\left(\frac{1}{\left(s-\frac{3}{2}\right)}\right) = \frac{7}{2s-3}$ , then by using definition 1.3.1, this means that $\mathcal{L}^{-1}\left\{\frac{7}{2s-3}\right\} = \frac{7}{2}e^{\frac{3}{2}t}$.

**Example 1.3.7** Find $\mathcal{L}^{-1}\left\{\frac{1+3s}{s^2+9}\right\}$.

**Solution:** $\mathcal{L}^{-1}\left\{\frac{1+3s}{s^2+9}\right\} = \mathcal{L}^{-1}\left\{\frac{1}{s^2+9} + \frac{3s}{s^2+9}\right\}$. Now, by using section 3 in table 1.1.1, we need to make it look like $\frac{c}{s^2+c^2}$ and $\frac{s}{s^2+c^2}$ because $\mathcal{L}\{\sin ct\} = \frac{c}{s^2+c^2}$ and $\mathcal{L}\{\cos ct\} = \frac{s}{s^2+c^2}$ .

Therefore, we do the following:
$\mathcal{L}^{-1}\left\{\frac{1}{s^2+9} + \frac{3s}{s^2+9}\right\} = \frac{1}{3}\mathcal{L}^{-1}\left\{\frac{1(3)}{s^2+9}\right\} + 3\mathcal{L}^{-1}\left\{\frac{s}{s^2+9}\right\}$.

Since $\frac{1}{3}\mathcal{L}\{\sin(3t)\} + 3\mathcal{L}\{\cos(3t)\} = \left(\frac{1}{s^2+9} + \frac{3s}{s^2+9}\right)$, then by using definition 1.3.1, this means that $\mathcal{L}^{-1}\left\{\frac{1+3s}{s^2+9}\right\} = \mathcal{L}^{-1}\left\{\frac{1}{s^2+9}\right\} + \mathcal{L}^{-1}\left\{\frac{3s}{s^2+9}\right\} = \frac{1}{3}\sin(3t) + 3\cos(3t)$.

**Example 1.3.8** Find $\mathcal{L}^{-1}\left\{\frac{6}{s^4}\right\}$.

**Solution:** Using definition 1.3.1 and the right side of section 2 in table 1.1.1, we find the following:

$\mathcal{L}\{t^m\} = \frac{m!}{s^{m+1}}$ where $m$ is a positive integer.

$\frac{6}{s^4}$ is an equivalent to $\frac{3!}{s^{3+1}}$

Since $\mathcal{L}\{t^3\} = \frac{3!}{s^{3+1}} = \frac{6}{s^4}$, then this means that $\mathcal{L}^{-1}\left\{\frac{6}{s^4}\right\} = t^3$.

# 1.4 Initial Value Problems

In this section, we will introduce the main theorem of differential equations known as Initial Value Problems (IVP), and we will use it with what we have learned from sections 1.2 and 1.3 to find the largest interval on the x-axis.

**Definition 1.4.1** Given $a_n(x)y^{(n)} + a_{n-1}(x)y^{(n-1)} + \cdots + a_0(x)y(x) = K(x)$. Assume that $a_n(m) \neq 0$ for every $m \in I$ where $I$ is some interval, and $a_n(x), a_{n-1}(x), \ldots, a_0(x), K(x)$ are continuous on $I$. Suppose that $y(w), y'(w), \ldots, y^{(n-1)}(w)$ for some $w \in I$. Then, the solution to the differential equations is unique which means that there exists exactly one $y(x)$ in terms of $x$, and this type of mathematical problems is called Initial Value Problems (IVP).

**Example 1.4.1** Find the largest interval on the $x - axis$ so that $\frac{x-3}{x+2}y^{(3)}(x) + 2y^{(2)}(x) + \sqrt{x+1}\,y'(x) = 5x + 7$ has a solution. Given $y'(5) = 10, y(5) = 2, y^{(2)}(5) = -5$.

**Solution:** Finding largest interval on the $x - axis$ means that we need to find the domain for the solution of the above differential equation in other words we need to find for what values of $x$ the solution of the above differential equation holds. Therefore, we do the following:

$\boxed{\frac{x-3}{x+2}}y^{(3)}(x) + \boxed{2}y^{(2)}(x) + \boxed{\sqrt{x+1}}\,y'(x) = \boxed{5x + 7}$

Using definition 1.4.1, we also suppose the following:

$a_3(x) = \dfrac{x-3}{x+2}$

$a_2(x) = 2$

$a_1(x) = \sqrt{x+1}$

$K(x) = 5x + 7$

Now, we need to determine the interval of each coefficient above as follows:

$a_3(x) = \frac{x-3}{x+2}$ has a solution which is continuous everywhere ($\Re$) except $x = -2$ and $x = 3$.

$a_2(x) = 2$ has a solution which is continuous everywhere ($\Re$).

$a_1(x) = \sqrt{x+1}$ has a solution which is continuous on the interval $[-1, \infty)$.

$K(x) = 5x + 7$ has a solution which is continuous everywhere ($\Re$).

Thus, the largest interval on the $x - axis$ is $(3, \infty)$.

**Example 1.4.2** Find the largest interval on the $x - axis$ so that $(x^2 + 2x - 3)y^{(2)}(x) + \frac{1}{x+3}y(x) = 10$ has a solution. Given $y'(2) = 10, y(2) = 4$.

**Solution:** Finding largest interval on the $x - axis$ means that we need to find the domain for the solution of the above differential equation in other words we need to find for what values of $x$ the solution of the above differential equation holds. Therefore, we do the following:

$$\boxed{(x^2 + 2x - 3)}y^{(2)}(x) + \boxed{\frac{1}{x+3}}y(x) = \boxed{10}$$

Using definition 1.4.1, we also suppose the following:

$a_2(x) = (x^2 + 2x - 3)$

$a_1(x) = 0$

$a_0(x) = \dfrac{1}{x+3}$

$K(x) = 10$

Now, we need to determine the interval of each coefficient above as follows:

$a_2(x) = (x^2 + 2x - 3) = (x-1)(x+3)$ has a solution which is continuous everywhere ($\Re$) except $x = 1$ and $x = -3$.

$a_0(x) = \frac{1}{x+3}$ has a solution which is continuous everywhere ($\Re$) except $x = -3$.

$K(x) = 10$ has a solution which is continuous everywhere ($\Re$).

Thus, the largest interval on the $x - axis$ is $(1, \infty)$.

**Example 1.4.3** Solve the following Initial Value Problem (IVP): $y'(x) + 3y(x) = 0$. Given $y(0) = 3$.

**Solution:** $y'(x) + 3y(x) = 0$ is a linear differential equation of order 1. First, we need to find the domain for the solution of the above differential equation in other words we need to find for what values of $x$ the solution of the above differential equation holds. Therefore, we do the following:

$\boxed{1}y'(x) + \boxed{3}y(x) = \boxed{0}$

Using definition 1.4.1, we also suppose the following:

$a_1(x) = 1$

$a_0(x) = 3$

$K(x) = 0$

The domain of solution is $(-\infty, \infty)$.

Now, to find the solution of the above differential equation, we need to take the laplace transform for both sides as follows

$\mathcal{L}\{y'(x)\} + \mathcal{L}\{3y(x)\} = \mathcal{L}\{0\}$

$\mathcal{L}\{y'(x)\} + 3\mathcal{L}\{y(x)\} = 0$ because $(\mathcal{L}\{0\} = 0)$.

$\big(sY(s) - y(0)\big) + 3Y(s) = 0$    from result 1.2.2.

$sY(s) - y(0) + 3Y(s) = 0$

$Y(s)(s + 3) = y(0)$

We substitute $y(0) = 3$ because it is given in the question itself.

$Y(s)(s + 3) = 3$

$Y(s) = \frac{3}{(s+3)}$

To find a solution, we need to find the inverse laplace transform as follows:

$\mathcal{L}^{-1}\{Y(s)\} = \mathcal{L}^{-1}\left\{\frac{3}{(s+3)}\right\} = 3\mathcal{L}^{-1}\left\{\frac{1}{(s+3)}\right\}$ and we use table 1.1.1 section 4.

$y(x) = 3e^{-3x}$ (It is written in terms of $x$ instead of $t$ because we need it in terms of $x$).

Then, we will find $y'(x)$ by finding the derivative of what we got above ($y(x) = 3e^{-3x}$) as follows:

$y'(x) = (3e^{-3x})' = (3)(-3)e^{-3x} = -9e^{-3x}$.

Finally, to check our solution if it is right, we substitute what we got from $y(x) = 3e^{-3x}$ and $y'(x) = -9e^{-3x}$ in

$y'(x) + 3y(x) = 0$ as follows:

$-9e^{-3x} + 3(3e^{-3x}) = -9e^{-3x} + (9e^{-3x}) = 0$

Thus, our solution is correct which is $y(x) = 3e^{-3x}$ and $y'(x) = -9e^{-3x}$.

**Example 1.4.4** Solve the following Initial Value Problem (IVP): $y^{(2)}(x) + 3y(x) = 0$. Given $y(0) = 0$, and $y'(0) = 1$.

**Solution:** $y^{(2)}(x) + 3y(x) = 0$ is a linear differential equation of order 2. First, we need to find the domain for the solution of the above differential equation in other words we need to find for what values of $x$ the solution of the above differential equation holds. Therefore, we do the following:

$\boxed{1}y^{(2)}(x) + \boxed{3}y(x) = \boxed{0}$

Using definition 1.4.1, we also suppose the following:

$a_2(x) = 1$

$a_1(x) = 0$

$a_0(x) = 3$

$K(x) = 0$

The domain of solution is $(-\infty, \infty)$.

Now, to find the solution of the above differential equation, we need to take the laplace transform for both sides as follows

$\mathcal{L}\{y^{(2)}(x)\} + \mathcal{L}\{3y(x)\} = \mathcal{L}\{0\}$

$\mathcal{L}\{y^{(2)}(x)\} + 3\mathcal{L}\{y(x)\} = 0$ because $(\mathcal{L}\{0\} = 0)$.

$(s^2 Y(s) - sy(0) - y'(0)) + 3Y(s) = 0$  from result 1.2.2.

$s^2 Y(s) - sy(0) - y'(0) + 3Y(s) = 0$

$Y(s)(s^2 + 3) = sy(0) + y'(0)$

We substitute $y(0) = 0$, and $y'(0) = 1$ because it is given in the question itself.

$Y(s)(s^2 + 3) = (s)(0) + 1$

$Y(s)(s^2 + 3) = 0 + 1$

$Y(s)(s^2 + 3) = 1$

$Y(s) = \dfrac{1}{(s^2 + 3)}$

To find a solution, we need to find the inverse laplace transform as follows:

$\mathcal{L}^{-1}\{Y(s)\} = \mathcal{L}^{-1}\left\{\dfrac{1}{(s^2+3)}\right\} = \dfrac{1}{\sqrt{3}}\mathcal{L}^{-1}\left\{\dfrac{1}{\left(s^2+(\sqrt{3})^2\right)}\right\}$   and   we   use

table 1.1.1 at the left side of section 3.

$y(x) = \frac{1}{\sqrt{3}}\sin(\sqrt{3}\,x)$ (It is written in terms of $x$ instead of $t$ because we need it in terms of $x$).

Then, we will find $y'(x)$ by finding the derivative of what we got above $\left(y(x) = \frac{1}{\sqrt{3}}\sin(\sqrt{3}\,x)\right)$ as follows:

$$y'(x) = \left(\frac{1}{\sqrt{3}}\sin(\sqrt{3}\,x)\right)' = \frac{1}{\sqrt{3}}\left(\sqrt{3}\right)\cos(\sqrt{3}\,x) = \cos(\sqrt{3}\,x)$$

Now, we will find $y^{(2)}(x)$ by finding the derivative of what we got above as follows:

$$y^{(2)}(x) = \left(\cos(\sqrt{3}\,x)\right)' = -\sqrt{3}\sin(\sqrt{3}x)$$

Finally, to check our solution if it is right, we substitute what we got from $y(x)$ and $y^{(2)}(x)$ in $y^{(2)}(x) + 3y(x) = 0$ as follows:

$$-\sqrt{3}\sin(\sqrt{3}x) + \sqrt{3}\,sin(\sqrt{3}\,x) = 0$$

Thus, our solution is correct which is

$y(x) = \frac{1}{\sqrt{3}}\sin(\sqrt{3}\,x)$ and $y^{(2)}(x) = -\sqrt{3}\sin(\sqrt{3}x)$.

# 1.5 Properties of Laplace Transforms

In this section, we discuss several properties of laplace transforms such as shifting, unit step function, periodic function, and convolution.

We start with some examples of shifting property.

**Example 1.5.1** Find $\mathcal{L}\{e^{3x}x^3\}$.

**Solution:** By using shifting property at the left side of section 5 in table 1.1.1, we obtain:

$\mathcal{L}\{e^{bx}f(x)\} = F(s) \mid s \to s - b$

Let $b = 3$, and $f(x) = x^3$.

$F(s) = \mathcal{L}\{x^3\}$.

Hence, $\mathcal{L}\{e^{3x}x^3\} = \mathcal{L}\{x^3\} \mid s \to s - 3 = \frac{3!}{(s)^{3+1}} \mid s \to s - 3 =$

$\frac{3!}{(s)^4} \mid s \to s - 3$.

Now, we need to substitute $s$ with $s - 3$ as follows:

$$\frac{3!}{(s-3)^4} = \frac{6}{(s-3)^4}$$

Thus, $\mathcal{L}\{e^{3x}x^3\} = \frac{3!}{(s-3)^4} = \frac{6}{(s-3)^4}$.

**Example 1.5.2** Find $\mathcal{L}\{e^{-2x}\sin(4x)\}$.

**Solution:** By using shifting property at the left side of section 5 in table 1.1.1, we obtain:

$\mathcal{L}\{e^{bx}f(x)\} = F(s) \mid s \to s - b$

Let $b = -2$, and $f(x) = \sin(4x)$.

$F(s) = \mathcal{L}\{\sin(4x)\}$.

Hence, $\mathcal{L}\{e^{-2x}\sin(4x)\} = \mathcal{L}\{\sin(4x)\} \mid s \to s + 2 =$

$$\frac{4}{s^2 + 16} \mid s \to s + 2$$

Now, we need to substitute $s$ with $s + 2$ as follows:

$$\frac{4}{(s+2)^2 + 16}$$

Thus, $\mathcal{L}\{e^{-2x}\sin(4x)\} = \frac{4}{(s+2)^2+16}$.

**Example 1.5.3** Find $\mathcal{L}^{-1}\left\{\frac{s}{(s-2)^2+4}\right\}$.

**Solution:** Since we have a shift such as $s - 2$, we need to do the following:

$$\mathcal{L}^{-1}\left\{\frac{s}{(s-2)^2+4}\right\} = \mathcal{L}^{-1}\left\{\frac{s-2+2}{(s-2)^2+4}\right\}$$

$$= \mathcal{L}^{-1}\left\{\frac{s-2}{(s-2)^2+4} + \frac{2}{(s-2)^2+4}\right\}$$

$$\mathcal{L}^{-1}\left\{\frac{s}{(s-2)^2+4}\right\} = \mathcal{L}^{-1}\left\{\frac{s-2}{(s-2)^2+4}\right\} + \mathcal{L}^{-1}\left\{\frac{2}{(s-2)^2+4}\right\}$$

By using shifting property at the right side of section 5 in table 1.1.1, we obtain:

$\mathcal{L}^{-1}\{F(s) \mid s \to s - b\} = e^{bx} f(x)$

Let $b = 2$, $F_1(s) = \mathcal{L}^{-1}\left\{\frac{s-2}{(s-2)^2+4}\right\}$, and $F_2(s) = \mathcal{L}^{-1}\left\{\frac{2}{(s-2)^2+4}\right\}$.

$\mathcal{L}^{-1}\left\{\frac{s}{(s-2)^2+4}\right\} = \mathcal{L}^{-1}\{F_1(s) \mid s \to s - 2\} + \mathcal{L}^{-1}\{F_2(s) \mid s \to s - 2\}$

Thus, $\mathcal{L}^{-1}\left\{\frac{s}{(s-2)^2+4}\right\} = e^{2x}\cos(2x) + e^{2x}\sin(2x)$.

**Example 1.5.4** Find $\mathcal{L}^{-1}\left\{\frac{s}{(s+2)^3}\right\}$.

**Solution:** Since we have a shift such as $s + 2$, we need to do the following:

$$\mathcal{L}^{-1}\left\{\frac{s}{(s+2)^3}\right\} = \mathcal{L}^{-1}\left\{\frac{s+2-2}{(s+2)^3}\right\} = \mathcal{L}^{-1}\left\{\frac{s+2}{(s+2)^3} - \frac{2}{(s+2)^3}\right\}$$

$$\mathcal{L}^{-1}\left\{\frac{s}{(s+2)^3}\right\} = \mathcal{L}^{-1}\left\{\frac{s+2}{(s+2)^3}\right\} + \mathcal{L}^{-1}\left\{\frac{-2}{(s+2)^3}\right\}$$

$$\mathcal{L}^{-1}\left\{\frac{s}{(s+2)^3}\right\} = \mathcal{L}^{-1}\left\{\frac{1}{(s+2)^2}\right\} + \mathcal{L}^{-1}\left\{\frac{-2}{(s+2)^3}\right\}$$

By using shifting property at the right side of section 5 in table 1.1.1, we obtain:

$\mathcal{L}^{-1}\{F(s) \mid s \to s - b\} = e^{bx} f(x)$

Let $b = -2$, $F_1(s) = \mathcal{L}^{-1}\left\{\frac{1}{(s+2)^2}\right\}$, and $F_2(s) = \mathcal{L}^{-1}\left\{\frac{-2}{(s+2)^3}\right\}$.

$\mathcal{L}^{-1}\left\{\frac{s}{(s+2)^3}\right\} = \mathcal{L}^{-1}\{F_1(s) \,|\, s \to s+2\} + \mathcal{L}^{-1}\{F_2(s) \,|\, s \to s+2\}$

Thus, $\mathcal{L}^{-1}\left\{\frac{s}{(s+2)^3}\right\} = e^{-2x}x - e^{-2x}x^2$.

**Example 1.5.5** Find $\mathcal{L}^{-1}\left\{\frac{1}{s^2-4}\right\}$.

**Solution:** $\mathcal{L}^{-1}\left\{\frac{1}{s^2-4}\right\} = \mathcal{L}^{-1}\left\{\frac{1}{(s-2)(s+2)}\right\}$.

Since the numerator has a polynomial of degree 0 ($x^0 = 1$), and the denominator a polynomial of degree 2, then this means the degree of numerator is less than the degree of denominator. Thus, in this case, we need to use the partial fraction as follows:

$$\frac{1}{(s-2)(s+2)} = \frac{a}{(s-2)} + \frac{b}{(s+2)}$$

It is easier to use a method known as cover method than using the traditional method that takes long time to finish it. In the cover method, we cover the original, say $(s-2)$, and substitute $s = 2$ in $\frac{1}{(s-2)(s+2)}$ to find the value of $a$. Then, we cover the original, say $(s+2)$, and substitute $s = -2$ in $\frac{1}{(s-2)(s+2)}$ to find the value of $b$.

Thus, $a = \frac{1}{4}$ and $b = -\frac{1}{4}$. This implies that

$$\frac{1}{(s-2)(s+2)} = \frac{\frac{1}{4}}{(s-2)} + \frac{-\frac{1}{4}}{(s+2)}$$

Now, we need to do the following:

$$\mathcal{L}^{-1}\left\{\frac{1}{s^2-4}\right\} = \mathcal{L}^{-1}\left\{\frac{\frac{1}{4}}{(s-2)} + \frac{-\frac{1}{4}}{(s+2)}\right\}$$

$$\mathcal{L}^{-1}\left\{\frac{1}{s^2-4}\right\} = \mathcal{L}^{-1}\left\{\frac{\frac{1}{4}}{(s-2)}\right\} + \mathcal{L}^{-1}\left\{\frac{-\frac{1}{4}}{(s+2)}\right\}$$

$$\mathcal{L}^{-1}\left\{\frac{1}{s^2-4}\right\} = \frac{1}{4}\mathcal{L}^{-1}\left\{\frac{1}{(s-2)}\right\} - \frac{1}{4}\mathcal{L}^{-1}\left\{\frac{1}{(s+2)}\right\}$$

Thus, $\mathcal{L}^{-1}\left\{\frac{1}{s^2-4}\right\} = \frac{1}{4}e^{2x} - \frac{1}{4}e^{-2x}$.

Now, we will introduce a new property from table 1.1.1 in the following two examples.

**Example 1.5.6** Find $\mathcal{L}\{xe^x\}$.

**Solution:** By using the left side of section 9 in table 1.1.1, we obtain:

$$\mathcal{L}\{t^m f(t)\}(s) = (-1)^m \frac{d^m F(s)}{ds^m} = (-1)^m F^{(m)}(s)$$

where $\mathcal{L}\{f(x)\} = F(s)$, and $f(x) = e^x$.

Hence, $\mathcal{L}\{xe^x\} = (-1)^1 F^{(1)}(s) = -F^{(1)}(s)$

Now, we need to find $F^{(1)}(s)$ as follows: This means that we first need to find the laplace transform of $f(x)$, and then we need to find the first derivative of the result from the laplace transform.

$$F^{(1)}(s) = F'(s) = (\mathcal{L}\{f(x)\})' = (\mathcal{L}\{e^x\})' = \left(\frac{1}{s-1}\right)' = -\frac{1}{(s-1)^2}$$

Thus, $\mathcal{L}\{xe^x\} = (-1)^1 F^{(1)}(s) = -F^{(1)}(s) = -\left(-\frac{1}{(s-1)^2}\right) = \frac{1}{(s-1)^2}$.

**Example 1.5.7** Find $\mathcal{L}\{x^2\sin(x)\}$.

**Solution:** By using the left side of section 9 in table 1.1.1, we obtain:

$\mathcal{L}\{t^m f(t)\}(s) = (-1)^m \dfrac{d^m F(s)}{ds^m} = (-1)^m F^{(m)}(s)$

where $\mathcal{L}\{f(x)\} = F(s)$, and $f(x) = \sin(x)$.

Hence, $\mathcal{L}\{x^2 \sin(x)\} = (-1)^2 F^{(2)}(s) = F^{(2)}(s)$

Now, we need to find $F^{(2)}(s)$ as follows: This means that we first need to find the laplace transform of $f(x)$, and then we need to find the second derivative of the result from the laplace transform.

$F^{(2)}(s) = F''(s) = (\mathcal{L}\{f(x)\})'' = (\mathcal{L}\{\sin(x)\})'' = \left(\dfrac{1}{s^2+1}\right)'' =$

$\left(\dfrac{-2s}{(s^2+1)^2}\right)' = \dfrac{-2(s^2+1)^2 + 8s^2(s^2+1)}{(s^2+1)^4}$

Thus, $\mathcal{L}\{x^2 \sin(x)\} = \dfrac{-2(s^2+1)^2 + 8s^2(s^2+1)}{(s^2+1)^4}$ .

**Example 1.5.8** Solve the following Initial Value Problem (IVP): $y^{(2)}(x) + 5y^{(1)}(x) + 6y(x) = 1$. Given $y(0) = y'(0) = 0$.

**Solution:** $y^{(2)}(x) + 5y^{(1)}(x) + 6y(x) = 1$ is a linear differential equation of order 2. First, we need to find the domain for the solution of the above differential equation in other words we need to find for what values of $x$ the solution of the above differential equation holds. Therefore, we do the following:
$\boxed{1}y^{(2)}(x) + \boxed{5}y^{(1)}(x) + \boxed{6}y(x) = \boxed{1}$
Using definition 1.4.1, we also suppose the following:
$a_2(x) = 1$

$a_1(x) = 5$

$a_0(x) = 6$

$K(x) = 1$

The domain of solution is $(-\infty, \infty)$.

Now, to find the solution of the above differential equation, we need to take the laplace transform for both sides as follows

$\mathcal{L}\{y^{(2)}(x)\} + \mathcal{L}\{5y^{(1)}(x)\} + \mathcal{L}\{6y(x)\} = \mathcal{L}\{1\}$

$\mathcal{L}\{y^{(2)}(x)\} + 5\mathcal{L}\{y^{(1)}(x)\} + 6\mathcal{L}\{y(x)\} = \dfrac{1}{s}$

because $\left(\mathcal{L}\{1\} = \dfrac{1}{s}\right)$.

$(s^2 Y(s) - sy(0) - y'(0)) + 5(sY(s) - y(0)) + 6Y(s) = \dfrac{1}{s}$

from result 1.2.2. We substitute $y(0) = 0$, and $y'(0) = 0$ because it is given in the question itself.

$(s^2 Y(s) - 0 - 0) + 5(sY(s) - 0) + 6Y(s) = \dfrac{1}{s}$

$s^2 Y(s) + 5sY(s) + 6Y(s) = \dfrac{1}{s}$

$Y(s)(s^2 + 5s + 6) = \dfrac{1}{s}$

$Y(s) = \dfrac{1}{s(s^2 + 5s + 6)} = \dfrac{1}{s(s + 3)(s + 2)}$

To find a solution, we need to find the inverse laplace transform as follows:

$y(x) = \mathcal{L}^{-1}\{Y(s)\} = \mathcal{L}^{-1}\left\{\dfrac{1}{s(s+3)(s+2)}\right\}$

Since the numerator has a polynomial of degree 0 ($x^0 = 1$), and the denominator a polynomial of degree 3, then this means the degree of numerator is less than the degree of denominator. Thus, in this case, we need to use the partial fraction as follows:

$$\frac{1}{s(s+3)(s+2)} = \frac{a}{s} + \frac{b}{(s+3)} + \frac{c}{(s+2)}$$

Now, we use the cover method. In the cover method, we cover the original, say $s$, and substitute $s = 0$ in $\frac{1}{s(s+3)(s+2)}$ to find the value of $a$. We cover the original, say $(s+3)$, and substitute $s = -3$ in $\frac{1}{s(s+3)(s+2)}$ to find the value of $b$. Then, we cover the original, say $(s+2)$, and substitute $s = -2$ in $\frac{1}{s(s+3)(s+2)}$ to find the value of $c$. Thus, $= \frac{1}{6}$, $b = \frac{1}{3}$ and $c = -\frac{1}{2}$. This implies that

$$\frac{1}{s(s+3)(s+2)} = \frac{\frac{1}{6}}{s} + \frac{\frac{1}{3}}{(s+3)} + \frac{-\frac{1}{2}}{(s+2)}$$

Now, we need to do the following:

$$\mathcal{L}^{-1}\left\{\frac{1}{s(s+3)(s+2)}\right\} = \mathcal{L}^{-1}\left\{\frac{\frac{1}{6}}{s} + \frac{\frac{1}{3}}{(s+3)} + \frac{-\frac{1}{2}}{(s+2)}\right\}$$

$$\mathcal{L}^{-1}\left\{\frac{1}{s(s+3)(s+2)}\right\} = \frac{1}{6}\mathcal{L}^{-1}\left\{\frac{1}{s}\right\} + \frac{1}{3}\mathcal{L}^{-1}\left\{\frac{1}{(s+3)}\right\}$$

$$-\frac{1}{2}\mathcal{L}^{-1}\left\{\frac{1}{(s+2)}\right\}$$

Thus, $y(x) = \mathcal{L}^{-1}\left\{\frac{1}{s(s+3)(s+2)}\right\} = \frac{1}{6} + \frac{1}{3}e^{-3x} - \frac{1}{2}e^{-2x}$.

**Definition 1.5.1** Given $a > 0$. Unit Step Function is defined as follows: $U(x - a) = \begin{cases} 0 & if \ 0 \le x < a \\ 1 & if \ a \le x < \infty \end{cases}$

**Result 1.5.1** Given $a \ge 0$. Then, we obtain:

  a) $U(x - 0) = U(x) = 1$ for every $0 \le x \le \infty$.

  b) $U(x - \infty) = 0$ for every $0 \le x \le \infty$.

**Example 1.5.9** Find $U(x - 7)|x = 8$.

**Solution:** Since $x = 8$ is between $a = 7$ and $\infty$, then by using definition 1.5.1, we obtain:

$U(x - 7)|x = 8 = U(8 - 7) = 1$.

**Example 1.5.10** Find $U(x - 3)|x = 1$.

**Solution:** Since $x = 1$ is between $0$ and $= 3$ , then by using definition 1.5.1, we obtain:

$U(x - 3)|x = 1 = U(1 - 7) = 0$.

**Example 1.5.11** Find $\mathcal{L}\{U(x - 3)\}$.

**Solution:** By using the left side of section 7 in table 1.1.1, we obtain:

$\mathcal{L}\{U(t - b)\} = \dfrac{e^{-bs}}{s}$

Thus, $\mathcal{L}\{U(x - 3)\} = \dfrac{e^{-3s}}{s}$.

**Example 1.5.12** Find $\mathcal{L}^{-1}\left\{\dfrac{e^{-2s}}{s}\right\}$.

**Solution:** By using the right side of section 7 in table 1.1.1, we obtain:

$\mathcal{L}^{-1}\left\{\dfrac{e^{-bs}}{s}\right\} = U(t - b)$

Thus, $\mathcal{L}^{-1}\left\{\dfrac{e^{-2s}}{s}\right\} = U(x - 2) = \begin{cases} 0 & if\ 0 \leq x < 2 \\ 1 & if\ 2 \leq x < \infty \end{cases}$

**Example 1.5.13** Given $f(x) = \begin{cases} 3 & if\ 1 \leq x < 4 \\ -x & if\ 4 \leq x < 10 \\ (x + 1) & if\ 10 \leq x < \infty \end{cases}$

Rewrite $f(x)$ in terms of *Unit Step Function*.

**Solution:** To re-write $f(x)$ in terms of

*Unit Step Function*, we do the following:

$f(x) = 3\big(U(x-1) - U(x-4)\big)$

$\qquad\qquad + (-x)\big(U(x-4) - U(x-10)\big)$

$\qquad\qquad + (x+1)(U(x-10) - 0).$

Now, we need to check our unit step functions as follows: We choose $x = 9$.

$f(9) = 3\big(U(9-1) - U(9-4)\big)$

$\qquad\qquad + (-9)\big(U(9-4) - U(9-10)\big)$

$\qquad\qquad + (9+1)(U(9-10) - 0)$

$f(9) = 3(1-1) + (-9)(1-0) + (10)(0-0)$

$f(9) = 3(0) + (-9)(1) + (10)(0)$

$f(9) = 0 + (-9)(1) + 0 = -9.$

Thus, our unit step functions are correct.

**Example 1.5.14** Find $\mathcal{L}\{xU(x-2)\}$.

**Solution:** By using the upper side of section 6 in table

1.1.1, we obtain:

$\mathcal{L}\{h(x)U(x-b)\} = e^{-bs}\mathcal{L}\{h(x+b)\}$

where $b = 2$, and $h(x) = x$.

Hence, $\mathcal{L}\{xU(x-2)\} = e^{-2s}\mathcal{L}\{h(x+2)\} = e^{-2s}\mathcal{L}\{x+2\} =$

$e^{-2s}\left(\frac{1}{s^2} + \frac{2}{s}\right).$

**Definition 1.5.2** Convolution, denoted by $*$, is defined

as follows:

$$f(x) * h(x) = \int\limits_0^x f(\psi)h(x-\psi)d\psi$$

where $f(x)$ and $h(x)$ are functions. (Note: do not confuse between multiplication and convolution).

**Result 1.5.2** $\mathcal{L}\{f(x) * h(x)\} = \mathcal{L}\{h(x) * f(x)\} = F(s) \cdot H(s)$ where $f(x)$ and $h(x)$ are functions. (The proof of this result left as an exercise 16 in section 1.7).

**Result** 1.5.3 $\mathcal{L}\{f(x) \cdot h(x)\} \neq \mathcal{L}\{h(x)\} \cdot \mathcal{L}\{f(x)\}$.

**Example 1.5.15** Use definition 1.5.2 to find $\mathcal{L}\{\int_0^x sin(\psi)d\psi\}$.

**Solution:** By using definition 1.5.2  and section 10 in table 1.1.1, we obtain:

$$\mathcal{L}\left\{\int_0^x sin(\psi)d\psi\right\} = \mathcal{L}\{1 * sin(x)\} = \mathcal{L}\{1\} \cdot \mathcal{L}\{sin(x)\} = \frac{1}{s} \cdot \frac{1}{(s^2 + 1)}$$

$$= \frac{1}{s(s^2 + 1)}$$

Thus, $\mathcal{L}\{\int_0^x sin(\psi)d\psi\} = \frac{1}{s(s^2+1)}$.

**Definition 1.5.3** $f(x)$ is a periodic function on $[0, \infty)$ if $f(x)$ has a period $P$ such that $f(b) = f(b - P)$ for every $b \geq P$.

**Example 1.5.16** Given $f(x)$ is periodic on $[0, \infty)$ such that the first period of $f(x)$ is given by the following piece-wise continuous function:
$$\begin{cases} 3 & if \ 0 \leq x < 2 \\ -2 & if \ 2 \leq x < 8 \end{cases}$$
  a)  Find the 8th period of this function.
  b)  Suppose $P = 8$. Find $f(10)$.
  c)  Suppose $P = 8$. Find $f(30)$.

**Solution: Part a:** By using section 13 in table 1.1.1, we obtain:

$$\mathcal{L}\{f(x)\} = \frac{1}{1 - e^{-Ps}} \int\limits_0^P e^{-st} f(x) dx$$

Since we need find the 8$^{th}$ period, then this means that $P = 8$, and we can apply what we got above as follows:

$$\mathcal{L}\{f(x)\} = \frac{1}{1 - e^{-8s}} \int\limits_0^8 e^{-st} f(x) dx = \frac{1}{1 - e^{-8s}} \int\limits_0^8 f(x) e^{-st} dx$$

Using the given first period function, we obtain:

$$\mathcal{L}\{f(x)\} = \frac{1}{1 - e^{-8s}} \int\limits_0^8 f(x) e^{-st} dx$$

$$= \frac{1}{1 - e^{-8s}} \left[ 3 \int\limits_0^2 e^{-st} dx - 2 \int\limits_2^8 e^{-st} dx \right]$$

$$= \frac{1}{1 - e^{-8s}} \left[ -\frac{3}{s} e^{-st} \Big|_{x=0}^{x=2} + \frac{2}{s} e^{-st} \Big|_{x=2}^{x=8} \right]$$

$$= \frac{1}{1 - e^{-8s}} \left[ -\frac{3}{s}(e^{-2t} - 1) + \frac{2}{s}(e^{-8t} - e^{-2t}) \right]$$

**Part b:** By using definition 1.5.3, we obtain:

$f(10) = f(10 - 8) = f(2) = -2$ from the given first period function.

**Part c:** By using definition 1.5.3, we obtain:

$f(30) = f(30 - 8) = f(22)$.

$f(22) = f(22 - 8) = f(14)$.

$f(14) = f(14 - 8) = f(6) = -2$ from the given first period function.

**Definition 1.5.4** Suppose that $i > 0$ is fixed, and $\delta < j < i$ is chosen arbitrary. Then, we obtain:

$$\delta_j(t-i) = \begin{cases} 0 & if\ 0 \leq t < (i-j) \\ \dfrac{1}{2j} & if\ (i-j) \leq t < (i+j) \\ 0 & t \geq (i+j) \end{cases}$$

$\delta_j(t-i)$ is called *Dirac Delta Function*.

**Result 1.5.4** $\delta(t-i) = \lim\limits_{i \to 0^+} \delta_j(t-i)$.

**Example 1.5.17** Find $\mathcal{L}\{\delta(t+12)\}$.

**Solution:** By using the right side of section 12 in table 1.1.1, we obtain: $\mathcal{L}\{\delta(t-b)\} = e^{-bs}$

Thus, $\mathcal{L}\{\delta(t+12)\} = e^{12s}$.

# 1.6 Systems of Linear Equations

Most of the materials of this section are taken from section 1.8 in my published book titled *A First Course in Linear Algebra: Study Guide for the Undergraduate Linear Algebra Course, First Edition[1],* because it is very important to give a review from linear algebra about Cramer's rule, and how some concepts of linear algebra can be used to solve some problems in differential equations. In this section, we discuss how to use what we have learned from previous sections such as initial value problems (IVP), and how to use Cramer's rule to solve systems of linear equations.

**Definition 1.6.1** Given $n \times n$ system of linear equations.

Let WX = A be the matrix form of the given system:

$$W \begin{bmatrix} x_1 \\ x_2 \\ x_3 \\ \vdots \\ x_n \end{bmatrix} = \begin{bmatrix} a_1 \\ a_2 \\ a_3 \\ \vdots \\ a_n \end{bmatrix}$$

The system has a unique solution if and only if

$\det(W) \neq 0$. Cramer's Rule tells us how to

find $x_1, x_2, \ldots, x_n$ as follows:

Let $W = \begin{bmatrix} 1 & 3 & 4 \\ 1 & 2 & 1 \\ 7 & 4 & 3 \end{bmatrix}$ Then, the solutions for the system of

linear equations are:

$$x_1 = \frac{\det \begin{bmatrix} a_1 & 3 & 4 \\ \vdots & 2 & 1 \\ a_n & 4 & 3 \end{bmatrix}}{\det(W)}$$

$$x_2 = \frac{\det \begin{bmatrix} 1 & a_1 & 4 \\ 1 & \vdots & 1 \\ 7 & a_n & 3 \end{bmatrix}}{\det(W)}$$

$$x_3 = \frac{\det \begin{bmatrix} 1 & 3 & a_1 \\ 1 & 2 & \vdots \\ 7 & 4 & a_n \end{bmatrix}}{\det(W)}$$

**Example 1.6.1** Solve the following system of linear equations using Cramer's Rule:

$$\begin{cases} 2x_1 + 7x_2 = 13 \\ -10x_1 + 3x_2 = -4 \end{cases}$$

**Solution:** First of all, we write $2 \times 2$ system in the form WX = A according to definition 1.6.1.

$$\begin{bmatrix} 2 & 7 \\ -10 & 3 \end{bmatrix} \begin{bmatrix} x_1 \\ x_2 \end{bmatrix} = \begin{bmatrix} 13 \\ -4 \end{bmatrix}$$

Since W in this form is $\begin{bmatrix} 2 & 7 \\ -10 & 3 \end{bmatrix}$, then

$\det(W) = (2 \cdot 3) - \big(7 \cdot (-10)\big) = 6 - (-70) = 76 \neq 0.$

The solutions for this system of linear equations are:

$$x_1 = \frac{\det \begin{bmatrix} 13 & 7 \\ -4 & 3 \end{bmatrix}}{\det(W)} = \frac{\det \begin{bmatrix} 13 & 7 \\ -4 & 3 \end{bmatrix}}{76} = \frac{67}{76}$$

$$x_2 = \frac{\det \begin{bmatrix} 2 & 13 \\ -10 & -4 \end{bmatrix}}{\det(W)} = \frac{\det \begin{bmatrix} 2 & 13 \\ -10 & -4 \end{bmatrix}}{76} = \frac{122}{76}$$

Thus, the solutions are $x_1 = \frac{67}{76}$ and $x_2 = \frac{122}{76}$.

**Example 1.6.2** Solve for $n(t)$ and $m(t)$:
$$\begin{cases} n''(t) - 4m(t) = 0 \dots \dots \dots \dots \dots \dots Equation\ 1 \\ n(t) + 2m'(t) = 5e^{2t} \dots \dots \dots \dots Equation\ 2 \end{cases}$$
Given that $n(0) = 1, n'(0) = 2,$ and $m(0) = 1.$

**Solution:** First, we need to take the laplace transform of both sides for each of the above two equations.

For *Equation* 1: We take the laplace transform of both sides:

$\mathcal{L}\{n''(t)\} + \mathcal{L}\{-4m(t)\} = \mathcal{L}\{0\}$

$\mathcal{L}\{n''(t)\} - 4\mathcal{L}\{m(t)\} = \mathcal{L}\{0\}$

$(s^2 N(s) - sn(0) - n'(0)) - 4M(s) = 0$

Now, we substitute what is given in this question to obtain the following:

$(s^2 N(s) - s - 2) - 4M(s) = 0$

Thus, $s^2 N(s) - 4M(s) = s + 2.$

For *Equation* 2: We take the laplace transform of both sides:

$$\mathcal{L}\{n(t)\} + \mathcal{L}\{2m'(t)\} = \mathcal{L}\{5e^{2t}\}$$

$$\mathcal{L}\{n(t)\} + 2\mathcal{L}\{m'(t)\} = 5\mathcal{L}\{e^{2t}\}$$

$$N(s) + 2\big(sM(s) - m(0)\big) = \frac{5}{s-2}$$

Now, we substitute what is given in this question to obtain the following:

$$N(s) + 2(sM(s) - 1) = \frac{5}{s-2}$$

$$N(s) + 2sM(s) = \frac{5}{s-2} + 2 = \frac{2s+1}{s-2}$$

Thus, $N(s) + 2sM(s) = \frac{2s+1}{s-2}$.

From what we got from *Equation 1* and *Equation 2*, we need to find $N(s)$ and $M(s)$ as follows:

$$\begin{cases} s^2 N(s) - 4M(s) = s + 2 \\ N(s) + 2sM(s) = \dfrac{2s+1}{s-2} \end{cases}$$

Now, we use Cramer's rule as follows:

$$N(s) = \frac{\det\begin{bmatrix} s+2 & -4 \\ \dfrac{2s+1}{s-2} & 2s \end{bmatrix}}{\det\begin{bmatrix} s^2 & -4 \\ 1 & 2s \end{bmatrix}} = \frac{2s^2 + 4s + \dfrac{8s+4}{s-2}}{2s^3 + 4}$$

$$= \frac{(s-2)(2s^2 + 4s) + 8s + 4}{(s-2)(2s^3 + 4)}$$

$$= \frac{2s^3 + 4s^2 - 4s^2 - 8s + 8s + 4}{(s-2)(2s^3 + 4)} = \frac{1}{s-2}$$

Hence, $n(t) = \mathcal{L}^{-1}\{N(s)\} = \mathcal{L}^{-1}\left\{\frac{1}{s-2}\right\} = e^{2t}$.

Further, we can use one of the given equations to find $m(t)$ as follows: $n''(t) - 4m(t) = 0$

We need to find the second derivative of $n(t)$.

$n'(t) = 2e^{2t}$.

$n''(t) = 4e^{2t}$.

Now, we can find $m(t)$ as follows:

$n''(t) - 4m(t) = 0 \rightarrow m(t) = \dfrac{n''(t)}{4} = \dfrac{4e^{2t}}{4} = e^{2t}$.

Thus, $m(t) = e^{2t}$.

# 1.7 Exercises

1. Find $\mathcal{L}^{-1}\left\{\dfrac{10}{(s-4)^4}\right\}$.

2. Find $\mathcal{L}^{-1}\left\{\dfrac{s+5}{(s-1)^2+16}\right\}$.

3. Find $\mathcal{L}^{-1}\left\{\dfrac{s+5}{(s+3)^4}\right\}$.

4. Find $\mathcal{L}^{-1}\left\{\dfrac{5}{3s-10}\right\}$.

5. Find $\mathcal{L}^{-1}\left\{\dfrac{2}{s^2-6s+13}\right\}$.

6. Find $\mathcal{L}^{-1}\left\{\dfrac{5}{s^2-7s-8}\right\}$.

7. Solve the following Initial Value Problem (IVP): $2y'(x) + 6y(x) = 0$. Given $y(0) = -4$.

8. Solve the following Initial Value Problem (IVP): $y''(x) + 4y(x) = 0$. Given $y(0) = 2$, and $y'(0) = 0$.

9. Find $\mathcal{L}^{-1}\left\{\dfrac{4}{(s-1)^2(s+3)}\right\}$.

10. Find $\mathcal{L}\{U\left(x - \frac{\pi}{2}\right)\sin(x)\}$.

11. Find $\mathcal{L}\{U(x-2)e^{3x}\}$.

12. Given $f(x) = \begin{cases} 4 & if\ 0 \le x < 3 \\ e^{2x} & if\ 3 \le x < \infty \end{cases}$

Rewrite $f(x)$ in terms of *Unit Step Function*.

13. Find $\mathcal{L}^{-1}\left\{\dfrac{se^{-4x}}{s^2+4}\right\}$.

**14.** Solve the following Initial Value Problem (IVP):

$y''(x) + 3y'(x) - 4y(x) = f(0)$. Given $f(x) =$

$\begin{cases} 1 & if\ 0 \le x < 3 \\ 0 & Otherwise \end{cases}$

$y(0) = y'(0) = 0$.

**15.** Solve the following Initial Value Problem (IVP):

$y''(x) + 7y'(x) - 8y(x) = f(x)$ where $f(x) =$

$\begin{cases} 3 & if\ 0 \le x < 5 \\ -2 & if\ 5 \le x < \infty \end{cases}$

$y(0) = y'(0) = 0$.

**16.** Prove result 1.5.2.

**17.** Solve the following Initial Value Problem (IVP):

$y^{(3)}(x) + y'(x) = U(x - 3)$. Given $y(0) = y'(0) = y''(0) = 0$.

**18.** Find $\mathcal{L}\{\int_0^x cos(2\psi)e^{(3x+2\psi)}d\psi\}$.

**19.** Find $\mathcal{L}^{-1}\left\{\frac{2s}{(s^2+4)^2}\right\}$.

**20.** Solve for $r(t)$ and $k(t)$:
$\begin{cases} r'(t) - 2k(t) = 0 \dots\dots\dots\dots\dots. Equation\ 1 \\ r''(t) - 2k'(t) = 0 \dots\dots\dots\dots. Equation\ 2 \end{cases}$
Given that $k(0) = 1, r'(0) = 2,$ and $r(0) = 0$.

**21.** Solve for $w(t)$ and $h(t)$:

$\begin{cases} w(t) - \int_0^t h(\psi)d\psi = 1 \dots\dots\dots\dots. Equation\ 1 \\ w''(t) + h'(t) = 4 \dots\dots\dots\dots. Equation\ 2 \end{cases}$

Given that $w(0) = 1, w'(0) = h(0) = 0$.

**22.** Find $f(t)$ such that $f(t) = e^{-3s} + \int_0^t f(\psi)d\psi$.

(Hint: Use laplace transform to solve this problem)

# Chapter 2

# Systems of Homogeneous Linear Differential Equations (HLDE)

In this chapter, we start introducing the homogeneous linear differential equations (HLDE) with constant coefficients. In addition, we discuss how to find the general solution of HLDE. At the end of this chapter, we introduce a new method called Undetermined Coefficient Method.

## 2.1 HLDE with Constant Coefficients

In this section, we discuss how to find the general solution of the homogeneous linear differential equations (HLDE) with constant coefficients.

To give an introduction about HLDE, it is important to start with the definition of homogeneous system.

**Definition 2.1.1** *Homogeneous System* is defined as a $m \times n$ system of linear equations that has all zero constants. (i.e. the following is an example of homogeneous system): $\begin{cases} 2a + b - c + d = 0 \\ 3a + 5b + 3c + 4d = 0 \\ -b + c - d = 0 \end{cases}$

*Definition 2.1.1 is taken from section 3.1 in my published book titled *A First Course in Linear Algebra: Study Guide for the Undergraduate Linear Algebra Course, First Edition*[1].

**Example 2.1.1** Describe the following differential equation: $y''(x) + 3y'(x) = 0$.

**Solution:** Since the above differential equation has a zero constant, then according to definition 2.1.1, it is a homogeneous differential equation. In addition, it is linear because the dependent variable $y$ and all its derivatives are to the power 1. For the order of this homogeneous differential equation, since the highest derivative is 2, then the order is 2. Thus,

$y''(x) + 3y'(x) = 0$ is a homogeneous linear differential equation of order 2.

**Example 2.1.2** Describe the following differential equation: $3y^{(3)}(x) - 2y'(x) + 7y(x) = 0$.

**Solution:** Since the above differential equation has a zero constant, then according to definition 2.1.1, it is a homogeneous differential equation. In addition, it is linear because the dependent variable $y$ and all its derivatives are to the power 1. For the order of this homogeneous differential equation, since the highest derivative is 3, then the order is 3.

Thus, $3y^{(3)}(x) - 2y'(x) + 7y(x) = 0$ is a homogeneous linear differential equation of order 3.

**Example 2.1.3** Describe the following differential equation: $3y^{(2)}(x) - 2y'(x) = 12$.

**Solution:** Since the above differential equation has a nonzero constant, then according to definition 2.1.1, it is a non-homogeneous differential equation. In addition, it is linear because the dependent variable $y$ and all its derivatives are to the power 1. For the order of this non-homogeneous differential equation, since the highest derivative is 2, then the order is 2.

Thus, $3y^{(2)}(x) - 2y'(x) = 12$ is a non-homogeneous linear differential equation of order 2.

**Example 2.1.4** Find the general solution the following Initial Value Problem (IVP): $2y'(x) + 4y(x) = 0$. Given: $y'(0) = 1$. (Hint: Use the concepts of section 1.4)

**Solution:** $2y'(x) + 4y(x) = 0$ is a homogeneous linear differential equation of order 1. First, we need to find the domain for the solution of the above differential equation in other words we need to find for what values of $x$ the solution of the above differential equation holds. Therefore, we do the following: $\boxed{2}y'(x) + \boxed{4}y(x) = \boxed{0}$

Using definition 1.4.1, we also suppose the following:

$a_1(x) = 2$

$a_0(x) = 4$

$K(x) = 0$

The domain of solution is $(-\infty, \infty)$.

Now, to find the solution of the above differential equation, we need to take the laplace transform for both sides as follows

$\mathcal{L}\{2y'(x)\} + \mathcal{L}\{4y(x)\} = \mathcal{L}\{0\}$

$2\mathcal{L}\{y'(x)\} + 4\mathcal{L}\{y(x)\} = 0$ because $(\mathcal{L}\{0\} = 0)$.

$2(sY(s) - y'(0)) + 4Y(s) = 0$　 from result 1.2.2.

$2sY(s) - 2y'(0) + 4Y(s) = 0$

$Y(s)(2s + 4) = 2y'(0)$

$Y(s)(2s + 4) = 2(1)$

$Y(s) = \dfrac{2}{2s + 4} = \dfrac{2}{2(s + 2)} = \dfrac{1}{s + 2}$

To find a solution, we need to find the inverse laplace transform as follows:

$y(x) = \mathcal{L}^{-1}\{Y(s)\} = \mathcal{L}^{-1}\left\{\dfrac{1}{s+2}\right\} = e^{-2x}$.

Thus, the general solution $y(x) = ce^{-2x}$, for some constant $c$. Here, $c = 1$.

**Result 2.1.1** Assume that $m_1 y^{(n)}(x) + m_2 y(x) = 0$ is a homogeneous linear differential equation of order $n$ with constant coefficients $m_1$ and $m_2$. Then,

a)  $m_1 y^{(n)}(x) + m_2 y(x) = 0$ must have exactly $n$ independent solutions, say $f_1(x), f_2(x), \dots, f_n(x)$.

b)  Every solution of $m_1 y^{(n)}(x) + m_2 y(x) = 0$ is of the form: $c_1 f_1(x) + c_2 f_2(x) + \cdots + c_n f_n(x)$, for some constants $c_1, c_2, \dots, c_n$.

**Result 2.1.2** Assume that $y_1 = e^{k_1 x}$, $y_2 = e^{k_2 x}$, ..., $y_n = e^{k_n x}$ are independent if and only if $k_1, k_2, ..., k_n$ are distinct real numbers.

**Example 2.1.5** Given $y''(x) - 4y(x) = 0$,

$y(0) = 4, y'(0) = 10$. Find the general solution for $y(x)$. (Hint: Use results 2.1.1 and 2.1.2)

**Solution:** $y''(x) - 4y(x) = 0$ is a homogeneous linear differential equation of order 2. In this example, we will use a different approach from example 2.1.4 (laplace transform approach) to solve it. Since $y''(x) - 4y(x) = 0$ is HLDE with constant coefficients, then we will do the following: Let $y(x) = e^{kx}$, we need to find $k$.

First of all, we will find the first and second derivatives as follows:

$$y'(x) = ke^{kx}$$
$$y''(x) = k^2 e^{kx}$$

Now, we substitute $y(x) = e^{kx}$ and $y'(x) = ke^{kx}$ in $y''(x) - 4y(x) = 0$ as follows:

$$k^2 e^{kx} - 4e^{kx} = 0$$
$$e^{kx}(k^2 - 4) = 0$$
$$e^{kx}(k - 2)(k + 2) = 0$$

Thus, $k = 2$ *or* $k = -2$. Then, we use our values to substitute $k$ in our assumption which is $y(x) = e^{kx}$:

at $k = 2$, $y_1(x) = e^{2x}$

at $k = -2$, $y_2(x) = e^{-2x}$

Hence, using result 2.1.1, the general solution for $y(x)$ is: $y_{homo}(x) = c_1 e^{2x} + c_2 e^{-2x}$, for some $c_1, c_2 \in \Re$. (Note: *homo* denotes to homogeneous). Now, we need to find the values of $c_1$ and $c_2$ as follows:

at $x = 0$, $\qquad y_{homo}(0) = c_1 e^{2(0)} + c_2 e^{-2(0)}$

$$y_{homo}(0) = c_1 e^0 + c_2 e^0$$

$$y_{homo}(0) = c_1 + c_2$$

Since $y(0) = 4$, then $c_1 + c_2 = 4 \ldots \ldots \ldots \ldots \ldots \ldots \ldots \ldots (1)$

$$y_{homo}'(x) = 2c_1 e^{2x} - 2c_2 e^{-2x}$$

at $x = 0$, $\qquad y_{homo}'(0) = 2c_1 e^{2(0)} - 2c_2 e^{-2(0)}$

$$y_{homo}'(0) = 2c_1 e^0 - 2c_2 e^0$$

$$y_{homo}'(0) = 2c_1 - 2c_2$$

Since $y'(0) = 10$, then $2c_1 - 2c_2 = 10 \ldots \ldots \ldots \ldots \ldots \ldots (2)$

From (1) and (2), $c_1 = 4.5$ and $c_2 = -0.5$.

Thus, the general solution is:

$y_{homo}(x) = 4.5 e^{2x} - 0.5 e^{-2x}$.

**Example 2.1.6** Given $y^{(3)}(x) - 5y^{(2)}(x) + 6y'(x) = 0$. Find the general solution for $y(x)$. (Hint: Use results 2.1.1 and 2.1.2, and in this example, no need to find the values of $c_1, c_2,$ and $c_3$)

**Solution:** $y^{(3)}(x) - 5y^{(2)}(x) + 6y'(x) = 0$ is a homogeneous linear differential equation of order 3. Since $y^{(3)}(x) - 5y^{(2)}(x) + 6y'(x) = 0$ is HLDE with constant coefficients, then we will do the following: Let $y(x) = e^{kx}$, we need to find $k$.

First of all, we will find the first, second, and third derivatives as follows:

$$y'(x) = ke^{kx}$$
$$y''(x) = k^2 e^{kx}$$
$$y^{(3)}(x) = k^3 e^{kx}$$

Now, we substitute $y'(x) = ke^{kx}$, $y''(x) = k^2 e^{kx}$, and $y^{(3)}(x) = k^3 e^{kx}$ in $y^{(3)}(x) - 5y^{(2)}(x) + 6y'(x) = 0$ as follows:

$$k^3 e^{kx} - 5k^2 e^{kx} + 6ke^{kx} = 0$$
$$e^{kx}(k^3 - 5k^2 + 6k) = 0$$
$$e^{kx}(k(k^2 - 5k + 6)) = 0$$
$$e^{kx}(k(k-2)(k-3)) = 0$$

Thus, $k = 0, k = 2 \; or \; k = 3$. Then, we use our values to substitute $k$ in our assumption which is $y(x) = e^{kx}$:

at $k = 0$, $y_1(x) = e^{(0)x} = 1$

at $k = 2$, $y_2(x) = e^{2x}$

at $k = 3$, $y_3(x) = e^{3x}$

Thus, using result 2.1.1, the general solution for $y(x)$ is: $y_{homo}(x) = c_1 + c_2 e^{2x} + c_3 e^{3x}$, for some $c_1, c_2, c_3 \in \Re$. (Note: *homo* denotes to homogeneous).

**Example 2.1.7** Given $y^{(5)}(x) - y^{(4)}(x) - 2y^{(3)}(x) = 0$.

Find the general solution for $y(x)$. (Hint: Use results 2.1.1 and 2.1.2, and in this example, no need to find the values of $c_1, c_2, c_3, c_4$ and $c_5$)

**Solution:** $y^{(5)}(x) - y^{(4)}(x) - 2y^{(3)}(x) = 0$ is a homogeneous linear differential equation of order 5. Since $y^{(5)}(x) - y^{(4)}(x) - 2y^{(3)}(x) = 0$ is HLDE with constant coefficients, then we will do the following: Let $y(x) = e^{kx}$, we need to find $k$.

First of all, we will find the first, second, third, fourth, and fifth derivatives as follows:

$$y'(x) = ke^{kx}$$
$$y''(x) = k^2 e^{kx}$$
$$y^{(3)}(x) = k^3 e^{kx}$$
$$y^{(4)}(x) = k^4 e^{kx}$$
$$y^{(5)}(x) = k^5 e^{kx}$$

Now, we substitute $y^{(5)}(x) = k^5 e^{kx}$, $y^{(4)}(x) = k^4 e^{kx}$, and $y^{(3)}(x) = k^3 e^{kx}$ in $y^{(5)}(x) - y^{(4)}(x) - 2y^{(3)}(x) = 0$ as follows:

$$k^5 e^{kx} - k^4 e^{kx} - 2k^3 e^{kx} = 0$$
$$e^{kx}(k^5 - k^4 - 2k^3) = 0$$
$$e^{kx}(k^3(k^2 - k - 2)) = 0$$
$$e^{kx}(k^3(k - 2)(k + 1)) = 0$$

Thus, $k = 0, k = 0, k = 0, k = 2$ *or* $k = -1$. Then, we use our values to substitute $k$ in our assumption which is $y(x) = e^{kx}$:

at $k = 0$, $y_1(x) = e^{(0)x} = 1$

at $k = 0$, $y_2(x) = e^{(0)x} = 1 \cdot x$

at $k = 0$, $y_3(x) = e^{(0)x} = 1 \cdot x^2$

because $k^3 = Span\{1, x, x^2\} = 0$ (Note: $Span\{1, x, x^2\} = 0$ means $a \cdot 1 + b \cdot x + c \cdot x^2 = 0$)

In other words $Span\{1, x, x^2\}$ is the set of all linear combinations of $1$, $x$, and $x^2$.

at $k = 2$, $y_4(x) = e^{2x}$

at $k = -1$, $y_5(x) = e^{-x}$

Thus, using result 2.1.1, the general solution for $y(x)$ is: $y_{homo}(x) = c_1 + c_2 x + c_3 x^2 + c_4 e^{2x} + c_5 e^{-x}$, for some $c_1, c_2, c_3, c_4, c_5 \in \Re$. (Note: $homo$ denotes to homogeneous).

**Example 2.1.8** Given $y_1(x) = e^{3x}, y_2(x) = e^{-3x}$, and $y_3(x) = e^x$. Are $y_1(x), y_3(x)$, and $y_3(x)$ independent?

**Solution:** We cannot write $y_3(x)$ as a linear combination of $y_1(x)$ and $y_2(x)$ as follows:

$e^x \neq (Fixed\ Number) \cdot e^{3x} + (Fixed\ Number) \cdot e^{-3x}$

Thus, $y_1(x), y_3(x)$, and $y_3(x)$ are independent.

**Example 2.1.9** Given $y_1(x) = e^{(x+3)}, y_2(x) = e^3$, and $y_3(x) = e^x$. Are $y_1(x), y_3(x)$, and $y_3(x)$ independent?

**Solution:** We can write $y_1(x)$ as a linear combination of $y_2(x)$ and $y_3(x)$ as follows:

$e^{(x+3)} = e^x \cdot e^3$. Thus, $y_1(x), y_3(x)$, and $y_3(x)$ are dependent (not independent).

# 2.2 Method of Undetermined Coefficients

In this section, we discuss how to use what we have learned from section 2.1 to combine it with what we will learn from section 2.2 in order to find the general solution using a method known as undetermined coefficients method. In this method, we will find a general solution consisting of homogeneous solution and particular solution together.

We give the following examples to introduce the undetermined coefficient method.

**Example 2.2.1** Given $y'(x) + 3y(x) = x$, $y(0) = 1$. Find the general solution for $y(x)$. (Hint: No need to find the value of $c_1$)

**Solution:** Since $y'(x) + 3y(x) = x$ does not have a constant coefficient, then we need to use the undetermined coefficients method as follows:

**Step 1:** We need to find the homogeneous solution by letting $y'(x) + 3y(x)$ equal to zero as follows:

$y'(x) + 3y(x) = 0$. Now, it is a homogeneous linear differential equation of order 1.

Since $y'(x) + 3y(x) = 0$ is a HLDE with constant coefficients, then we will do the following:

Let $y(x) = e^{kx}$, we need to find $k$.

First of all, we will find the first, second, third, fourth, and fifth derivatives as follows:

$$y'(x) = ke^{kx}$$

Now, we substitute $y'(x) = ke^{kx}$ in $y'(x) + 3y(x) = 0$ as follows:

$$ke^{kx} + 3e^{kx} = 0$$
$$e^{kx}(k + 3) = 0$$

Thus, $k = -3$. Then, we use our value to substitute $k$ in our assumption which is $y(x) = e^{kx}$:

at $k = -3$, $y_1(x) = e^{(-3)x} = e^{-3x}$

Thus, using result 2.1.1, the general solution for $y(x)$ is: $y_{homo}(x) = c_1 e^{-3x}$, for some $c_1 \in \Re$. (Note: *homo* denotes to homogeneous).

**Step 2:** We need to find the particular solution as follows: Since $y'(x) + 3y(x)$ equals $x$, then the particular solution should be in the following form: $y_{particular}(x) = a + bx$ because $x$ is a polynomial of the first degree, and the general form for first degree polynomial is $a + bx$.

Now, we need to find $a$ and $b$ as follows:

$y_{particular}'(x) = b$

We substitute $y_{particular}'(x) = b$ in $y'(x) + 3y(x) = x$.

$$b + 3(a + bx) = x$$
$$b + 3a + 3bx = x$$

at $x = 0$, we obtain: $b + 3a + (3b)(0) = 0$

$b + 3a + 0 = 0$

$b + 3a = 0$ … … … … … … … … … … … … … … … … … … … … (1)

at $x = 1$, we obtain: $b + 3a + (3b)(1) = 1$

$b + 3a + 3b = 1$

$4b + 3a = 1$ … … … … … … … … … … . … … … … … … … … … (2)

From (1) and (2), we get: $a = -\frac{1}{9}$ and $b = \frac{1}{3}$.

Thus, $y_{particular}(x) = -\frac{1}{9} + \frac{1}{3}x$.

**Step 3:** We need to find the general solution as follows:

$y_{general}(x) = y_{homo}(x) + y_{particular}(x)$

Thus, $y_{general}(x) = c_1 e^{-3x} + \left(-\frac{1}{9} + \frac{1}{3}x\right)$.

**Example 2.2.2** Given $y'(x) + 3y(x) = e^{-3x}$, $y(0) = 1$. Find the general solution for $y(x)$. (Hint: No need to find the value of $c_1$)

**Solution:** In this example, we will have the same homogeneous solution as we did in example 2.2.1 but the only difference is the particular solution. We will repeat some steps in case you did not read example 2.2.1. Since $y'(x) + 3y(x) = e^{-3x}$ does not have a constant coefficient, then we need to use the undetermined coefficients method as follows:

**Step 1:** We need to find the homogeneous solution by letting $y'(x) + 3y(x)$ equal to zero as follows:

$y'(x) + 3y(x) = 0$. Now, it is a homogeneous linear differential equation of order 1.

Since $y'(x) + 3y(x) = 0$ is a HLDE with constant coefficients, then we will do the following:

Let $y(x) = e^{kx}$, we need to find $k$.

First of all, we will find the first, second, third, fourth, and fifth derivatives as follows:

$$y'(x) = ke^{kx}$$

Now, we substitute $y'(x) = ke^{kx}$ in $y'(x) + 3y(x) = 0$ as follows:

$$ke^{kx} + 3e^{kx} = 0$$

$$e^{kx}(k + 3) = 0$$

Thus, $k = -3$. Then, we use our value to substitute $k$ in our assumption which is $y(x) = e^{kx}$:

at $k = -3$, $y_1(x) = e^{(-3)x} = e^{-3x}$

Thus, using result 2.1.1, the general solution for $y(x)$ is: $y_{homo}(x) = c_1 e^{-3x}$, for some $c_1 \in \Re$. (Note: $homo$ denotes to homogeneous).

**Step 2:** We need to find the particular solution as follows: Since $y'(x) + 3y(x)$ equals $e^{-3x}$, then the particular solution should be in the following form:

$y_{particular}(x) = (ae^{-3x})x$

Now, we need to find $a$ as follows:

$y_{particular}'(x) = (ae^{-3x}) - 3(axe^{-3x})$

We substitute $y_{particular}'(x) = (ae^{-3x}) - 3(axe^{-3x})$ in $y'(x) + 3y(x) = e^{-3x}$.

$[(ae^{-3x}) - 3(axe^{-3x})] + 3(ae^{-3x})x = e^{-3x}$

$ae^{-3x} = e^{-3x}$

$a = 1$

Thus, $y_{particular}(x) = (1e^{-3x})x = xe^{-3x}$.

**Step 3:** We need to find the general solution as follows:

$y_{general}(x) = y_{homo}(x) + y_{particular}(x)$

Thus, $y_{general}(x) = c_1e^{-3x} + xe^{-3x}$.

**Result 2.2.1** Suppose that you have a linear differential equation with the least derivative, say $m$, and this differential equation equals to a polynomial of degree $w$. Then, we obtain the following:

$y_{particular} = [Polynomial\ of\ Degree\ w]x^m$ .

**Result 2.2.2** Suppose that you have a linear differential equation, then the general solution is always written as: $y_{general}(x) = y_{homogeneous}(x) + y_{particular}(x)$.

**Example 2.2.3** Given $y^{(4)}(x) - 7y^{(3)}(x) = x^2$. Describe $y_{particular}(x)$ but do not find it.

**Solution:** To describe $y_{particular}(x)$, we do the following:

By using result 2.2.1, we obtain the following:

$y_{particular} = [a + bx + cx^2]x^3$.

**Example 2.2.4** Given $y^{(4)}(x) - 7y^{(3)}(x) = 3$. Describe $y_{particular}(x)$ but do not find it.

**Solution:** To describe $y_{particular}(x)$, we do the following:

By using result 2.2.1, we obtain the following:

$y_{particular} = [a]x^3 = ax^3$.

**Example 2.2.5** Given $y^{(2)}(x) - 3y(x) = x^2 e^x$. Describe $y_{particular}(x)$ but do not find it.

**Solution:** To describe $y_{particular}(x)$, we do the following: In this example, we look at $x^2$, and we write it as: $a + bx + cx^2$, and then we multiply it by $e^x$.

$y_{particular} = [a + bx + cx^2]e^x$.

**Example 2.2.6** Given $y^{(2)}(x) - 3y(x) = \sin(3x)e^x$. Describe $y_{particular}(x)$ but do not find it.

**Solution:** To describe $y_{particular}(x)$, we do the following: In this example, we look at $\sin(3x)$, and we write it as: $(a\sin(3x) + b\sin(3x))$, and then we multiply it by $e^x$.

$y_{particular} = [a\sin(3x) + b\sin(3x)]e^x$.

# 2.3 Exercises

**1.** Given $y^{(2)}(x) + 2y'(x) + y(x) = 0$. Find the general solution for $y(x)$. (Hint: Use results 2.1.1 and 2.1.2, and in this exercise, no need to find the values of $c_1$, and $c_2$)

**2.** Given $y^{(3)}(x) - y^{(2)}(x) = 3x$. Find the general solution for $y(x)$. (Hint: No need to find the value of $c_1, c_2$, and $c_3$)

**3.** Given $y^{(4)}(x) - y^{(3)}(x) = 3x^2$. Find the general solution for $y(x)$. (Hint: No need to find the value of $c_1, c_2, c_3$, and $c_4$)

**4.** Given $y^{(3)}(x) - y^{(2)}(x) = e^x$. Find the general solution for $y(x)$. (Hint: No need to find the value of $c_1, c_2,$ and $c_3$)

**5.** Given $y^{(3)}(x) - y^{(2)}(x) = \sin(2x)$. Find the general solution for $y(x)$. (Hint: No need to find the value of $c_1, c_2,$ and $c_3$)

**6.** Given $y^{(2)}(x) - 3y(x) = 3x\sin(5x)$. Describe $y_{particular}(x)$ but do not find it.

**7.** Given $y^{(2)}(x) - 3y(x) = e^{x^2}$. Describe $y_{particular}(x)$ but do not find it.

# Chapter 3

# Methods of First and Higher Orders Differential Equations

In this chapter, we introduce two new methods called Variation Method and Cauchy-Euler Method in order to solve first and higher orders differential equations. In addition, we give several examples about these methods, and the difference between them and the previous methods in chapter 2.

## 3.1 Variation Method

In this section, we discuss how to find the particular solution using Variation Method. For the homogeneous solution, it will be similar to what we learned in chapter 2.

**Definition 3.1.1** Given $a_2(x)y^{(2)} + a_1(x)y' = K(x)$ is a linear differential equation of order 2. Assume that $y_1(x)$ and $y_2(x)$ are independent solution to the homogeneous solution. Then, the particular solution using *Variation Method* is written as:

$y_{particular}(x) = h_1(x)y_1(x) + h_2(x)y_2(x)$. To find $h_1(x)$ and $h_2(x)$, we need to solve the following two equations:

$$h_1'(x)y_1(x) + h_2'(x)y_2(x) = 0$$
$$h_1'(x)y_1'(x) + h_2'(x)y_2'(x) = \frac{K(x)}{a_2(x)}$$

**Definition 3.1.2** Given $a_3(x)y^{(3)} + \cdots + a_1(x)y' = K(x)$ is a linear differential equation of order 3. Assume that $y_1(x), y_2(x)$ and $y_3(x)$ are independent solution to the homogeneous solution. Then, the particular solution using *Variation Method* is written as:

$y_{particular}(x) = h_1(x)y_1(x) + h_2(x)y_2(x) + h_3(x)y_3(x)$. To find $h_1(x), h_2(x)$ and $h_3(x)$, we need to solve the following three equations:

$$h_1'(x)y_1(x) + h_2'(x)y_2(x) + h_3'(x)y_3(x) = 0$$
$$h_1'(x)y_1'(x) + h_2'(x)y_2'(x) + h_3'(x)y_3'(x) = 0$$
$$h_1^{(2)}(x)y_1^{(2)}(x) + h_2^{(2)}(x)y_2^{(2)}(x) + h_3^{(2)}(x)y_3^{(2)}(x) = \frac{K(x)}{a_3(x)}$$

**Example 3.1.1** Given $y^{(2)} + 3y' = \frac{1}{x}$. Find the general solution for $y(x)$. (Hint: No need to find the values of $c_1$ and $c_2$)

**Solution:** Since $y^{(2)} + 3y' = \frac{1}{x}$ does not have a constant coefficient, then we need to use the variation method as follows:

**Step 1:** We need to find the homogeneous solution by letting $y^{(2)} + 3y'$ equal to zero as follows:

$y^{(2)} + 3y' = 0$. Now, it is a homogeneous linear differential equation of order 2.

Since $y^{(2)} + 3y' = 0$ is a HLDE with constant coefficients, then we will do the following:

Let $y(x) = e^{kx}$, we need to find $k$.

First of all, we will find the first and second derivatives as follows:

$$y'(x) = ke^{kx}$$

$$y''(x) = k^2 e^{kx}$$

Now, we substitute $y'(x) = ke^{kx}$ and $y''(x) = k^2 e^{kx}$ in $y^{(2)} + 3y' = 0$ as follows:

$$k^2 e^{kx} + 3ke^{kx} = 0$$

$$e^{kx}(k^2 + 3k) = 0$$

$$e^{kx}(k(k + 3)) = 0$$

Thus, $k = 0$ and $k = -3$. Then, we use our values to substitute $k$ in our assumption which is $y(x) = e^{kx}$:

at $k = 0$, $y_1(x) = e^{(0)x} = e^0 = 1$

at $k = -3$, $y_2(x) = e^{(-3)x} = e^{-3x}$

Notice that $y_1(x)$ and $y_2(x)$ are independent.

Thus, using result 2.1.1, the general homogenous solution for $y(x)$ is: $y_{homo}(x) = c_1 + c_2 e^{-3x}$, for some $c_1$ and $c_2 \in \Re$. (Note: *homo* denotes to homogeneous).

**Step 2:** We need to find the particular solution using definition 3.1.1 as follows: Since $y^{(2)} + 3y'$ equals $\frac{1}{x}$, then the particular solution should be in the following form: $y_{particular}(x) = h_1(x)y_1(x) + h_2(x)y_2(x)$. To find $h_1(x)$ and $h_2(x)$, we need to solve the following two equations:

$$h_1{}'(x)y_1(x) + h_2{}'(x)y_2(x) = 0$$

$$h_1{}'(x)y_1{}'(x) + h_2{}'(x)y_2{}'(x) = \frac{K(x)}{a_2(x)}$$

$y_1(x) = 1 \quad \cdots\cdots\rightarrow \quad y_1{}'(x) = 0$

$y_2(x) = e^{-3x} \quad \cdots\cdots\rightarrow \quad y_2{}'(x) = -3e^{-3x}$

Now, we substitute what we got above in the particular solution form as follows:

$$h_1{}'(x)(1) + h_2{}'(x)(e^{-3x}) = 0$$

$$h_1{}'(x)(0) + h_2{}'(x)(-3e^{-3x}) = \frac{\frac{1}{x}}{1}$$

$h_1{}'(x)(1) + h_2{}'(x)(e^{-3x}) = 0 \dots\dots\dots\dots\dots\dots\dots (1)$

$h_2{}'(x)(-3e^{-3x}) = \frac{1}{x} \dots\dots\dots\dots\dots\dots\dots\dots (2)$

By solving (1) and (2), $h_2{}'(x) = -\frac{1}{3x}e^{3x}$ and

$h_1{}'(x) = \frac{1}{3x}e^{3x}(e^{-3x}) = \frac{1}{3x}$.

Since it is impossible to integrate $h_2{}'(x) = -\frac{1}{3x}e^{3x}$ to find $h_2(x)$, then it is enough to write as:

$h_2(x) = \int_0^x -\frac{1}{3t}e^{3t}\, dt$.

Since it is possible to integrate $h_1{}'(x) = \frac{1}{3x}$ to find $h_1(x)$, then we do the following:

$h_1(x) = \int \frac{1}{3x}\, dx = \frac{1}{3}(\ln|x|), x > 0$.

Thus, we write the particular solution as follows:

$$y_{particular}(x) = \frac{1}{3}(\ln|x|) + e^{-3x}\left(\int_0^x -\frac{1}{3t}e^{3t}\, dt\right)$$

**Step 3:** We need to find the general solution as follows:

$y_{general}(x) = y_{homo}(x) + y_{particular}(x)$

Thus, $y_{general}(x) = (c_1 + c_2 e^{-3x}) + \frac{1}{3}(\ln|x|) +$

$e^{-3x}\left(\int_0^x -\frac{1}{3t} e^{3t}\, dt\right)$, for some $c_1$ and $c_2 \in \Re$.

**Example 3.1.2** Given $y^{(2)} + 6y' + 8y = e^{-4x}$. Find the general solution for $y(x)$. (Hint: No need to find the values of $c_1$ and $c_2$)

**Solution:** Since $y^{(2)} + 6y' + 8y = e^{-4x}$ does not have a constant coefficient, then we need to use the variation method as follows:

**Step 1:** We need to find the homogeneous solution by letting $y^{(2)} + 6y' + 8y$ equal to zero as follows:

$y^{(2)} + 6y' + 8y = 0$. Now, it is a homogeneous linear differential equation of order 2.

Since $y^{(2)} + 6y' + 8y = 0$ is a HLDE with constant coefficients, then we will do the following:

Let $y(x) = e^{kx}$, we need to find $k$.

First of all, we will find the first and second derivatives as follows:

$$y'(x) = ke^{kx}$$

$$y''(x) = k^2 e^{kx}$$

Now, we substitute $y(x) = e^{kx}$, $y'(x) = ke^{kx}$ and $y''(x) = k^2 e^{kx}$ in $y^{(2)} + 6y' + 8y = 0$ as follows:

$$k^2 e^{kx} + 6ke^{kx} + 8e^{kx} = 0$$

$$e^{kx}(k^2 + 6k + 8) = 0$$

$$e^{kx}((k+2)(k+4)) = 0$$

Thus, $k = -2$ and $k = -4$. Then, we use our values to substitute $k$ in our assumption which is $y(x) = e^{kx}$:

at $k = -2$, $y_1(x) = e^{(-2)x} = e^{-2x}$

at $k = -4$, $y_2(x) = e^{(-4)x} = e^{-4x}$

Notice that $y_1(x)$ and $y_2(x)$ are independent.

Thus, using result 2.1.1, the general homogenous solution for $y(x)$ is: $y_{homo}(x) = c_1 e^{-2x} + c_2 e^{-4x}$, for some $c_1$ and $c_2 \in \Re$. (Note: $homo$ denotes to homogeneous).

**Step 2:** We need to find the particular solution using definition 3.1.1 as follows: Since $y^{(2)} + 6y' + 8y$ equals $e^{-4x}$, then the particular solution should be in the following form: $y_{particular}(x) = h_1(x)y_1(x) + h_2(x)y_2(x)$. To find $h_1(x)$ and $h_2(x)$, we need to solve the following two equations:

$$h_1'(x)y_1(x) + h_2'(x)y_2(x) = 0$$

$$h_1'(x)y_1'(x) + h_2'(x)y_2'(x) = \frac{K(x)}{a_2(x)}$$

$y_1(x) = e^{-2x}$ $\text{------>}$ $y_1'(x) = -2e^{-2x}$

$y_2(x) = e^{-4x}$ $\text{------>}$ $y_2'(x) = -4e^{-4x}$

Now, we substitute what we got above in the particular solution form as follows:

$$h_1'(x)(e^{-2x}) + h_2'(x)(e^{-4x}) = 0$$

$$h_1'(x)(-2e^{-2x}) + h_2'(x)(-4e^{-4x}) = \frac{e^{-4x}}{1}$$

$h_1'(x)(e^{-2x}) + h_2'(x)(e^{-4x}) = 0 \ldots\ldots\ldots\ldots\ldots\ldots\ldots (1)$

$h_1'(x)(-2e^{-2x}) + h_2'(x)(-4e^{-4x}) = e^{-4x} \ldots\ldots\ldots\ldots\ldots (2)$

By solving (1) and (2), and using Cramer's rule, we obtain:

$$h_1'(x) = x_1 = \frac{\det\begin{bmatrix} 0 & e^{-4x} \\ e^{-4x} & -4e^{-4x} \end{bmatrix}}{\det\begin{bmatrix} e^{-2x} & e^{-4x} \\ -2e^{-2x} & -4e^{-4x} \end{bmatrix}}$$

$$= \frac{-e^{-8x}}{-4e^{-6x} + 2e^{-6x}} = \frac{-e^{-8x}}{-2e^{-6x}}$$

$$= \frac{1}{2}e^{-2x}$$

By substituting $h_1'(x)$ in (1) to find $h_2'(x)$ as follows:

$$\left(\frac{1}{2}e^{-2x}\right)(e^{-2x}) + h_2'(x)(e^{-4x}) = 0$$

$$h_2'(x) = -\frac{1}{2}$$

Since it is possible to integrate $h_1'(x) = \frac{1}{2}e^{-2x}$ to find $h_1(x)$, then we do the following:

$$h_1(x) = \int \frac{1}{2}e^{-2x}\, dx = -\frac{1}{4}e^{-2x}.$$

Since it is possible to integrate $h_2'(x) = -\frac{1}{2}$ to find $h_2(x)$, then we do the following:

$$h_2(x) = \int -\frac{1}{2}dx = -\frac{1}{2}x.$$

Thus, we write the particular solution as follows:

$$y_{particular}(x) = -\frac{1}{4}e^{-2x}(e^{-2x}) - \frac{1}{2}x(e^{-4x})$$

**Step 3:** We need to find the general solution as follows:

$$y_{general}(x) = y_{homo}(x) + y_{particular}(x)$$

Thus, $y_{general}(x) = (c_1 e^{-2x} + c_2 e^{-4x}) + \left(-\frac{1}{4}e^{-2x}(e^{-2x}) - \frac{1}{2}x(e^{-4x})\right)$, for some $c_1$ and $c_2 \in \Re$.

# 3.2 Cauchy-Euler Method

In this section, we will show how to use Cauchy-Euler Method to find the general solution for differential equations that do not have constant coefficients.

To introduce this method, we start with some examples as follows:

**Example 3.2.1** Given $xy^{(2)} - y' + \frac{1}{x}y = 0$. Find the general solution for $y(x)$. (Hint: No need to find the values of $c_1$ and $c_2$)

**Solution:** Since $xy^{(2)} - y' + \frac{1}{x}y = 0$ does not have constant coefficients, then we need to use the Cauchy-Euler method by letting $y = x^k$, and after substitution all terms must be of the same degree as follows:

First of all, we will find the first and second derivatives as follows:

$$y' = kx^{k-1}$$
$$y'' = k(k-1)x^{k-2}$$

Now, we substitute $y = x^k,\ y' = kx^{k-1}$ and $y'' = k(k-1)x^{k-2}$ in $xy^{(2)} - y' + \frac{1}{x}y = 0$ as follows:

$$xk(k-1)x^{k-2} - kx^{k-1} + \frac{1}{x}x^k = 0$$
$$k(k-1)x^{k-1} - kx^{k-1} + x^{k-1} = 0$$
$$x^{k-1}(k(k-1) - k + 1) = 0$$
$$x^{k-1}(k^2 - k - k + 1) = 0$$

$$x^{k-1}(k^2 - 2k + 1) = 0$$

$$x^{k-1}((k-1)(k-1)) = 0$$

Thus, $k = 1$ and $k = 1$. Then, we use our values to substitute $k$ in our assumption which is $y = x^k$:

at $k = 1$, $y_1 = x^1 = x$

at $k = 1$, $y_2 = x^1 = x \cdot \ln(x)$

In the above case, we multiplied $x$ by $\ln(x)$ because we had a repeating for $x$, and in Cauchy-Euler Method, we should multiply any repeating by natural logarithm. Thus, the general solution for $y(x)$ is:

$y(x) = c_1 x + c_2 x ln(x)$, for some $c_1$ and $c_2 \in \Re$.

**Example 3.2.2** Given $x^3 y^{(2)} - x^2 y' + xy = 0$. Find the general solution for $y(x)$. (Hint: No need to find the values of $c_1$ and $c_2$)

**Solution:** Since $x^3 y^{(2)} + x^2 y' + xy = 0$ does not have constant coefficients, then we need to use the Cauchy-Euler method by letting $y = x^k$, and after substitution all terms must be of the same degree as follows:

First of all, we will find the first, second and third derivatives as follows:

$$y' = kx^{k-1}$$

$$y'' = k(k-1)x^{k-2}$$

$$y''' = k(k-1)(k-2)x^{k-3}$$

Now, we substitute $y = x^k$, $y' = kx^{k-1}$, and $y'' = k(k-1)x^{k-2}$ in $x^3 y^{(2)} + x^2 y' + xy = 0$ as follows:

$$x^3(k(k-1)x^{k-2}) + x^2(kx^{k-1}) + x(x^k) = 0$$
$$(k(k-1)x^{k+1}) + (kx^{k+1}) + (x^{k+1}) = 0$$
$$x^{k+1}(k^2 - k + k + 1) = 0$$
$$x^{k+1}(k^2 + 1) = 0$$

Thus, $k = \pm\sqrt{1} = \pm i = 0 \pm (1)(i)$.

Then, we use our values to substitute $k$ in our assumption which is $y = x^k$:

Since we have two parts (real and imaginary), then by using the Cauchy-Euler Method, we need to write our solution as follows:

$y_1 = x^{(real\ part)}\cos(Imaginary\ part \cdot \ln(x)$

$$= x^{(0)}\cos(1 \cdot \ln(x)) = \cos(\ln(x))$$

$y_2 = x^{(real\ part)}\sin(Imaginary\ part \cdot \ln(x)$

$$= x^{(0)}\sin(1 \cdot \ln(x)) = \sin(\ln(x))$$

Thus, the general solution for $y(x)$ is:

$y(x) = c_1\cos(\ln(x)) + c_2\sin(\ln(x))$, for some $c_1$ and $c_2 \in \Re$.

# 3.3 Exercises

1. Given $y^{(2)} + y' + 4y = 0$. Find the general solution for $y(x)$. (Hint: No need to find the values of $c_1$ and $c_2$)

2. Given $y^{(2)} + 5y' + 7y = 0$. Find the general solution for $y(x)$. (Hint: No need to find the values of $c_1$ and $c_2$)

3. Given $x^3y^{(3)} + xy' = 0$. Find the general solution for $y(x)$. (Hint: No need to find the values of $c_1$, $c_2$ and $c_3$)

**4.** Given $x^3 y^{(3)} - 2xy' = 0$. Find the general solution for $y(x)$. (Hint: No need to find the values of $c_1$, $c_2$ and $c_3$)

**5.** Given $x^2 y^{(2)} + y' = 2x^2$. Is it possible to find the general solution for $y(x)$ using Cauchy-Euler Method? Why?

# Chapter 4

# Extended Methods of First and Higher Orders Differential Equations

In this chapter, we discuss some new methods such as Bernulli Method, Separable Method, Exact Method, Reduced to Separable Method and Reduction of Order Method. We use  these methods to solve first and higher orders linear and non-linear differential equations. In addition, we give examples about these methods, and the differences between them and the previous methods in chapter 2 and chapter 3.

## 4.1 Bernoulli Method

In this section, we start with two examples about using integral factor to solve first order linear differential equations. Then, we introduce Bernulli Method to solve some examples of first order non-linear differential equations.

**Definition 4.1.1** Given $a_1(x)y' + a_2(x)y = K(x)$, where $a_1(x) \neq 0$ is a linear differential equation of order 1. Dividing both sides by $a_1(x)$, we obtain:

$$\frac{a_1(x)}{a_1(x)}y' + \frac{a_2(x)}{a_1(x)}y = \frac{K(x)}{a_1(x)}$$

$$y' + \frac{a_2(x)}{a_1(x)}y = \frac{K(x)}{a_1(x)}$$

Assume that $g(x) = \frac{a_2(x)}{a_1(x)}$ and $F(x) = \frac{K(x)}{a_1(x)}$, then:

$$y' + g(x)y = F(x) \dots \dots \dots \dots \dots \dots . (1)$$

Thus, the solution using *Integral Factor Method* is written in the following steps:

**Step 1:** Multiply both sides of (1) by letting $I = e^{\int g(x)dx}$:

$$y'e^{\int g(x)dx} + g(x)ye^{\int g(x)dx} = F(x)e^{\int g(x)dx}$$

**Step 2:** $y'e^{\int g(x)dx} + g(x)ye^{\int g(x)dx} = \left[y \cdot e^{\int g(x)dx}\right]' \dots (2)$

**Step 3:** $\left[y \cdot e^{\int g(x)dx}\right]' = F(x)e^{\int g(x)dx} \dots \dots \dots \dots \dots \dots . (3)$

**Step 4:** Integrate both sides of (3), we obtain:

$$\int \left[y \cdot e^{\int g(x)dx}\right]' dx = \int \left(F(x)e^{\int g(x)dx}\right) dx$$

$$y \cdot e^{\int g(x)dx} = \int \left(F(x)e^{\int g(x)dx}\right) dx$$

**Step 5: By** solving for $y$, and substituting $I = e^{\int g(x)dx}$ we obtain:

$$y \cdot I = \int (F(x))(I)\, dx$$

$$y = \frac{\int (F(x))(I)\, dx}{I}$$

$$y = \frac{\int I \cdot F(x)\, dx}{I}$$

Thus, the final solution is:

$$y = \frac{\int I \cdot F(x)\, dx}{I}$$

**Example 4.1.1** Given $x^2 y' - 2xy = 4x^3$. Find the general solution for $y(x)$. (Hint: Use integral factor method and no need to find the value of $c$)

**Solution:** Since $x^2 y' - 2xy = 4x^3$ does not have constant coefficients, and it is a first order non-linear differential equation, then by using definition 4.1.1, we need to use the integral factor method by letting

$I = e^{\int g(x)dx}$, where $g(x) = \frac{a_2(x)}{a_1(x)} = -\frac{2x}{x^2} = -\frac{2}{x}$ and $F(x) = \frac{K(x)}{a_1(x)} = \frac{4x^3}{x^2} = 4x$

Hence, $I = e^{\int g(x)dx} = e^{\int -\frac{2}{x}dx} = e^{-2\ln(x)} = e^{\ln(x^{-2})} = \frac{1}{x^2}$

The general solution is written as follows:

$$y = \frac{\int I \cdot F(x)\,dx}{I}$$

$$y = \frac{\int \frac{1}{x^2} \cdot 4x\,dx}{\frac{1}{x^2}}$$

$$y = \frac{\int \frac{4}{x}\,dx}{\frac{1}{x^2}}$$

$$y = \frac{4\ln(x) + c}{\frac{1}{x^2}}$$

$$y = 4x^2 \ln(x) + cx^2$$

Thus, the general solution is: $y = 4x^2 \ln(x) + cx^2$ for some $c \in \Re$.

**Example 4.1.2** Given $(x + 1)y' + y = 5$. Find the general solution for $y(x)$. (Hint: Use integral factor method and no need to find the value of $c$)

**Solution:** Since $(x + 1)y' + y = 5$ does not have constant coefficients, and it is a first order non-linear differential equation, then by using definition 4.1.1, we need to use the integral factor method by letting

$I = e^{\int g(x)dx}$, where $g(x) = \frac{a_2(x)}{a_1(x)} = \frac{1}{(x+1)}$ and

$$F(x) = \frac{K(x)}{a_1(x)} = \frac{5}{(x + 1)}$$

Hence, $I = e^{\int g(x)dx} = e^{\int \frac{1}{(x+1)} dx} = e^{\ln(x+1)} = (x + 1)$

The general solution is written as follows:

$$y = \frac{\int I \cdot F(x) \, dx}{I}$$

$$y = \frac{\int (x + 1) \cdot \frac{5}{(x + 1)} dx}{(x + 1)}$$

$$y = \frac{\int 5 \, dx}{(x + 1)}$$

$$y = \frac{5x + c}{(x + 1)}$$

$$y = \frac{5x}{(x + 1)} + \frac{c}{(x + 1)}$$

Thus, the general solution is: $y = \frac{5x}{(x+1)} + \frac{c}{(x+1)}$

for some $c \in \Re$.

**Definition 4.1.2** Given $y' + g(x)y = f(x)y^n$ where $n \in \Re$ and $n \neq 0$ and $n \neq 1$ is a non-linear differential equation of order 1. Thus, the solution using *Bernoulli Method* is written in the following steps:

**Step 1:** Change it to first order linear differential equation by letting $w = y^{1-n}$.

**Step 2:** Find the derivative of both sides for $w = y^{1-n}$ as follows:

$$\frac{dw}{dx} = (1-n)y^{1-n-1} \cdot \frac{dy}{dx}$$

$$\frac{dw}{dx} = (1-n)y^{-n} \cdot \frac{dy}{dx} \dots \dots \dots \dots \dots \dots \dots \dots \dots \dots (1)$$

**Step 3:** Solve (1) for $\frac{dy}{dx}$ as follows:

$$\frac{dy}{dx} = \frac{1}{1-n}y^n \cdot \frac{dw}{dx} \dots \dots \dots \dots \dots \dots \dots \dots \dots \dots (2)$$

**Step 4:** Since we assumed that $w = y^{1-n}$, then $y = w^{\left(\frac{1}{1-n}\right)}$, and hence $y^n = w^{\left(\frac{n}{1-n}\right)}$.

**Step 5:** Substitute what we got above in $\frac{dy}{dx} + g(x)y = f(x)y^n$ as follows:

$$\frac{1}{1-n}y^n \cdot \frac{dw}{dx} + g(x)w^{\left(\frac{1}{1-n}\right)} = f(x)w^{\left(\frac{n}{1-n}\right)} \dots \dots \dots \dots \dots (3)$$

**Step 6:** Divide (3) by $\frac{1}{1-n}y^n$ as follows:

$$\frac{dw}{dx} + \frac{g(x)w^{\left(\frac{1}{1-n}\right)}}{\frac{1}{1-n}y^n} = \frac{f(x)w^{\left(\frac{n}{1-n}\right)}}{\frac{1}{1-n}y^n}$$

$$\frac{dw}{dx} + g(x)(1-n)w^{\left(\frac{1}{1-n}\right)}y^{-n} = f(x)(1-n)w^{\left(\frac{n}{1-n}\right)}y^{-n}$$

**Step 7:** After substitution, we obtain:

$$\frac{dw}{dx} + g(x)(1-n)yy^{-n} = f(x)(1-n)y^ny^{-n}$$

$$\frac{dw}{dx} + g(x)(1-n)y^{1-n} = f(x)(1-n)y^{n-n}$$

$$\frac{dw}{dx} + g(x)(1-n)y^{1-n} = f(x)(1-n)y^0$$

$$\frac{dw}{dx} + g(x)(1-n)y^{1-n} = f(x)(1-n)$$

**Step 8:** We substitute $w = y^{1-n}$ in the above equation as follows:

$$\frac{dw}{dx} + g(x)(1-n)w = f(x)(1-n)$$

Thus, the final solution is:
$$\frac{dw}{dx} + (1-n)g(x)w = (1-n)f(x)$$
In the following example, we will show how to use Bernoulli Method, and we will explore the relationship between Bernoulli Method and Integral Factor Method.

**Example 4.1.3** Given $xy' + 3x^2y = (6x^2)y^3$. Find the general solution for $y(x)$. (Hint: Use Bernoulli method and no need to find the value of $c$)

**Solution:** Since $xy' + 3x^2y = (6x^2)y^3$ does not have constant coefficients, and it is a first order non-linear differential equation, then by using definition 4.1.2, we need to do the following by letting $w = y^{1-n}$, where in this example $n = 3$, and $g(x) = 3x^2$ and $f(x) = 6x^2$.

Since we assumed that $w = y^{1-3} = y^{-2}$, then $y = w^{\left(\frac{1}{1-n}\right)} = w^{\left(\frac{1}{1-3}\right)} = w^{-\frac{1}{2}} = \frac{1}{\sqrt{w}}$, and $\frac{dy}{dx} = -\frac{1}{2}y^3 \cdot \frac{dw}{dx}$

We substitute what we got above in $xy' + 3x^2y = (6x^2)y^3$ as follows:

$$x\left(-\frac{1}{2}y^3 \cdot \frac{dw}{dx}\right) + 3x^2\left(\frac{1}{\sqrt{w}}\right) = (6x^2)\left(\frac{1}{w\sqrt{w}}\right) \ldots \ldots \ldots \ldots \ldots (1)$$

Now, we divide (1) by $-\frac{1}{2}xy^3$ as follows:

$$\frac{dw}{dx} + \frac{3x^2\left(\frac{1}{\sqrt{w}}\right)}{-\frac{1}{2}xy^3} = \frac{(6x^2)y^3}{-\frac{1}{2}xy^3}$$

$$\frac{dw}{dx} + (-2)3x\left(\frac{1}{\sqrt{w}}\right)y^{-3} = (-2)6x \ldots \ldots \ldots \ldots \ldots \ldots \ldots \ldots (2)$$

Then, we substitute $y = \frac{1}{\sqrt{w}}$ in (2) as follows:

$$\frac{dw}{dx} + (-2)3x(y)y^{-3} = (-2)6x$$

$$\frac{dw}{dx} + (-2)3xy^{1-3} = (-2)6x$$

$$\frac{dw}{dx} + (-2)3xy^{-2} = (-2)6x \dots \dots \dots \dots \dots \dots \dots \dots \dots \dots \dots (3)$$

Now, we substitute $w = y^{-2}$ in (3) as follows:

$$\frac{dw}{dx} + (-2)3xw = (-2)6x$$

$$\frac{dw}{dx} - 6xw = -12x \dots \dots \dots \dots \dots \dots \dots \dots \dots \dots \dots \dots. \dots. (4)$$

Then, we solve (4) for $w(x)$ as follows:

To solve (4), we need to use the integral factor method:

Hence, $I = e^{\int g(x)dx} = e^{\int -6x\,dx} = e^{-\frac{6x^2}{2}} = e^{-3x^2}$

The general solution for $w(x)$ is written as follows:

$$w = \frac{\int I \cdot F(x)\,dx}{I}$$

$$w = \frac{\int I \cdot (-12x)\,dx}{I}$$

$$w = \frac{\int e^{-3x^2} \cdot (-12x)\,dx}{e^{-3x^2}}$$

$$w = \frac{2\int e^{-3x^2} \cdot (-6x)\,dx}{e^{-3x^2}}$$

$$w = \frac{2e^{-3x^2} + c}{e^{-3x^2}}$$

$$w = \frac{2e^{-3x^2}}{e^{-3x^2}} + \frac{c}{e^{-3x^2}}$$

$$w = 2 + \frac{c}{e^{-3x^2}}$$

$$w = 2 + ce^{3x^2}$$

The general solution for $w(x)$ is: $w(x) = 2 + ce^{3x^2}$.

Thus, the general solution for $y(x)$ is:

$$y(x) = \frac{1}{\sqrt{w(x)}} = \frac{1}{\sqrt{2 + ce^{3x^2}}}$$

for some $c \in \Re$.

# 4.2 Separable Method

In this section, we will solve some differential equations using a method known as Separable Method. This method is called separable because we separate two different terms from each other.

**Definition 4.2.1** The standard form of *Separable Method* is written as follows:

$(All\ in\ terms\ of\ x)dx - (All\ in\ terms\ of\ y)dy = 0$

**Note:** it does not matter whether it is the above form or in the following form:

$(All\ in\ terms\ of\ y)dy - (All\ in\ terms\ of\ x)dx = 0$

**Example 4.2.1** Solve the following differential equation: $\frac{dy}{dx} = \frac{y^3}{(x+3)}$

**Solution:** By using definition 4.2.1, we need to rewrite the above equation in a way that each term is separated from the other term as follows:

$$\frac{dy}{dx} = \frac{y^3}{(x+3)} = \frac{\frac{1}{(x+3)}}{\frac{1}{y^3}} \ldots\ldots\ldots\ldots\ldots\ldots\ldots\ldots\ldots\ldots\ldots. (1)$$

Now, we need to do a cross multiplication for (1) as follows:

$$\frac{1}{y^3} dy = \frac{1}{(x+3)} dx$$

$$\frac{1}{y^3} dy - \frac{1}{(x+3)} dx = 0 \ldots\ldots\ldots\ldots\ldots\ldots\ldots\ldots.\ldots\ldots\ldots (2)$$

Then, we integrate both sides of (2) as follows:

$$\int \left( \frac{1}{y^3} dy - \frac{1}{(x+3)} dx \right) = \int 0$$

$$\int \left( \frac{1}{y^3} \right) dy - \int \left( \frac{1}{(x+3)} \right) dx = c$$

$$\int (y^{-3}) dy - \int \left( \frac{1}{(x+3)} \right) dx = c$$

$$-\frac{1}{2} y^{-2} - \ln(|(x+3)|) = c$$

Thus, the general solution is :

$$-\frac{1}{2} y^{-2} - \ln(|(x+3)|) = c$$

**Example 4.2.2** Solve the following differential

equation: $\frac{dy}{dx} = e^{3y+2x}$

**Solution:** By using definition 4.2.1, we need to rewrite the above equation in a way that each term is separated from the other term as follows:

$$\frac{dy}{dx} = e^{3y+2x} = e^{3y} \cdot e^{2x} = \frac{e^{2x}}{e^{-3y}} \dots \dots \dots \dots \dots \dots \dots \dots \dots (1)$$

Now, we need to do a cross multiplication for (1) as follows:

$$e^{-3y} dy = e^{2x} dx$$

$$e^{-3y} dy - e^{2x} dx = 0 \dots \dots \dots \dots \dots \dots \dots \dots \dots \dots \dots \dots \dots (2)$$

Then, we integrate both sides of (2) as follows:

$$\int (e^{-3y} dy - e^{2x} dx) = \int 0$$

$$\int (e^{-3y}) dy - \int (e^{2x}) \, dx = c$$

$$-\frac{1}{3} e^{-3y} - \frac{1}{2} e^{2x} = c$$

Thus, the general solution is :

$$-\frac{1}{3} e^{-3y} - \frac{1}{2} e^{2x} = c$$

# 4.3 Exact Method

In this section, we will solve some differential equations using a method known as Exact Method. In other words, this method is called the Anti-Implicit Derivative Method.

**Definition 4.3.1** The standard form of *Exact Method* is written as follows:

$$\frac{dy}{dx} = -\frac{F_x}{F_y} \dots \dots \dots \dots \dots \dots \dots \dots \dots \dots \dots \dots \dots \dots \dots \dots \dots (1)$$

Then, we solve (1) to find $F(x, y)$, and our general solution will be as follows: $F(x, y) = c$ for some constant $c \in \Re$. In other words, the standard form for exact first order differential equation is: $F_y dy + F_x dx = 0$, and it is considered exact if $F_{xy} = F_{yx}$.

**Note:** $F_x(x, y)$ is defined as the first derivative with respect to $x$ and considering $y$ as a constant, while $F_y(x, y)$ is defined as the first derivative with respect to $y$ and considering $x$ as a constant.

**Example 4.3.1** Given $x^2 + y^2 - 4 = 0$. Find $\frac{dy}{dx}$.

**Solution:** By using definition 4.3.1, we first find $F_x(x, y)$ by finding the first derivative with respect to $x$ and considering $y$ as a constant as follows: $F_x(x, y) = 2x$. Then, we find $F_y(x, y)$ by finding the first derivative with respect to $y$ and considering $x$ as a constant as follows: $F_y(x, y) = 2y$.

Thus, $\frac{dy}{dx} = -\frac{F_x}{F_y} = -\frac{2x}{2y} = -\frac{x}{y}$.

**Example 4.3.2** Given $y^3 e^x + 3xy^2 - x^3 + xy - 13 = 0$. Find $\frac{dy}{dx}$.

**Solution:** By using definition 4.3.1, we first find $F_x(x, y)$ by finding the first derivative with respect to $x$ and considering $y$ as a constant as follows: $F_x(x, y) = y^3 e^x + 3y^2 - 3x^2 + y$.

Then, we find $F_y(x, y)$ by finding the first derivative with respect to $y$ and considering $x$ as a constant as follows: $F_y(x, y) = 3y^2 e^x + 6xy + x$.

Thus, $\frac{dy}{dx} = -\frac{F_x}{F_y} = -\frac{(y^3 e^x + 3y^2 - 3x^2 + y)}{(3y^2 e^x + 6xy + x)}$.

**Example 4.3.3** Solve the following differential

equation: $(-3x + y)dy - (5x + 3y)dx = 0$

**Solution:** First of all, we need to check for the exact method as follows: We rewrite the above differential equation according to definition 4.3.1:

$$(-3x + y)dy + -(5x + 3y)dx = 0$$

$(-3x + y)dy + (-5x - 3y)dx = 0 \dots \dots \dots \dots \dots \dots \dots (1)$

Thus, from the above differential equation, we obtain:

$F_x(x, y) = (-5x - 3y)$ and $F_y(x, y) = (-3x + y)$

Now, we need to check for the exact method by finding

$F_{xy}(x, y) = F_{yx}(x, y)$ as follows:

We first find $F_{xy}(x, y)$ by finding the first derivative of $F_x(x, y)$ with respect to $y$ and considering $x$ as a constant as follows:

$F_{xy}(x, y) = -3$

Then, we find $F_{yx}(x, y)$ by finding the first derivative of $F_y(x, y)$ with respect to $x$ and considering $y$ as a constant as follows:

$F_{yx}(x, y) = -3$

Since $F_{xy}(x, y) = F_{yx}(x, y) = -3$, then we can use the exact method.

Now, we choose either $F_x(x, y) = (-5x - 3y)$ or $F_y(x, y) = (-3x + y)$, and then we integrate. We will choose $F_y(x, y) = (-3x + y)$ and we will integrate it as follows:

$$\int F_y(x, y)dy = \int (-3x + y)dy = -3xy + \frac{1}{2}y^2 + D(x) \dots (2)$$

$F(x, y) = -3xy + \frac{1}{2}y^2 + D(x)$

We need to find $D(x)$ as follows:

Since we selected $F_y(x, y) = (-3x + y)$ previously for integration, then we need to find $F_x(x, y)$ for (2) as follows:

$F_x(x, y) = -3y + D'(x) \ldots \ldots \ldots \ldots \ldots \ldots \ldots \ldots \ldots \ldots (3)$

Now, we substitute $F_x(x, y) = (-5x - 3y)$ in (3) as follows:

$(-5x - 3y) = -3y + D'(x)$
$D'(x) = -5x - 3y + 3y = -5x \ldots \ldots \ldots \ldots \ldots \ldots \ldots (4)$

Then, we integrate both sides of (4) as follows:

$$\int D'(x) dx = \int -5x dx$$

$$D(x) = \int -5x dx = -\frac{5}{2}x^2$$

Thus, the general solution of the exact method is :

$$F(x, y) = c$$

$$-3xy + \frac{1}{2}y^2 + -\frac{5}{2}x^2 = c$$

# 4.4 Reduced to Separable Method

In this section, we will solve some differential equations using a method known as Reduced to Separable Method.

**Definition 4.4.1** The standard form of *Reduced to Separable Method* is written as follows:

$\frac{dy}{dx} = f(ax + by + c)$ where $a, b \neq 0$.

**Example 4.4.1** Solve the following differential equation: $\frac{dy}{dx} = \frac{\sin(5x+y)}{\cos(5x+y) - 2\sin(5x+y)} - 5$.

**Solution:** By using definition 4.4.1, we first let $u = 5x + y$, and then, we need to find the first derivative of both sides of $u = 5x + y$.

$$\frac{du}{dx} = 5 + \frac{dy}{dx} \ldots \ldots \ldots \ldots \ldots \ldots \ldots \ldots \ldots \ldots \ldots \ldots (1)$$

Now, we solve (1) for $\frac{dy}{dx}$ as follows:

$$\frac{dy}{dx} = \frac{du}{dx} - 5 \ldots \ldots \ldots \ldots \ldots \ldots \ldots \ldots \ldots \ldots \ldots (2)$$

Then, we substitute $u = 5x + y$ and (2) in

$\frac{dy}{dx} = \frac{\sin(5x+y)}{\cos(5x+y) - 2\sin(5x+y)} - 5$ as follows:

$$\frac{du}{dx} - 5 = \frac{\sin(u)}{\cos(u) - 2\sin(u)} - 5$$

$$\frac{du}{dx} = \frac{\sin(u)}{\cos(u) - 2\sin(u)} \ldots \ldots \ldots \ldots \ldots \ldots \ldots \ldots \ldots (3)$$

Now, we can use the separable method to solve (3) as follows:

By using definition 4.2.1, we need to rewrite (3) in a way that each term is separated from the other term as follows:

$$\frac{du}{dx} = \frac{\sin(u)}{\cos(u) - 2\sin(u)} = \frac{1}{\dfrac{\cos(u) - 2\sin(u)}{\sin(u)}} \ldots \ldots \ldots \ldots \ldots (4)$$

Now, we need to do a cross multiplication for (4) as follows:

$$\frac{\cos(u) - 2\sin(u)}{\sin(u)} du = 1 dx$$

$$\frac{\cos(u) - 2\sin(u)}{\sin(u)} du - 1 dx = 0 \ldots \ldots \ldots \ldots \ldots \ldots \ldots \ldots (5)$$

Then, we integrate both sides of (5) as follows:

$$\int \left( \frac{\cos(u) - 2\sin(u)}{\sin(u)} du - 1 dx \right) = \int 0$$

$$\int \left(\frac{cos(u) - 2\sin(u)}{\sin(u)}\right) du - \int (1)\, dx = c$$

$$\int \left(\frac{cos(u)}{\sin(u)} - \frac{2\sin(u)}{\sin(u)}\right) du - x = c$$

$$\int \left(\frac{cos(u)}{\sin(u)} - 2\right) du - x = c$$

$$\ln(|\sin(u)|) - 2u - x = c \dots \dots \dots \dots \dots \dots \dots \dots \dots \dots \dots . (6)$$

Now, we substitute $u = 5x + y$ in (6) as follows:

$$\ln(|\sin(5x + y)|) - 2(5x + y) - x = c$$

Thus, the general solution is :

$$\ln(|\sin(5x + y)|) - 2(5x + y) - x = c$$

# 4.5 Reduction of Order Method

In this section, we will solve differential equations using a method called Reduction of Order Method.

**Definition 4.5.1** *Reduction of Order Method* is valid method only for second order differential equations, and one solution to the homogenous part must be given. For example, given $(x + 1)y^{(2)}(x) - y'(x) = 0$, and $y_1(x) = 1$. To find $y_2(x)$, the differential equation must be written in the standard form (Coefficient of $y^{(2)}$ must be 1) as follows:

$$\frac{(x+1)}{(x+1)}y^{(2)}(x) - \frac{y'(x)}{(x+1)} = \frac{0}{(x+1)}$$

$$1y^{(2)}(x) - \frac{y'(x)}{(x+1)} = 0$$

$$y^{(2)}(x) + \frac{-1}{(x+1)}y'(x) = 0 \ldots\ldots\ldots\ldots\ldots\ldots\ldots\ldots\ldots\ldots\ldots (1)$$

Now, let $M(x) = \frac{-1}{(x+1)}$, and substitute it in (1) as follows:

$$y^{(2)}(x) + M(x)y'(x) = 0$$

Hence, $y_2(x)$ is written as follows:

$$y_2(x) = y_1(x) \cdot \int \frac{e^{\int -M(x)dx}}{y_1{}^2(x)} dx$$

In our example, $y_2(x) = 1 \cdot \int \frac{e^{\int \frac{1}{(x+1)}dx}}{(1)^2} dx =$

$1 \cdot \int \frac{e^{\ln(x+1)}}{1} dx = \int e^{\ln(x+1)} dx = \int (x+1)\,dx = \frac{1}{2}x^2 + x$.

Thus, the homogenous solution is written as follows:

$y_{homogenous}(x) = c_1 + c_2\left(\frac{1}{2}x^2 + x\right)$, for some $c_1, c_2 \in \Re$.

**Example 4.5.1** Given the following differential equation: $xy^{(2)}(x) + (x+1)y'(x) - (2x+1)y = xe^{7x}$, and $y_1(x) = e^x$ is a solution to the associated homogenous part. Find $y_{homogenous}(x)$? (Hint: Find first $y_2(x)$, and then write $y_{homogenous}(x)$)

**Solution:** By using definition 4.5.1, To find $y_2(x)$, the differential equation must be equal to zero and must also be written in the standard form (Coefficient of $y^{(2)}$ must be 1) as follows:

We divide both sides of

$xy^{(2)}(x) + (x + 1)y'(x) - (2x + 1)y = 0$ by $x$ as follows:

$$\frac{x}{x}y^{(2)}(x) + \frac{(x + 1)}{x}y'(x) - \frac{2x + 1}{x}y = \frac{0}{x}$$

$1y^{(2)}(x) + \frac{(x + 1)}{x}y'(x) - \frac{2x + 1}{x}y = 0 \ldots \ldots \ldots \ldots \ldots \ldots (1)$

Now, let $(x) = \frac{(x+1)}{x} = \left(1 + \frac{1}{x}\right)$, and substitute it in (1) as follows:

$$y^{(2)}(x) + M(x)y'(x) - \frac{2x + 1}{x}y = 0$$

Hence, $y_2(x)$ is written as follows:

$$y_2(x) = y_1(x) \cdot \int \frac{e^{\int -M(x)dx}}{y_1{}^2(x)} dx$$

In example 4.5.1, $y_2(x) = e^x \cdot \int \frac{e^{\int \left(1 + \frac{1}{x}\right)dx}}{(e^x)^2} dx =$

$$e^x \cdot \int \frac{e^{-x} \cdot e^{\ln\left(\frac{1}{x}\right)}}{e^{2x}} dx = e^x \cdot \int \frac{e^{-3x}}{x} dx$$

Since it is impossible to integrate $e^x \cdot \int \frac{e^{-3x}}{x} dx$, then it is enough to write it as: $e^x \cdot \int_0^x \frac{e^{-3t}}{t} dt$.

Therefore, $y_2(x) = e^x \cdot \int_0^x \frac{e^{-3t}}{t} dt$.

Thus, the homogenous solution is written as follows:

$y_{homogenous}(x) = c_1 e^x + c_2 \left(e^x \cdot \int_0^x \frac{e^{-3t}}{t} dt\right)$, for some $c_1, c_2 \in \Re$.

# 4.6 Exercises

**1.** Given $(x + 1)y' + xy = \frac{(x+1)^4}{y^2}$. Find the general solution for $y(x)$. (Hint: Use Bernoulli method and no need to find the value of $c$)

**2.** Given $x' + 3yx = 3y^3$. Find the general solution for $y(x)$. (Hint: Use Bernoulli method and no need to find the value of $c$)

**3.** Solve the following differential equation: $\frac{dy}{dx} = \frac{1+y^2}{1+x^2}$

**4.** Solve the following differential equation: $\frac{dy}{dx} = \frac{1}{3x+x^2 y}$

**5.** Solve the following differential equation:

$$\frac{dy}{dx} = 3xe^{(x+5y)}$$

**6.** Solve the following differential equation:

$$(e^x y + 3yx - 2)dy + \left(\frac{1}{2}e^x y^2 + \frac{3}{2}y^2 + x^2\right)dx = 0$$

**7.** Solve the following differential equation:

$$\frac{dy}{dx} = \frac{\sin(5x + y)}{\cos(5x + y) - 2\sin(5x + y)} - 5$$

**8.** Given the following differential equation:
$(x + 1)y^{(2)}(x) - y'(x) = 10$, and $y_1(x) = 1$ is a solution to the associated homogenous part. Find $y_{particular}(x)$?

# Chapter 5

# Applications of Differential Equations

In this chapter, we give examples of three different applications of differential equations: temperature, growth and decay, and water tank. In each section, we give one example of each of the above applications, and we discuss how to use what we have learned previously in this book to solve each problem.

## 5.1 Temperature Application

In this section, we give an example of temperature application, and we introduce how to use one of the differential equations methods to solve it.

**Example 5.1.1** Thomas drove his car from Pullman, WA to Olympia, WA, and the outside air temperature was constant 104°F. During his trip, he took a break at Othello, WA gas station, and then he switched off the engine of his car, and checked his car temperature gauge, and it was 144°F. After ten minutes, Thomas checked his car temperature gauge, and it was 136°F.

a) How long will it take for the temperature of the engine to cool to 117°F?

b) What will be the temperature of the engine 30 minutes from now?

**Solution: Part a:** To determine how long will it take for the temperature of the engine to cool to 117°F, we need to do the following:

Assume that $T(t)$ is the temperature of engine at the time $t$, and $T_0$ is the constant outside air temperature.

Now, we need to write the differential equation for this example as follows:

$$\frac{dT}{dt} = \beta(T - T_0) \ldots \ldots \ldots \ldots \ldots \ldots \ldots \ldots \ldots \ldots \ldots . \ldots \ldots \ldots \ldots \ldots \ldots (1)$$

where $\beta$ is a constant.

From (1), we can write as follows:

$$T' = \beta T - \beta T_0$$

$$T' - \beta T = -\beta T_0 \ldots \ldots \ldots \ldots \ldots \ldots \ldots \ldots \ldots \ldots \ldots \ldots \ldots . (2)$$

From this example, it is given the following:

$T(0) = 144°F$, $T(10) = 136°F$, and $T_0 = 104°F$

From (2), $-\beta T_0$ is constant, and the dependent variable is $T$, while the independent variable is the time $t$.

By substituting $T_0 = 104°F$ in (2), we obtain:

$$T' - \beta T = -104\beta \ldots \ldots \ldots \ldots \ldots \ldots \ldots \ldots \ldots \ldots \ldots \ldots \ldots . (3)$$

Since (3) is a first order linear differential equation, then by using definition 4.1.1, we need to use the integral factor method by letting $I = e^{\int g(t)dt}$, where $g(t) = -\beta$ and $F(t) = -104\beta$.

Hence, $I = e^{\int g(t)dt} = e^{\int -\beta \, dt} = e^{-\beta t}$.

The general solution is written as follows:

$$T(t) = \frac{\int I \cdot F(t) \, dt}{I}$$

$$T(t) = \frac{\int e^{-\beta t} \cdot (-104\beta) \, dt}{e^{-\beta t}}$$

$$T(t) = \frac{\int (-104\beta) e^{-\beta t} \, dt}{e^{-\beta t}}$$

$$T(t) = \frac{104 e^{-\beta t} + c}{e^{-\beta t}}$$

$$T(t) = \frac{104 e^{-\beta t}}{e^{-\beta t}} + \frac{c}{e^{-\beta t}}$$

$$T(t) = 104 + \frac{c}{e^{-\beta t}}$$

$T(t) = 104 + ce^{\beta t} \dots \dots \dots \dots \dots \dots \dots \dots \dots \dots \dots \dots \dots \dots (4)$

The general solution is: $T(t) = 104 + ce^{\beta t}$ for some $c \in \Re$.

Now, we need to find $c$ by substituting $T(0) = 144°F$ in (4) as follows:

$$T(0) = 104 + ce^{\beta(0)}$$

$$144 = 104 + ce^0$$

$$144 = 104 + c(1)$$

$$144 = 104 + c$$

$$c = 144 - 104 = 40$$

Thus, $T(t) = 104 + 40e^{\beta t} \dots \dots \dots \dots \dots \dots \dots \dots \dots \dots \dots \dots (5)$

By substituting $T(10) = 136°F$ in (5), we obtain:

$$T(10) = 104 + 40e^{\beta(10)}$$

$$136 = 104 + 40e^{\beta(10)}$$

$$136 - 104 = 40e^{\beta(10)}$$

$$32 = 40e^{\beta(10)}$$

$$e^{\beta(10)} = \frac{32}{40} = 0.8 \ldots\ldots\ldots\ldots\ldots\ldots\ldots\ldots\ldots\ldots\ldots. (6)$$

By taking the natural logarithm for both sides of (6), we obtain:

$$ln\big(e^{\beta(10)}\big) = ln(0.8)$$

$$\beta(10) = ln(0.8)$$

$$\beta = \frac{l\,n(0.8)}{10} = -0.0223 \ldots\ldots\ldots\ldots\ldots\ldots\ldots. (7)$$

Now, we substitute (7) in (5) as follows:

$$T(t) = 104 + 40e^{-0.0223t} \ldots\ldots\ldots\ldots\ldots\ldots. (8)$$

Then, we need to find the time $t$ when $T(t) = 117°F$ by substituting it in (8) as follows:

$$117 = 104 + 40e^{-0.0223(t)}$$

$$117 - 104 = 40e^{-0.0223(t)}$$

$$3 = 40e^{-0.0223(t)}$$

$$e^{-0.0223(t)} = \frac{3}{40} = 0.075 \ldots\ldots\ldots\ldots\ldots\ldots. (9)$$

By taking the natural logarithm for both sides of (9), we obtain:

$$ln\big(e^{-0.0223(t)}\big) = ln(0.075)$$

$$-0.0223(t) = ln(0.075)$$

$$t = \frac{ln(0.075)}{-0.0223} \approx 116.16 \text{ minutes}$$

Thus, the temperature of the engine will take approximately 116.16 minutes to cool to 117°F

**Part b:** To determine what will be the temperature of the engine 30 minutes from now, we need to do the following:

We assume that $t = 30$, and then we substitute it in (8) as follows:

$$T(30) = 104 + 40e^{-0.0223(30)}$$

$$T(30) \approx 124.49°F$$

Thus, the temperature of the engine 30 minutes from now will be approximately 124.49°F.

# 5.2 Growth and Decay Application

In this section, we give an example of growth and decay application, and we introduce how to use one of the differential equations methods to solve it.

**Example 5.2.1** The rate change of number of students at *Washington State University* (WSU) is proportional to the square root of the number of students at any time $t$. If the number of WSU students in 2013 was 28,686 students[2], and suppose that the number of students at WSU after one year was 32,000 students.

a) How long will it take to double number of WSU students in 2013?

b) What will be the number of WSU students in 2018?

**Solution: Part a:** To determine how long will it take to double number of WSU students in 2013, we need to do the following:

Assume that $W(t)$ is the number of WSU students at any time $t$.

Now, we need to write the differential equation for this example as follows:

$$\frac{dW}{dt} = \beta \sqrt{W(t)} \dots \dots \dots \dots \dots \dots \dots \dots \dots \dots \dots \dots \dots \dots \dots (1)$$

where $\beta$ is a constant.

From (1), we can write as follows:

$$W' = \beta \sqrt{W(t)} \dots \dots \dots \dots \dots \dots \dots \dots \dots \dots \dots \dots \dots \dots \dots (2)$$

From this example, it is given the following:

$W(0) = 28,686$, and $W(1) = 32,000$.

From (2), the dependent variable is $W$, while the independent variable is the time $t$.

To solve (1), we need to use separable method as follows:

By using definition 4.2.1, we need to rewrite (1) in a way that each term is separated from the other term as follows:

$$\frac{dW}{dt} = \beta\sqrt{W(t)} = \frac{\beta}{W^{-\frac{1}{2}}} \ldots \ldots \ldots \ldots \ldots \ldots \ldots \ldots \ldots \ldots \ldots \ldots (3)$$

Now, we need to do a cross multiplication for (3) as follows:

$$\left(W^{-\frac{1}{2}}\right)dW = \beta dt$$

$$\left(W^{-\frac{1}{2}}\right)dW - \beta dt = 0 \ldots \ldots \ldots \ldots \ldots \ldots \ldots \ldots \ldots . \ldots \ldots \ldots (4)$$

Then, we integrate both sides of (4) as follows:

$$\int\left(\left(W^{-\frac{1}{2}}\right)dW - \beta dt\right) = \int 0$$

$$\int\left(W^{-\frac{1}{2}}\right)dW - \int(\beta)\,dt = c$$

$$2W^{\frac{1}{2}} - \beta t = c$$

Thus, the general solution is :

$$2W^{\frac{1}{2}} - \beta t = c \ldots \ldots \ldots \ldots \ldots \ldots \ldots \ldots \ldots \ldots \ldots \ldots (5)$$

for some $c \in \Re$.

Then, we rewrite (5) as follows:

$$2W^{\frac{1}{2}} = c + \beta t$$

$$W^{\frac{1}{2}} = \frac{c + \beta t}{2} \ldots \ldots \ldots \ldots \ldots \ldots \ldots \ldots \ldots \ldots . \ldots \ldots (6)$$

We square both sides of (6) as follows:

$$W(t) = \left(\frac{c + \beta t}{2}\right)^2 \ldots \ldots \ldots \ldots \ldots \ldots \ldots . \ldots \ldots \ldots . \ldots (7)$$

Now, we need to find $c$ by substituting $W(0) = 28{,}686$ in (7) as follows:

$$W(0) = \left(\frac{c + \beta(0)}{2}\right)^2$$

$$28{,}686 = \left(\frac{c + 0}{2}\right)^2$$

$$28{,}686 = \left(\frac{c}{2}\right)^2$$

$$28{,}686 = \frac{c^2}{4}$$

$$c^2 = 4(28{,}686)$$

$$c = \sqrt{4(28{,}686)}$$

$$c \approx 338.74$$

Thus, $W(t) = \left(\frac{338.74 + \beta t}{2}\right)^2 \ldots \ldots \ldots \ldots \ldots \ldots \ldots . \ldots \ldots \ldots \ldots . (8)$

By substituting $W(1) = 32{,}000$ in (8), we obtain:

$$W(t) = \left(\frac{338.74 + \beta}{2}\right)^2$$

$$32{,}000 = \left(\frac{(338.74)^2 + 2(338.74)\beta + \beta^2}{4}\right)$$

$$(338.74)^2 + 2(338.74)\beta + \beta^2 = 4(32{,}000)$$

$$2(338.74)\beta + \beta^2 = 4(32{,}000) - (338.74)^2$$

$$\beta^2 + 2(338.74)\beta - 13{,}255.2124 = 0$$

Thus, $\beta \approx 63{,}661.26 \ldots \ldots \ldots \ldots \ldots \ldots \ldots \ldots . \ldots \ldots \ldots . (9)$

Now, we substitute (9) in (8) as follows:

$$W(t) = \left(\frac{338.74 + (63{,}661.26)t}{2}\right)^2 \ldots \ldots . \ldots \ldots \ldots \ldots . \ldots \ldots . (10)$$

Then, we need to find the time $t$ when

$W(t) = 2(28,686) = 57,372$ by substituting it in (10) as follows:

$$57,372 = \left(\frac{338.74 + (63,661.26)t}{2}\right)^2$$

$$(338.74)^2 + 2(338.74)(63,661.26)t + t^2 = 4(57,372)$$

$$t^2 + 2(338.74)(63,661.26)t - 114,743.2124 = 0$$

$$t \approx 0.00266 \text{ years}$$

Thus, it will take approximately 0.00266 years to double the number of WSU students in 2013.

**Part b:** To determine what will be number of WSU students in 2018, we need to do the following:

We assume that $t = 2018$, and then we substitute it in (10) as follows:

$$W(2018) = \left(\frac{338.74 + (63,661.26)(2018)}{2}\right)^2$$

$$W(2018) \approx 4.126 \times 10^{15} \text{ students}$$

Thus, the number of WSU students will be approximately $4.126 \times 10^{15}$ students in 2018.

# 5.3 Water Tank Application

In this section, we give an example of water tank application, and we introduce how to use one of the differential equations methods to solve it.

**Example 5.3.1** One of the most beautiful places at *Washington State University* campus is known as *WSU Water Tower*. Assume that *WSU Water Tower* has a tank that contains initially 350 gallons of purified water, given that when $t = 0$, the amount of minerals is 5 bound. Suppose that there is a mixture of minerals containing 0.2 bound of minerals per gallon is poured into the tank at rate of 5 gallons per minute, while the mixture of minerals goes out of the tank at rate of 2 gallons per minute.

a) What is the amount of minerals in the tank of WSU Water Tower at any time $t$?

b) What is the concentration of minerals in the tank of WSU Water Tower at $t = 34$ minutes?

**Solution: Part a:** To determine the amount of minerals in the tank of WSU Water Tower at any time $t$, we need to do the following:

Assume that $W(t)$ is the amount of minerals at any time $t$, and $M(t)$ is the concentration of minerals in the tank at any time $t$. $M(t)$ is written in the following form:

$$M(t) = \frac{The\ amount\ of\ Minerals}{The\ Volume\ of\ Purified\ Water}$$

$$= \frac{W(t)}{Purified\ Water + ((Inner\ Rate - Outer\ Rate)t)} \dots (1)$$

From this example, it is given the following:

$W(0) = 5\ bounds$, $Inner\ Rate = 5$ gallons/minute, and $Outer\ Rate = 2$ gallons/minute.

Now, we need to rewrite our previous equation for this example by substituting what is given in the example itself in (1) as follows:

$$M(t) = \frac{W(t)}{350 + ((5-3)t)} = \frac{W(t)}{350 + 2t} \ldots \ldots \ldots \ldots \ldots, \ldots \ldots (2)$$

From (2), we can write the differential equation as follows:

$$\frac{dW}{dt} = 0.2 \cdot Inner\ Rate - M(t) \cdot Outer\ Rate$$

$$\frac{dW}{dt} = 0.2 \cdot (5) - \left(\frac{W(t)}{350 + 2t}\right) \cdot (2) \ldots \ldots \ldots \ldots \ldots \ldots \ldots . (3)$$

From (3), the dependent variable is $W$, while the independent variable is the time $t$. Then, we rewrite (3) as follows:

$$W'(t) = 0.2 \cdot (5) - \left(\frac{W(t)}{350 + 2t}\right) \cdot (2)$$

$$W'(t) = 1 - (2)\left(\frac{W(t)}{350 + 2t}\right)$$

$$W'(t) + \left(\frac{2}{350 + 2t}\right)W(t) = 1 \ldots \ldots \ldots \ldots \ldots \ldots . . \ldots \ldots \ldots . (4)$$

Since (4) is a first order linear differential equation, then by using definition 4.1.1, we need to use the integral factor method by letting $I = e^{\int g(t)dt}$, where $g(t) = \left(\frac{2}{350+2t}\right)$ and $F(t) = 1$.

$$I = e^{\int g(t)dt} = e^{\int \left(\frac{2}{350+2t}\right)dt} = e^{\ln(350+2t)} = (350 + 2t).$$

The general solution is written as follows:

$$W(t) = \frac{\int I \cdot F(t)\, dt}{I}$$

$$W(t) = \frac{\int (350 + 2t) \cdot (1)\, dt}{(350 + 2t)}$$

$$W(t) = \frac{\int (350 + 2t)\, dt}{(350 + 2t)}$$

$$W(t) = \frac{350t + t^2 + c}{(350 + 2t)}$$

$$W(t) = \left( \frac{350t}{(350 + 2t)} + \frac{t^2}{(350 + 2t)} + \frac{c}{(350 + 2t)} \right) \dots\dots (5)$$

The general solution is:

$W(t) = \left( \frac{350t}{(350+2t)} + \frac{t^2}{(350+2t)} + \frac{c}{(350+2t)} \right)$ for some $c \in \Re$.

Now, we need to find $c$ by substituting

$W(0) = 5$ bounds in (5) as follows:

$$W(0) = \left( \frac{350(0)}{(350 + 2(0))} + \frac{(0)^2}{(350 + 2(0))} + \frac{c}{(350 + 2(0))} \right)$$

$$5 = \left( 0 + 0 + \frac{c}{(350 + 2(0))} \right)$$

$$5 = \left( \frac{c}{350} \right)$$

$$c = (5)(350) = 1750$$

Thus, $W(t) = \left( \frac{350t}{(350+2t)} + \frac{t^2}{(350+2t)} + \frac{1750}{(350+2t)} \right). \dots\dots\dots\dots (6)$

The amount of minerals in the tank of WSU Water Tower at any time $t$ is:

$$W(t) = \left( \frac{350t}{(350 + 2t)} + \frac{t^2}{(350 + 2t)} + \frac{1750}{(350 + 2t)} \right)$$

**Part b:** To determine the concentration of minerals in the tank of WSU Water Tower at $t = 34$ minutes, we need to do the following:

We substitute $t = 34$ minutes in (6) as follows:

$$W(34) = \left( \frac{350(34) + (34)^2 + 1750}{(350 + 2(34))} \right) \approx 35.42$$

Thus, the concentration of minerals in the tank of WSU Water Tower at $t = 34$ minutes is approximately 35.42 minutes.

# Appendices

# Review of Linear Algebra

## Appendix A: Determinants*

*The materials of appendix A are taken from section 1.7 in my published book titled *A First Course in Linear Algebra: Study Guide for the Undergraduate Linear Algebra Course, First Edition*[1].

In this section, we introduce step by step for finding determinant of a certain matrix. In addition, we discuss some important properties such as invertible and non-invertible. In addition, we talk about the effect of row-operations on determinants.

**Definition A.1** Determinant is a square matrix. Given $M_2(\mathbb{R}) = \mathbb{R}^{2\times2} = \mathbb{R}_{2\times2}$, let $A \in M_2(\mathbb{R})$ where A is $2 \times 2$ matrix, $A = \begin{bmatrix} a_{11} & a_{12} \\ a_{21} & a_{22} \end{bmatrix}$. The determinant of A is represented by $\det(A)$ or $|A|$.

Hence, $\det(A) = |A| = a_{11}a_{22} - a_{12}a_{21} \in \mathbb{R}$. (Warning: this definition works only for $2 \times 2$ matrices).

**Example A.1** Given the following matrix:
$$A = \begin{bmatrix} 3 & 2 \\ 5 & 7 \end{bmatrix}$$
Find the determinant of A.

**Solution:** Using definition A.1, we do the following:

$\det(A) = |A| = (3)(7) - (2)(5) = 21 - 10 = 11$.
Thus, the determinant of A is 11.

**Example A.2** Given the following matrix:
$$A = \begin{bmatrix} 1 & 0 & 2 \\ 3 & 1 & -1 \\ 1 & 2 & 4 \end{bmatrix}$$
Find the determinant of A.

**Solution:** Since A is $3 \times 3$ matrix such that
$A \in M_3(\mathbb{R}) = \mathbb{R}^{3 \times 3}$, then we cannot use definition A.1
because it is valid only for $2 \times 2$ matrices. Thus, we
need to use the following method to find the
determinant of A.

**Step 1:** Choose any row or any column. It is
recommended to choose the one that has more zeros.
In this example, we prefer to choose the second column
or the first row. Let's choose the second column as
follows:
$$A = \begin{bmatrix} 1 & \boxed{0} & 2 \\ 3 & \boxed{1} & -1 \\ 1 & \boxed{2} & 4 \end{bmatrix}$$
$a_{12} = 0, a_{22} = 1$ and $a_{32} = 2$.

**Step 2:** To find the determinant of A, we do the
following: For $a_{12}$, since $a_{12}$ is in the first row and
second column, then we virtually remove the first row
and second column.
$$A = \begin{bmatrix} 1 & \boxed{0} & 2 \\ 3 & \boxed{1} & -1 \\ 1 & \boxed{2} & 4 \end{bmatrix}$$
$$(-1)^{1+2} a_{12} \det \begin{bmatrix} 3 & -1 \\ 1 & 4 \end{bmatrix}$$

For $a_{22}$, since $a_{22}$ is in the second row and second
column, then we virtually remove the second row and
second column.

$A = \begin{bmatrix} 1 & \boxed{0} & 2 \\ 3 & \boxed{1} & -1 \\ 1 & \boxed{2} & 4 \end{bmatrix}$

$(-1)^{2+2} a_{22} \det \begin{bmatrix} 1 & 2 \\ 1 & 4 \end{bmatrix}$

For $a_{32}$, since $a_{32}$ is in the third row and second column, then we virtually remove the third row and second column.

$A = \begin{bmatrix} 1 & \boxed{0} & 2 \\ 3 & \boxed{1} & -1 \\ 1 & \boxed{2} & 4 \end{bmatrix}$

$(-1)^{3+2} a_{32} \det \begin{bmatrix} 1 & 2 \\ 3 & -1 \end{bmatrix}$

**Step 3:** Add all of them together as follows:

$\det(A) = (-1)^{1+2} a_{12} \det \begin{bmatrix} 3 & -1 \\ 1 & 4 \end{bmatrix} + (-1)^{2+2} a_{22} \det \begin{bmatrix} 1 & 2 \\ 1 & 4 \end{bmatrix}$

$\qquad + (-1)^{3+2} a_{32} \det \begin{bmatrix} 1 & 2 \\ 3 & -1 \end{bmatrix}$

$\det(A) = (-1)^{3}(0) \det \begin{bmatrix} 3 & -1 \\ 1 & 4 \end{bmatrix} + (-1)^{4}(1) \det \begin{bmatrix} 1 & 2 \\ 1 & 4 \end{bmatrix}$

$\qquad + (-1)^{5}(2) \det \begin{bmatrix} 1 & 2 \\ 3 & -1 \end{bmatrix}$

$\det(A) = (-1)(0) \det \begin{bmatrix} 3 & -1 \\ 1 & 4 \end{bmatrix} + (1)(1) \det \begin{bmatrix} 1 & 2 \\ 1 & 4 \end{bmatrix}$

$\qquad + (-1)(2) \det \begin{bmatrix} 1 & 2 \\ 3 & -1 \end{bmatrix}$

$\det(A) = (-1)(0)(12 - -1) + (1)(1)(4 - 2) + (-1)(2)(-1 - 6)$

$\det(A) = 0 + 2 + 14 = 16.$

Thus, the determinant of A is 16.

**Result A.1** Let $A \in M_n(\mathbb{R})$. Then, A is invertible (non-singular) if and only if $\det(A) \neq 0$.

The above result means that if $\det(A) \neq 0$, then A is invertible (non-singular), and if A is invertible (non-singular), then $\det(A) \neq 0$.

**Example A.3** Given the following matrix:

$$A = \begin{bmatrix} 2 & 3 \\ 4 & 6 \end{bmatrix}$$

Is A invertible (non-singular)?

**Solution:** Using result A.1, we do the following:
$\det(A) = |A| = (2)(6) - (3)(4) = 12 - 12 = 0$.
Since the determinant of A is 0, then A is non-invertible (singular).
Thus, the answer is No because A is non-invertible (singular).

**Definition A.2** Given $A = \begin{bmatrix} a_{11} & a_{12} \\ a_{21} & a_{22} \end{bmatrix}$. Assume that $\det(A) \neq 0$ such that $\det(A) = a_{11}a_{22} - a_{12}a_{21}$. To find $A^{-1}$ (the inverse of A), we use the following format that applies only for $2 \times 2$ matrices:

$$A^{-1} = \frac{1}{\det(A)} \begin{bmatrix} a_{22} & -a_{12} \\ -a_{21} & a_{11} \end{bmatrix}$$

$$A^{-1} = \frac{1}{a_{11}a_{22} - a_{12}a_{21}} \begin{bmatrix} a_{22} & -a_{12} \\ -a_{21} & a_{11} \end{bmatrix}$$

**Example A.4** Given the following matrix:

$$A = \begin{bmatrix} 3 & 2 \\ -4 & 5 \end{bmatrix}$$

Is A invertible (non-singular)? If Yes, Find $A^{-1}$.

**Solution:** Using result A.1, we do the following:
$\det(A) = |A| = (3)(5) - (2)(-4) = 15 + 8 = 23 \neq 0$.
Since the determinant of A is not 0, then A is invertible (non-singular).

Thus, the answer is Yes, there exists $A^{-1}$ according to definition 1.7.2 as follows:

$$A^{-1} = \frac{1}{\det(A)}\begin{bmatrix} 5 & -2 \\ 4 & 3 \end{bmatrix} = \frac{1}{23}\begin{bmatrix} 5 & -2 \\ 4 & 3 \end{bmatrix} = \begin{bmatrix} \dfrac{5}{23} & -\dfrac{2}{23} \\ \dfrac{4}{23} & \dfrac{3}{23} \end{bmatrix}$$

**Result A.2** Let $A \in M_n(\mathbb{R})$ be a triangular matrix. Then, $\det(A) =$ multiplication of the numbers on the main diagonal of A.

There are three types of triangular matrix:

a) Upper Triangular Matrix: it has all zeros on the left side of the diagonal of $n \times n$ matrix.

(i.e. $A = \begin{bmatrix} 1 & 7 & 3 \\ 0 & 2 & 5 \\ 0 & 0 & 4 \end{bmatrix}$ is an Upper Triangular Matrix).

b) Diagonal Matrix: it has all zeros on both left and right sides of the diagonal of $n \times n$ matrix.

(i.e. $B = \begin{bmatrix} 1 & 0 & 0 \\ 0 & 2 & 0 \\ 0 & 0 & 4 \end{bmatrix}$ is a Diagonal Matrix).

c) Lower Triangular Matrix: it has all zeros on the right side of the diagonal of $n \times n$ matrix.

(i.e. $C = \begin{bmatrix} 1 & 0 & 0 \\ 5 & 2 & 0 \\ 1 & 9 & 4 \end{bmatrix}$ is a Diagonal Matrix).

**Fact A.1** Let $A \in M_n(\mathbb{R})$. Then, $\det(A) = \det(A^T)$.

**Fact A.2** Let $A \in M_n(\mathbb{R})$. If A is an invertible (non-singular) matrix, then $A^T$ is also an invertible (non-singular) matrix. (i.e. $(A^T)^{-1} = (A^{-1})^T$ ).

**Proof of Fact A.2** We will show that $(A^T)^{-1} = (A^{-1})^T$.

We know from previous results that $AA^{-1} = I_n$.

By taking the transpose of both sides, we obtain:

$(AA^{-1})^T = (I_n)^T$

Then, $(A^{-1})^T A^T = (I_n)^T$

Since $(I_n)^T = I_n$ , then $(A^{-1})^T A^T = I_n$.

Similarly, $(A^T)^{-1} A^T = (I_n)^T = I_n$.

Thus, $(A^T)^{-1} = (A^{-1})^T$. ∎

## The effect of Row-Operations on determinants:

Suppose $\propto$ is a non-zero constant, and *i and k* are row numbers in the augmented matrix.

* $\propto R_i$   , $\propto \neq 0$ (Multiply a row with a non-zero constant $\propto$).

i.e. $A = \begin{bmatrix} 1 & 2 & 3 \\ 0 & 4 & 1 \\ 2 & 0 & 1 \end{bmatrix} 3R_2 \dashrightarrow \begin{bmatrix} 1 & 2 & 3 \\ 0 & 12 & 3 \\ 2 & 0 & 1 \end{bmatrix} = B$

Assume that $\det(A) = \gamma$ where $\gamma$ is known, then $\det(B) = 3\gamma$.

Similarly, if $\det(B) = \beta$ where $\beta$ is known, then $\det(A) = \frac{1}{3}\beta$.

* $\propto R_i + R_k \dashrightarrow R_k$ (Multiply a row with a non-zero constant $\propto$, and add it to another row).

i.e. $A = \begin{bmatrix} 1 & 2 & 3 \\ 0 & 4 & 1 \\ 2 & 0 & 1 \end{bmatrix} \propto R_i + R_k \dashrightarrow R_k$

$\begin{bmatrix} 1 & 2 & 3 \\ 0 & 12 & 3 \\ 2 & 0 & 1 \end{bmatrix} = B$

Then, $\det(A) = \det(B)$.

* $R_i \leftrightarrow R_k$ (Interchange two rows). It has no effect on the determinants.

In general, the effect of Column-Operations on determinants is the same as for Row-Operations.

**Example A.5** Given the following $4 \times 4$ matrix A with some Row-Operations:

A $2R_1 \dashrightarrow A_1$ $3R_3 \dashrightarrow A_2$ -$2R_4 \dashrightarrow A_4$

If det(A) = 4, then find det($A_3$)

**Solution:** Using what we have learned from the effect of determinants on Row-Operations:

det($A_1$) = 2 * det(A) = 2 * 4 = 8 because $A_1$ has the first row of A multiplied by 2.

det($A_2$) = 3 * det($A_1$) = 3 * 8 = 24 because $A_2$ has the third row of $A_1$ multiplied by 3.

Similarly, det($A_3$) = $-2$ * det($A_2$) = $-2$ * 24 = $-48$ because $A_3$ has the fourth row of $A_2$ multiplied by -2.

**Result A.3** Assume A is $n \times n$ matrix with a given det(A) = $\gamma$ . Let $\alpha$ be a number. Then, det($\alpha$A) = $\alpha^n * \gamma$.

**Result A.4** Assume A and B are $n \times n$ matrices.

Then: a) det(A) = det(A) * det(B).

      b) Assume $A^{-1}$ exists and $B^{-1}$ exists.

      Then, $(AB)^{-1} = B^{-1}A^{-1}$.

      c) det(AB) = det(BA).

      d) det(A) = det($A^T$).

      e) If $A^{-1}$ exists, then det($A^{-1}$) = $\frac{1}{\det(A)}$.

**Proof of Result A.4 (b)** We will show that $(AB)^{-1} = B^{-1}A^{-1}$.

If we multiply $(B^{-1}A^{-1})$ by $(AB)$, we obtain:

$B^{-1}(A^{-1}A)B = B^{-1}(I_n)B = B^{-1}B = I_n$.

Thus, $(AB)^{-1} = B^{-1}A^{-1}$. ◻

**Proof of Result A.4 (e)** We will show that

$\det(A^{-1}) = \frac{1}{\det(A)}$.

Since $AA^{-1} = I_n$, then $\det(AA^{-1}) = \det(I_n) = 1$.

$\det(AA^{-1}) = \det(A) * \det(A^{-1}) = 1$.

Thus, $\det(A^{-1}) = \frac{1}{\det(A)}$. ◻

# Appendix B: Vector Spaces*

*The materials of appendix B are taken from chapter 2 in my published book titled *A First Course in Linear Algebra: Study Guide for the Undergraduate Linear Algebra Course, First Edition*[1].

We start this chapter reviewing some concepts of set theory, and we discuss some important concepts of vector spaces including span and dimension. In the remaining sections we introduce the concept of linear independence. At the end of this chapter we discuss other concepts such as subspace and basis.

# B.1 Span and Vector Spaces

In this section, we review some concepts of set theory, and we give an introduction to span and vector spaces including some examples related to these concepts.

Before reviewing the concepts of set theory, it is recommended to revisit section 1.4, and read the notations of numbers and the representation of the three sets of numbers in figure 1.4.1.

Let's explain some symbols and notations of set theory:

$3 \in \mathbb{Z}$ means that 3 is an element of $\mathbb{Z}$.

$\frac{1}{2} \notin \mathbb{Z}$ means that $\frac{1}{2}$ is not an element of $\mathbb{Z}$.

{ } means that it is a set.

{5} means that 5 is an element of $\mathbb{Z}$, and the set consists of exactly one element which is 5.

**Definition B.1.1** The span of a certain set is the set of all possible linear combinations of the subset of that set.

**Example B.1.1** Find Span{1}.

**Solution:** According to definition B.1.1, then the span of the set {1} is the set of all possible linear combinations of the subset of {1} which is 1.

Hence, Span{1} = $\mathbb{R}$.

**Example B.1.2** Find Span{(1,2),(2,3)}.

**Solution:** According to definition B.1.1, then the span of the set {(1,2),(2,3)} is the set of all possible linear combinations of the subsets of {(1,2),(2,3)} which are (1,2) and (2,3). Thus, the following is some possible linear combinations:

$(1,2) = 1 * (1,2) + 0 * (2,3)$

$(2,3) = 0 * (1,2) + 1 * (2,3)$

$(5,8) = 1 * (1,2) + 2 * (2,3)$

Hence, $\{(1,2), (2,3), (5,8)\} \in \text{Span}\{(1,2), (2,3)\}$.

**Example B.1.3** Find Span$\{0\}$.

**Solution:** According to definition B.1.1, then the span of the set $\{0\}$ is the set of all possible linear combinations of the subset of $\{0\}$ which is 0.

Hence, Span$\{0\} = 0$.

**Example B.1.4** Find Span$\{c\}$ where c is a non-zero integer.

**Solution:** Using definition B.1.1, the span of the set $\{c\}$ is the set of all possible linear combinations of the subset of $\{c\}$ which is $c \neq 0$.

Thus, Span$\{c\} = \mathbb{R}$.

**Definition B.1.2** $\mathbb{R}^n = \{(a_1, a_2, a_3, \ldots, a_n) | a_1, a_2, a_3, \ldots, a_n \in \mathbb{R}\}$ is a set of all points where each point has exactly $n$ coordinates.

**Definition B.1.3** $(V, +, \cdot)$ is a vector space if satisfies the following:

    a. For every $v_1, v_2 \in V$, $v_1 + v_2 \in V$.

    b. For every $\alpha \in \mathbb{R}$ and $v \in V$, $\alpha v \in V$.

(i.e. Given $Span\{x, y\}$ and $set\ \{x, y\}$, then $\sqrt{10}x + 2y \in Span\{x, y\}$. Let's assume that $v \in Span\{x, y\}$, then $v = c_1 x + c_2 y$ for some numbers $c_1$ and $c_2$).

# B.2 The Dimension of Vector Space

In this section, we discuss how to find the dimension of vector space, and how it is related to what we have learned in section B.1.

**Definition B.2.1** Given a vector space $V$, the dimension of $V$ is the number of minimum elements needed in $V$ so that their *Span* is equal to $V$, and it is denoted by $\dim(V)$. (i.e. $\dim(\mathbb{R}) = 1$ and $\dim(\mathbb{R}^2) = 2$).

**Result B.2.1** $\dim(\mathbb{R}^n) = n$.

**Proof of Result B.2.1** We will show that $\dim(\mathbb{R}^n) = n$.

Claim: $D = Span\{(1,0),(0.1)\} = \mathbb{R}^2$

$\alpha_1(1,0) + \alpha_2(0,1) = (\alpha_1, \alpha_2) \in \mathbb{R}^2$

Thus, $D$ is a subset of $\mathbb{R}^2$ ($D \subseteq \mathbb{R}^2$).

For every $x_1, y_1 \in \mathbb{R}$, $(x_1, y_1) \in \mathbb{R}^2$.

Therefore, $(x_1, y_1) = x_1(1,0) + y_1(0,1) \in D$.

We prove the above claim, and hence

$\dim(\mathbb{R}^n) = n$.                    $\square$

**Fact 2B.2.1** $Span\{(3,4)\} \neq \mathbb{R}^2$.

**Proof of Fact B.2.1** We will show that $Span\{(3,4)\} \neq \mathbb{R}^2$.

Claim: $F = Span\{(6,5)\} \neq \mathbb{R}^2$ where $(6,5) \in \mathbb{R}^2$.

We cannot find a number $\alpha$ such that $(6,5) = \alpha(3,4)$

We prove the above claim, and hence $Span\{(3,4)\} \neq \mathbb{R}^2$. ∎

**Fact B.2.2** $Span\{(1,0),(0,1)\} = \mathbb{R}^2$.

**Fact B.2.3** $Span\{(2,1),(1,0.5)\} \neq \mathbb{R}^2$.

# B.3 Linear Independence

In this section, we learn how to determine whether vector spaces are linearly independent or not.

**Definition B.3.1** Given a vector space $(V,+,\cdot)$, we say $v_1, v_2, \ldots, v_n \in V$ are linearly independent if none of them is a linear combination of the remaining $v_i's$.

(i.e. $(3,4),(2,0) \in \mathbb{R}$ are linearly independent because we cannot write them as a linear combination of each other, in other words, we cannot find a number $\alpha_1, \alpha_2$ such that $(3,4) = \alpha_1(2,0)$ and $(2,0) = \alpha_2(3,4)$).

**Definition B.3.2** Given a vector space $(V,+,\cdot)$, we say $v_1, v_2, \ldots, v_n \in V$ are linearly dependent if at least one of $v_i's$ is a linear combination of the others.

**Example B.3.1** Assume $v_1$ and $v_2$ are linearly independent. Show that $v_1$ and $3v_1 + v_2$ are linearly independent.

**Solution:** We will show that $v_1$ and $3v_1 + v_2$ are linearly independent. Using proof by contradiction, we assume that $v_1$ and $3v_1 + v_2$ are linearly dependent. For some non-zero number $c_1$, $v_1 = c_1(3v_1 + v_2)$.

Using the distribution property and algebra, we obtain:

$$v_1 = 3v_1 c_1 + v_2 c_1$$

$$v_1 - 3v_1 c_1 = v_2 c_1$$

$$v_1(1 - 3c_1) = v_2 c_1$$

$$\frac{(1 - 3c_1)}{c_1} v_1 = v_2$$

Thus, none of $v_1$ and $3v_1 + v_2$ is a linear combination of the others which means that $v_1$ and $3v_1 + v_2$ are linearly independent. This is a contradiction.

Therefore, our assumption that $v_1$ and $3v_1 + v_2$ were linearly dependent is false. Hence, $v_1$ and $3v_1 + v_2$ are linearly independent.                                ◻

**Example B.3.2** Given the following vectors:

$$v_1 = (1,0,-2)$$

$$v_2 = (-2,2,1)$$

$$v_3 = (-1,0,5)$$

Are these vectors independent elements?

**Solution:** First of all, to determine whether these vectors are independent elements or not, we need to write these vectors as a matrix.

$\begin{bmatrix} 1 & 0 & -2 \\ -2 & 2 & 1 \\ -1 & 0 & 5 \end{bmatrix}$ Each point is a row-operation. We need to reduce this matrix to Semi-Reduced Matrix.

**Definition B.3.3** Semi-Reduced Matrix is a reduced-matrix but the leader numbers can be any non-zero number.

Now, we apply the Row-Reduction Method to get the Semi-Reduced Matrix as follows:

$\begin{bmatrix} 1 & 0 & -2 \\ -2 & 2 & 1 \\ -1 & 0 & 5 \end{bmatrix} \begin{matrix} 2R_1 + R_2 \to R_2 \\ R_1 + R_3 \to R_3 \end{matrix} \begin{bmatrix} 1 & 0 & -2 \\ 0 & \boxed{2} & -3 \\ 0 & 0 & 3 \end{bmatrix}$ This is a Semi-Reduced Matrix.

Since none of the rows in the Semi-Reduced Matrix become zero-row, then the elements are independent because we cannot write at least one of them as a linear combination of the others.

**Example 2.3.3** Given the following vectors:

$$v_1 = (1, -2, 4, 6)$$
$$v_2 = (-1, 2, 0, 2)$$
$$v_3 = (1, -2, 8, 14)$$

Are these vectors independent elements?

**Solution:** First of all, to determine whether these vectors are independent elements or not, we need to write these vectors as a matrix.

$\begin{bmatrix} 1 & -2 & 4 & 6 \\ -1 & 2 & 0 & 2 \\ 1 & -2 & 8 & 14 \end{bmatrix}$ Each point is a row-operation. We need to reduce this matrix to Semi-Reduced Matrix.

Now, we apply the Row-Reduction Method to get the Semi-Reduced Matrix as follows:

$\begin{bmatrix} 1 & -2 & 4 & 6 \\ -1 & 2 & 0 & 2 \\ 1 & -2 & 8 & 14 \end{bmatrix} \begin{matrix} R_1 + R_2 \to R_2 \\ -R_1 + R_3 \to R_3 \end{matrix} \begin{bmatrix} 1 & -2 & 4 & 6 \\ 0 & 0 & \boxed{4} & 8 \\ 0 & 0 & 4 & 8 \end{bmatrix}$

$-R_2 + R_3 \rightarrow R_3$ $\begin{bmatrix} 1 & -2 & 4 & 6 \\ 0 & 0 & \boxed{4} & 8 \\ \boxed{0 \quad 0 \quad 0 \quad 0} \end{bmatrix}$ This is a Semi-Reduced Matrix.

Since there is a zero-row in the Semi-Reduced Matrix, then the elements are dependent because we can write at least one of them as a linear combination of the others.

# B.4 Subspace and Basis

In this section, we discuss one of the most important concepts in linear algebra that is known as subspace. In addition, we give some examples explaining how to find the basis for subspace.

**Definition B.4.1** Subspace is a vector space but we call it a subspace because it lives inside a bigger vector space. (i.e. Given vector spaces $V$ and $D$, then according to the figure 2.4.1, $D$ is called a subspace of $V$).

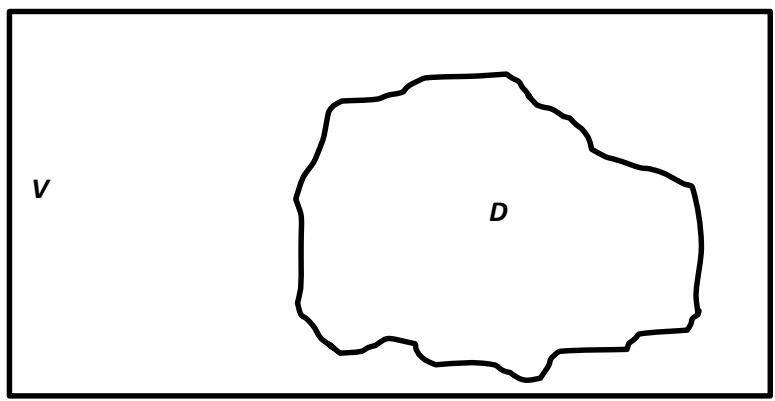

Figure B.4.1: Subspace of $V$

**Fact B.4.1** Every vector space is a subspace of itself.

**Example B.4.1** Given a vector space $L = \{(c, 3c)|c \in \mathbb{R}\}$.

    a.  Does $L$ live in $\mathbb{R}^2$?

    b.  Does $L$ equal to $\mathbb{R}^2$?

    c.  Is $L$ a subspace of $\mathbb{R}^2$?

    d.  Does $L$ equal to $Span\{(0,3)\}$?

    e.  Does $L$ equal to $Span\{(1,3),(2,6)\}$?

**Solution:** To answer all these questions, we need first to draw an equation from this vector space, say $y = 3x$. The following figure represents the graph of the above equation, and it passes through a point (1,3).

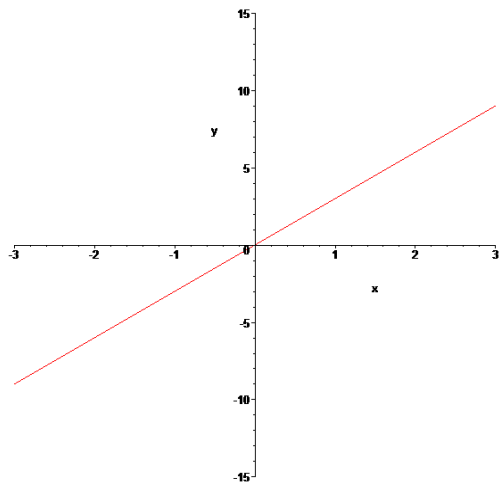

Figure B.4.2: Graph of $y = 3x$

Now, we can answer the given questions as follows:

**Part a:** Yes; $L$ lives in $\mathbb{R}^2$.

**Part b:** No; $L$ does not equal to $\mathbb{R}^2$. To show that we prove the following claim:

Claim: $L = Span\{(5,15)\} \neq \mathbb{R}^2$ where $(5,15) \in \mathbb{R}^2$.

It is impossible to find a number $\alpha = 3$ such that

$$(20,60) = \alpha(5,15)$$

because in this case $\alpha = 4$ where $(20,60) = 4(5,15)$.

We prove the above claim, and $Span\{(5,15)\} \neq \mathbb{R}^2$.

Thus, $L$ does not equal to $\mathbb{R}^2$

**Part c:** Yes; $L$ is a subspace of $\mathbb{R}^2$ because $L$ lives inside a bigger vector space which is $\mathbb{R}^2$.

**Part d:** No; according to the graph in figure 2.4.2, $(0,3)$ does not belong to $L$.

**Part e:** Yes; because we can write $(1,3)$ and $(2,6)$ as a linear combination of each other.

$\alpha_1(1,3) + \alpha_2(2,6) = \{(\alpha_1 + 2\alpha_2), (3\alpha_1 + 6\alpha_2)\}$

$\alpha_1(1,3) + \alpha_2(2,6) = \{(\alpha_1 + 2\alpha_2), 3(\alpha_1 + 2\alpha_2)\}$

Assume $c = (\alpha_1 + 2\alpha_2)$, then we obtain:

$\alpha_1(1,3) + \alpha_2(2,6) = \{(c, 3c)|c \in \mathbb{R}\} = L$.

Thus, $L = Span\{(1,3), (2,6)\}$.

**Result B.4.1** $L$ is a subspace of $\mathbb{R}^2$ if satisfies the following:

    a. $L$ lives inside $\mathbb{R}^2$.
    b. $L$ has only lines through the origin $(0,0)$.

**Example B.4.2** Given a vector space

$D = \{(a, b, 1) | a, b \in \mathbb{R}\}$.

    a.  Does $D$ live in $\mathbb{R}^3$?

    b.  Is $D$ a subspace of $\mathbb{R}^3$?

**Solution:** Since the equation of the above vector space is a three-dimensional equation, there is no need to draw it because it is difficult to draw it exactly. Thus, we can answer the above questions immediately.

**Part a:** Yes; $D$ lives inside $\mathbb{R}^3$.

**Part b:** No; since $(0,0,0) \notin D$, then $D$ is not a subspace of $\mathbb{R}^3$.

**Fact B.4.2** Assume $D$ lives inside $\mathbb{R}^n$. If we can write $D$ as a *Span*, then it is a subspace of $\mathbb{R}^n$.

**Fact B.4.3** Assume $D$ lives inside $\mathbb{R}^n$. If we cannot write $D$ as a *Span*, then it is not a subspace of $\mathbb{R}^n$.

**Fact B.4.4** Assume $D$ lives inside $\mathbb{R}^n$. If $(0,0,0,...,0)$ is in $D$, then $D$ is a subspace of $\mathbb{R}^n$.

**Fact B.4.5** Assume $D$ lives inside $\mathbb{R}^n$. If $(0,0,0,...,0)$ is not in $D$, then $D$ is not a subspace of $\mathbb{R}^n$.

Now, we list the main results on $\mathbb{R}^n$:

**Result B.4.2** Maximum number of independent points is $n$.

**Result B.4.3** Choosing any $n$ independent points in $\mathbb{R}^n$, say $Q_1, Q_2, ..., Q_n$, then $\mathbb{R}^n = Span\{Q_1, Q_2, ..., Q_n\}$.

**Result B.4.4** $\dim(\mathbb{R}^n) = n$.

Results B.4.3 and B.4.4 tell us the following: In order to get all $\mathbb{R}^n$, we need exactly $n$ independent points.

**Result B.4.5** Assume $\mathbb{R}^n = Span\{Q_1, Q_2, \dots, Q_k\}$, then $k \geq n$ ($n$ points of the $Q_k's$ are independents).

**Definition B.4.2** Basis is the set of points that is needed to $Span$ the vector space.

**Example B.4.3** Let $D = Span\{(1, -1,0), (2,2,1), (0,4,1)\}$.

   a.  Find dim($D$).

   b.  Find a basis for $D$.

**Solution:** First of all, we have infinite set of points, and $D$ lives inside $\mathbb{R}^3$. Let's assume the following:

$$v_1 = (1, -1,0)$$
$$v_2 = (2,2,1)$$
$$v_3 = (0,4,1)$$

**Part a:** To find dim($D$), we check whether $v_1, v_2$ and $v_3$ are dependent elements or not. Using what we have learned so far from section 2.3: We need to write these vectors as a matrix.

$\begin{bmatrix} 1 & -1 & 0 \\ 2 & 2 & 1 \\ 0 & 4 & 1 \end{bmatrix}$ Each point is a row-operation. We need to reduce this matrix to Semi-Reduced Matrix.

Now, we apply the Row-Reduction Method to get the Semi-Reduced Matrix as follows:

$\begin{bmatrix} 1 & -1 & 0 \\ 2 & 2 & 1 \\ 0 & 4 & 1 \end{bmatrix} -2R_1 + R_2 \to R_2 \begin{bmatrix} 1 & -1 & 0 \\ 0 & \boxed{4} & 1 \\ 0 & 4 & 1 \end{bmatrix} -R_2 + R_3 \to R_3$

$$\begin{bmatrix} 1 & -1 & 0 \\ 0 & \boxed{4} & 1 \\ \boxed{0 & 0 & 0} \end{bmatrix}$$ This is a Semi-Reduced Matrix.

Since there is a zero-row in the Semi-Reduced Matrix, then these elements are dependent because we can write at least one of them as a linear combination of the others. Only two points survived in the Semi-Reduced Matrix. Thus, $\dim(D) = 2$.

**Part b:** $D$ is a plane that passes through the origin $(0,0,0)$. Since $\dim(D) = 2$, then any two independent points in $D$ will form a basis for $D$. Hence, the following are some possible bases for $D$:

Basis for $D$ is $\{(1,-1,0),(2,2,1)\}$.

Another basis for $D$ is $\{(1,-1,0),(0,4,1)\}$.

**Result B.4.6** It is always true that $|Basis| = dim(D)$.

**Example B.4.4** Given the following:

$M = Span\{(-1,2,0,0),(1,-2,3,0),(-2,0,3,0)\}$.

Find a basis for $M$.

**Solution:** We have infinite set of points, and $M$ lives inside $\mathbb{R}^4$. Let's assume the following:

$$v_1 = (-1,2,0,0)$$
$$v_2 = (1,-2,3,0)$$
$$v_3 = (-2,0,3,0)$$

We check if $v_1, v_2$ and $v_3$ are dependent elements. Using what we have learned so far from section 2.3 and example 2.4.3: We need to write these vectors as a matrix.

$$\begin{bmatrix} -1 & 2 & 0 & 0 \\ 1 & -2 & 3 & 0 \\ -2 & 0 & 3 & 0 \end{bmatrix}$$ Each point is a row-operation. We need to reduce this matrix to Semi-Reduced Matrix.

Now, we apply the Row-Reduction Method to get the Semi-Reduced Matrix as follows:

$$\begin{bmatrix} -1 & 2 & 0 & 0 \\ 1 & -2 & 3 & 0 \\ -2 & 0 & 3 & 0 \end{bmatrix} \begin{matrix} R_1 + R_2 \to R_2 \\ -2R_1 + R_3 \to R_3 \end{matrix} \begin{bmatrix} -1 & 2 & 0 & 0 \\ 0 & 0 & \boxed{3} & 0 \\ 0 & -4 & 3 & 0 \end{bmatrix}$$

$$-R_2 + R_3 \to R_3 \begin{bmatrix} -1 & 2 & 0 & 0 \\ 0 & 0 & 3 & 0 \\ 0 & -4 & 0 & 0 \end{bmatrix}$$ This is a Semi-Reduced Matrix.

Since there is no zero-row in the Semi-Reduced Matrix, then these elements are independent. All the three points survived in the Semi-Reduced Matrix. Thus, $\dim(M) = 3$. Since $\dim(M) = 3$, then any three independent points in $M$ from the above matrices will form a basis for $M$. Hence, the following are some possible bases for $M$:

Basis for $M$ is $\{(-1,2,0,0), (0,0,3,0), (0,-4,0,0)\}$.

Another basis for $M$ is $\{(-1,2,0,0), (0,0,3,0), (0,-4,3,0)\}$.

Another basis for $M$ is $\{(-1,2,0,0), (1,-2,3,0), (-2,0,3,0)\}$.

**Example B.4.5** Given the following:

$W = Span\{(a, -2a + b, -a) | a, b \in \mathbb{R}\}$.

   a. Show that $W$ is a subspace of $\mathbb{R}^3$.

   b. Find a basis for $W$.

   c. Rewrite $W$ as a $Span$.

**Solution:** We have infinite set of points, and $W$ lives inside $\mathbb{R}^3$.

**Part a:** We write each coordinate of $W$ as a linear combination of the free variables $a$ and $b$.

$$a = 1 \cdot a + 0 \cdot b$$

$$-2a + b = -2 \cdot a + 1 \cdot b$$

$$-a = -1 \cdot a + 0 \cdot b$$

Since it is possible to write each coordinate of $W$ as a linear combination of the free variables $a$ and $b$, then we conclude that $W$ is a subspace of $\mathbb{R}^3$.

**Part b:** To find a basis for $W$, we first need to find $\dim(W)$. To find $\dim(W)$, let's play a game called (ON-OFF GAME) with the free variables $a$ and $b$.

| $a$ | $b$ | $Point$ |
|-----|-----|---------|
| 1 | 0 | $(1, -2, -1)$ |
| 0 | 1 | $(0,1,0)$ |

Now, we check for independency: We already have the Semi-Reduced Matrix: $\begin{bmatrix} 1 & -2 & -1 \\ 0 & 1 & 0 \end{bmatrix}$. Thus, $\dim(W) = 2$.

Hence, the basis for $W$ is $\{(1, -2, -1), (0,1,0)\}$.

**Part b:** Since we found the basis for $W$, then it is easy to rewrite $W$ as a $Span$ as follows:

$W = Span\{(1, -2, -1), (0,1,0)\}$.

**Fact B.4.6** $\dim(W) \leq Number of\ Free - Variables$.

**Example B.4.6** Given the following:

$H = Span\{(a^2, 3b + a, -2c, a + b + c) | a, b, c \in \mathbb{R}\}.$

Is $H$ a subspace of $\mathbb{R}^4$?

**Solution:** We have infinite set of points, and $H$ lives inside $\mathbb{R}^4$. We try write each coordinate of $H$ as a linear combination of the free variables $a, b$ and $c$.

$a^2 = Fixed\ Number \cdot a + Fixed\ Number \cdot b + Fixed\ Number \cdot c$

$a^2$ is not a linear combination of $a, b$ and $c$.

We assume that $w = (1,1,0,1) \in H$, and $a = 1, b = c = 0$.

If $\alpha = -2$, then $-2 \cdot w = -2 \cdot (1,1,0,1) = (-2, -2, 0, -2) \notin H$.

Since it is impossible to write each coordinate of $H$ as a linear combination of the free variables $a, b$ and $c$, then we conclude that $H$ is not a subspace of $\mathbb{R}^4$.

**Example B.4.7** Form a basis for $\mathbb{R}^4$.

**Solution:** We just need to select any random four independent points, and then we form a $4 \times 4$ matrix with four independent rows as follows:

$\begin{bmatrix} 2 & 3 & 0 & 4 \\ 0 & 5 & 1 & 1 \\ 0 & 0 & 2 & 3 \\ 0 & 0 & 0 & \pi^e \end{bmatrix}$ Note: $\pi^e$ is a number.

Let's assume the following:

$$v_1 = (2,3,0,4)$$
$$v_2 = (0,5,1,1)$$
$$v_3 = (0,0,2,3)$$
$$v_4 = (0,0,0,\pi^e)$$

Thus, the basis for $\mathbb{R}^4 = \{v_1, v_2, v_3, v_4\}$, and

$Span\{v_1, v_2, v_3, v_4\} = \mathbb{R}^4$.

**Example B.4.8** Form a basis for $\mathbb{R}^4$ that contains the following two independent points:

$(0,2,1,4)$ and $(0,-2,3,-10)$.

**Solution:** We need to add two more points to the given one so that all four points are independent. Let's assume the following:

$$v_1 = (0,2,1,4)$$
$$v_2 = (0,-2,3,-10)$$
$$v_3 = (0,0,4,-6) \text{ This is a random point.}$$
$$v_4 = (0,0,0,1000) \text{ This is a random point.}$$

Then, we need to write these vectors as a matrix.

$\begin{bmatrix} 0 & 2 & 1 & 4 \\ 0 & -2 & 3 & -10 \\ 0 & 0 & 4 & -6 \\ 0 & 0 & 0 & 1000 \end{bmatrix}$ Each point is a row-operation. We

need to reduce this matrix to Semi-Reduced Matrix.

Now, we apply the Row-Reduction Method to get the Semi-Reduced Matrix as follows:

$\begin{bmatrix} 0 & 2 & 1 & 4 \\ 0 & -2 & 3 & -10 \\ 0 & 0 & 4 & -6 \\ 0 & 0 & 0 & 1000 \end{bmatrix}$ $R_1 + R_2 \rightarrow R_2$ $\begin{bmatrix} 0 & 2 & 1 & 4 \\ 3 & 0 & 5 & 30 \\ 0 & 0 & 4 & -6 \\ 0 & 0 & 0 & 1000 \end{bmatrix}$

This is a Semi-Reduced Matrix.

Thus, the basis for $\mathbb{R}^4$ is

$\{(0,2,1,4), (0,-2,3,-10), (3,0,5,30), (0,0,0,1000)\}$.

**Example B.4.9** Given the following:

$D = Span\{(1,1,1,1), (-1,-1,0,0), (0,0,1,1)\}$

Is $(1,1,2,2) \in D$?

**Solution:** We have infinite set of points, and $D$ lives inside $\mathbb{R}^4$. There are two different to solve this example:

**The First Way:** Let's assume the following:

$$v_1 = (1,1,1,1)$$
$$v_2 = (-1,-1,0,0)$$
$$v_3 = (0,0,1,1)$$

We start asking ourselves the following question:

Question: Can we find $\alpha_1, \alpha_2$ and $\alpha_3$ such that $(1,1,2,2) = \alpha_1 \cdot v_1 + \alpha_2 \cdot v_2 + \alpha_3 \cdot v_3$?

Answer: Yes but we need to solve the following system of linear equations:

$$1 = \alpha_1 - \alpha_2 + 0 \cdot \alpha_3$$
$$1 = \alpha_1 - \alpha_2 + 0 \cdot \alpha_3$$
$$2 = \alpha_1 + \alpha_3$$
$$2 = \alpha_1 + \alpha_3$$

Using what we have learned from chapter 1 to solve the above system of linear equations, we obtain:

$$\alpha_1 = \alpha_2 = \alpha_3 = 1$$

Hence, Yes: $(1,1,2,2) \in D$.

**The Second Way (Recommended):** We first need to find $dim(D)$, and then a basis for $D$. We have to write $v_1, v_2$ and $v_3$ as a matrix.

$$\begin{bmatrix} 1 & 1 & 1 & 1 \\ -1 & -1 & 0 & 0 \\ 0 & 0 & 1 & 1 \end{bmatrix}$$ Each point is a row-operation. We need to reduce this matrix to Semi-Reduced Matrix.

Now, we apply the Row-Reduction Method to get the Semi-Reduced Matrix as follows:

$$\begin{bmatrix} 1 & 1 & 1 & 1 \\ -1 & -1 & 0 & 0 \\ 0 & 0 & 1 & 1 \end{bmatrix} R_1 + R_2 \to R_2 \begin{bmatrix} 1 & 1 & 1 & 1 \\ 0 & 0 & 1 & 1 \\ 0 & 0 & 1 & 1 \end{bmatrix}$$

$$-R_2 + R_3 \to R_3 \begin{bmatrix} 1 & 1 & 1 & 1 \\ 0 & 0 & 1 & 1 \\ 0 & 0 & 0 & 0 \end{bmatrix}$$ This is a Semi-Reduced Matrix.

Since there is a zero-row in the Semi-Reduced Matrix, then these elements are dependent. Thus, $\dim(D) = 2$.

Thus, Basis for $D$ is $\{(1,1,1,1),(0,0,1,1)\}$, and

$D = Span\{(1,1,1,1),(0,0,1,1)\}$.

Now, we ask ourselves the following question:

Question: Can we find $\alpha_1, \alpha_2$ and $\alpha_3$ such that $(1,1,2,2) = \alpha_1 \cdot (1,1,1,1) + \alpha_2 \cdot (0,0,1,1)$?

Answer: Yes:

$$1 = \alpha_1$$

$$1 = \alpha_1$$

$$2 = \alpha_1 + \alpha_2$$

$$2 = \alpha_1 + \alpha_2$$

Thus, $\alpha_1 = \alpha_2 = \alpha_3 = 1$. Hence, Yes: $(1,1,2,2) \in D$.

# Appendix C: Homogenous Systems*

*The materials of appendix C are taken from chapter 3 in my published book titled *A First Course in Linear Algebra: Study Guide for the Undergraduate Linear Algebra Course, First Edition*[1].

In this chapter, we introduce the homogeneous systems, and we discuss how they are related to what we have learned in chapter B. We start with an introduction to null space and rank. Then, we study one of the most important topics in linear algebra which is linear transformation. At the end of this chapter we discuss how to find range and kernel, and their relation to sections C.1 and C.2.

# C.1 Null Space and Rank

In this section, we first give an introduction to homogeneous systems, and we discuss how to find the null space and rank of homogeneous systems. In addition, we explain how to find row space and column space.

**Definition C.1.1** Homogeneous System is a $m \times n$ system of linear equations that has all zero constants.

(i.e. the following is an example of homogeneous system): $\begin{cases} 2x_1 + x_2 - x_3 + x_4 = 0 \\ 3x_1 + 5x_2 + 3x_3 + 4x_4 = 0 \\ -x_2 + x_3 - x_4 = 0 \end{cases}$

Imagine we have the following solution to the homogeneous system: $x_1 = x_2 = x_3 = x_4 = 0$.

Then, this solution can be viewed as a point of $\mathbb{R}^n$ (here is $\mathbb{R}^4$) : (0,0,0,0)

**Result C.1.1** The solution of a homogeneous system $m \times n$ can be written as

$\{(a_1, a_2, a_3, a_4, \dots, a_n | a_1, a_2, a_3, a_4, \dots, a_n \in \mathbb{R}\}$.

**Result C.1.2** All solutions of a homogeneous system $m \times n$ form a subset of $\mathbb{R}^n$, and it is equal to the number of variables.

**Result C.1.3** Given a homogeneous system $m \times n$. We write it in the matrix-form: $C \begin{bmatrix} x_1 \\ x_2 \\ x_3 \\ \vdots \\ x_n \end{bmatrix} = \begin{bmatrix} 0 \\ 0 \\ 0 \\ \vdots \\ 0 \end{bmatrix}$ where $C$ is a coefficient. Then, the set of all solutions in this system is a subspace of $\mathbb{R}^n$.

**Proof of Result C.1.3** We assume that $M_1 = (m_1, m_2, \dots, m_n)$ and $W_1 = (w, w_2, \dots, w_n)$ are two solutions to the above system. We will show that $M +$

$W$ is a solution. We write them in the matrix-form:

$$C\begin{bmatrix} m_1 \\ m_2 \\ m_3 \\ \vdots \\ m_n \end{bmatrix} = \begin{bmatrix} 0 \\ 0 \\ 0 \\ \vdots \\ 0 \end{bmatrix} \text{ and } \begin{bmatrix} w_1 \\ w_2 \\ w_3 \\ \vdots \\ w_n \end{bmatrix} = \begin{bmatrix} 0 \\ 0 \\ 0 \\ \vdots \\ 0 \end{bmatrix}$$

Now, using algebra: $M + W = C\begin{bmatrix} m_1 \\ m_2 \\ m_3 \\ \vdots \\ m_n \end{bmatrix} + C\begin{bmatrix} w_1 \\ w_2 \\ w_3 \\ \vdots \\ w_n \end{bmatrix} = \begin{bmatrix} 0 \\ 0 \\ 0 \\ \vdots \\ 0 \end{bmatrix}$

By taking $C$ as a common factor, we obtain:

$$C\left(\begin{bmatrix} m_1 \\ m_2 \\ m_3 \\ \vdots \\ m_n \end{bmatrix} + \begin{bmatrix} w_1 \\ w_2 \\ w_3 \\ \vdots \\ w_n \end{bmatrix}\right) = \begin{bmatrix} 0 \\ 0 \\ 0 \\ \vdots \\ 0 \end{bmatrix}$$

$$C\begin{bmatrix} m_1 + w_1 \\ m_2 + w_2 \\ m_3 + w_3 \\ \vdots \\ m_n + w_n \end{bmatrix} = \begin{bmatrix} 0 \\ 0 \\ 0 \\ \vdots \\ 0 \end{bmatrix}$$

Thus, $M + W$ is a solution.  □

**Fact C.1.1** If $M_1 = (m_1, m_2, \ldots, m_n)$ is a solution, and $\alpha \in \mathbb{R}$, then $\alpha M = (\alpha m_1, \alpha m_2, \ldots, \alpha m_n)$ is a solution.

**Fact C.1.2** The only system where the solutions form a vector space is the homogeneous system.

**Definition C.1.2** Null Space of a matrix, say $A$ is a set of all solutions to the homogeneous system, and it is denoted by $Null(A)$ or $N(A)$.

**Definition C.1.3** Rank of a matrix, say $A$ is the number of independent rows or columns of $A$, and it is denoted by $Rank(A)$.

**Definition C.1.4** Row Space of a matrix, say $A$ is the $Span$ of independent rows of $A$, and it is denoted by $Row(A)$.

**Definition C.1.5** Column Space of a matrix, say $A$ is the $Span$ of independent columns of $A$, and it is denoted by $Column(A)$.

**Example C.1.1** Given the following $3 \times 5$ matrix:

$$A = \begin{bmatrix} 1 & -1 & 2 & 0 & -1 \\ 0 & 1 & 2 & 0 & 2 \\ 0 & 0 & 0 & 1 & 0 \end{bmatrix}.$$

   a.  Find $Null(A)$.

   b.  Find $dim(Null(A))$.

   c.  Rewrite $Null(A)$ as $Span$.

   d.  Find $Rank(A)$.

   e.  Find $Row(A)$.

**Solution: Part a:** To find the null space of $A$, we need to find the solution of $A$ as follows:

**Step 1:** Write the above matrix as an Augmented-Matrix, and make all constants' terms zeros.

$$\begin{pmatrix} 1 & -1 & 2 & 0 & -1 & | & 0 \\ 0 & 1 & 2 & 0 & 2 & | & 0 \\ 0 & 0 & 0 & 1 & 0 & | & 0 \end{pmatrix}$$

**Step 2:** Apply what we have learned from chapter 1 to solve systems of linear equations use Row-Operation Method.

$$\begin{pmatrix} 1 & -1 & 2 & 0 & -1 & | & 0 \\ 0 & 1 & 2 & 0 & 2 & | & 0 \\ 0 & 0 & 0 & 1 & 0 & | & 0 \end{pmatrix} R_2 + R_1 \rightarrow R_1$$

$$\begin{pmatrix} 1 & 0 & 4 & 0 & 1 & | & 0 \\ 0 & 1 & 2 & 0 & 2 & | & 0 \\ 0 & 0 & 0 & 1 & 0 & | & 0 \end{pmatrix} \text{This is a Completely-Reduced}$$

Matrix.

**Step 3:** Read the solution for the above system of linear equations after using Row-Operation.

$$x_1 + 4x_3 + x_5 = 0$$
$$x_2 + 2x_3 + 2x_5 = 0$$
$$x_4 = 0$$

Free variables are $x_3$ and $x_5$.

Assuming that $x_3, x_5 \in \mathbb{R}$. Then, the solution of the above homogeneous system is as follows:

$$x_1 = -4x_3 - x_5$$
$$x_2 = -2x_3 - 2x_5$$
$$x_4 = 0$$

Thus, according to definition 3.1.2,

$Null(A) = \{(-4x_3 - x_5, -2x_3 - 2x_5, x_3, 0, x_5) | x_3, x_5 \in \mathbb{R}\}.$

**Part b:** It is always true that

$dim\big(Null(A)\big) = dim\big(N(A)\big) = The\ Number\ of\ Free\ Variables$

Here, $dim\big(Null(A)\big) = 2$.

**Definition C.1.6** The nullity of a matrix, say $A$ is the dimension of the null space of $A$, and it is denoted by $dim(Null(A))$ or $dim(N(A))$.

**Part c:** We first need to find a basis for $Null(A)$ as follows: To find a basis for $Null(A)$, we play a game called (ON-OFF GAME) with the free variables $x_3$ and $x_5$.

| $x_3$ | $x_5$ | $Point$ |
|-------|-------|---------|
| 1 | 0 | $(-4,-2,1,0,0)$ |
| 0 | 1 | $(-1,-2,0,0,1)$ |

The basis for $Null(A) = \{(-4,-2,1,0,0), (-1,-2,0,0,1)\}$.

Thus, $Null(A) = Span\{(-4,-2,1,0,0), (-1,-2,0,0,1)\}$.

**Part d:** To find the rank of matrix $A$, we just need to change matrix $A$ to the Semi-Reduced Matrix. We already did that in part a. Thus, $Rank(A) = 3$.

**Part e:** To find the row space of matrix $A$, we just need to write the $Span$ of independent rows. Thus, $Row(A) = Span\{(1,-1,2,0,-1), (0,1,2,0,2), (0,0,0,1,0)\}$. It is also a subspace of $\mathbb{R}^5$.

**Result C.1.4** Let $A$ be $m \times n$ matrix. Then, $Rank(A) + dim(N(A)) = n = Number\ of\ Columns\ of\ A$.

**Result C.1.5** Let $A$ be $m \times n$ matrix. The geometric meaning of $Row(A) = Span\{Independent\ Rows\}$ "lives" inside $\mathbb{R}^n$.

**Result C.1.6** Let $A$ be $m \times n$ matrix. The geometric meaning of $Column(A) = Span\{Independent\ Columns\}$ "lives" inside $\mathbb{R}^m$.

**Result C.1.7** Let $A$ be $m \times n$ matrix. Then,

$Rank(A) = dim\big(Row(A)\big) = dim(Column(A))$.

**Example C.1.2** Given the following $3 \times 5$ matrix:

$$B = \begin{bmatrix} 1 & 1 & 1 & 1 & 1 \\ -1 & -1 & -1 & 0 & 2 \\ 0 & 0 & 0 & 0 & 0 \end{bmatrix}.$$

   a. Find $Row(B)$.

   b. Find $Column(B)$.

   c. Find $Rank(B)$.

**Solution: Part a:** To find the row space of $B$, we need to change matrix $B$ to the Semi-Reduced Matrix as follows:

$$\begin{bmatrix} 1 & 1 & 1 & 1 & 1 \\ -1 & -1 & -1 & 0 & 2 \\ 0 & 0 & 0 & 0 & 0 \end{bmatrix} \begin{matrix} R_1 + R_2 \to R_2 \\ R_1 + R_3 \to R_3 \end{matrix} \begin{bmatrix} \boxed{1} & 1 & 1 & 1 & 1 \\ 0 & 0 & 0 & \boxed{1} & 3 \\ 0 & 0 & 0 & 0 & 0 \end{bmatrix}$$

This is a Semi-Reduced Matrix. To find the row space of matrix $B$, we just need to write the $Span$ of independent rows. Thus, $Row(B) = Span\{(1,1,1,1,1),(0,0,0,1,3)\}$.

**Part b:** To find the column space of $B$, we need to change matrix $B$ to the Semi-Reduced Matrix. We already did that in part a. Now, we need to locate the columns in the Semi-Reduced Matrix of $B$ that contain the leaders, and then we should locate them to the original matrix $B$.

$$\begin{array}{c}\downarrow \qquad\quad \downarrow \\ \begin{bmatrix} \boxed{1} & 1 & 1 & 1 & 1 \\ 0 & 0 & 0 & \boxed{1} & 3 \\ 0 & 0 & 0 & 0 & 0 \end{bmatrix}\end{array}$$ Semi-Reduced Matrix

$$\begin{array}{c}\downarrow \qquad\qquad\quad \downarrow \\ \begin{bmatrix} 1 & 1 & 1 & 1 & 1 \\ -1 & -1 & -1 & 0 & 2 \\ 0 & 0 & 0 & 0 & 0 \end{bmatrix}\end{array}$$ Matrix $B$

Each remaining columns is a linear combination of the first and fourth columns.

Thus, $Column(B) = Span\{(1,-1,0),(1,0,0)\}$.

**Part c:** To find the rank of matrix $B$, we just need to change matrix $A$ to the Semi-Reduced Matrix. We already did that in part a. Thus,

$Rank(A) = dim\big(Row(B)\big) = dim\big(Column(B)\big) = 2$.

# C.2 Linear Transformation

We start this section with an introduction to polynomials, and we explain how they are similar to $\mathbb{R}^n$ as vector spaces. At the end of this section we discuss a new concept called linear transformation.

Before discussing polynomials, we need to know the following mathematical facts:

**Fact C.2.1** $\mathbb{R}^{n \times m} = \mathbb{R}_{n \times m} = M_{n \times m}(\mathbb{R})$ is a vector space.

**Fact C.2.2** $\mathbb{R}^{2 \times 3}$ is equivalent to $\mathbb{R}^6$ as a vector space.

(i.e. $\begin{bmatrix} 1 & 2 & 3 \\ 0 & 1 & 1 \end{bmatrix}$ is equivalent to (1,2,3,0,1,1) ).

**Fact C.2.3** $\mathbb{R}^{3 \times 2}$ is equivalent to $\mathbb{R}^6$ as a vector space.

(i.e. $\begin{bmatrix} 1 & 2 \\ 3 & 0 \\ 1 & 1 \end{bmatrix}$ is equivalent to (1,2,3,0,1,1) ).

After knowing the above facts, we introduce polynomials as follows:

$P_n = Set\ of\ all\ polynomials\ of\ degree < n.$

The algebraic expression of polynomials is in the following from: $a_n x^n + a_{n-1} x^{n-1} + \cdots + a_1 x^1 + a_0$

$a_n, a_{n-1}$ and $a_1$ are coefficients.

$n$ and $n - 1$ are exponents that must be positive integers whole numbers.

$a_0$ is a constant term.

The degree of polynomial is determined by the highest power (exponent).

We list the following examples of polynomials:

- $P_2 = Set\ of\ all\ polynomials\ of\ degree < 2$ (i.e. $3x + 2 \in P_2$ , $0 \in P_2$, $10 \in P_2$, $\sqrt{3} \in P_2$ but $\sqrt{3}\sqrt{x} \notin P_2$).
- $P_4 = Set\ of\ all\ polynomials\ of\ degree < 4$ (i.e. $31x^2 + 4 \in P_4$).
- If $P(x) = 3$, then $deg\big(P(x)\big) = 0$.
- $\sqrt{x} + 3$ is not a polynomial.

**Result C.2.1** $P_n$ is a vector space.

**Fact C.2.4** $\mathbb{R}^{2\times 3} = M_{2\times 3}(\mathbb{R})$ as a vector space same as $\mathbb{R}^6$.

**Result C.2.2** $P_n$ is a vector space, and it is the same as $\mathbb{R}^n$. (i.e. $a_0 + a_1 x^1 + \cdots + a_{n-1} x^{n-1} \leftrightarrow (a_0, a_1, \dots, a_{n-1})$. Note: The above form is in an ascending order.

**Result C.2.3** $dim(P_n) = n$.

**Fact C.2.5** $P_3 = Span\{3\ Independent\ Polynomials, and\ Each\ of\ Degree < 3\}$. (i.e. $P_3 = Span\{1, x, x^2\}$).

**Example C.2.1** Given the following polynomials:

$3x^2 - 2, -5x, 6x^2 - 10x - 4$.

   a.  Are these polynomials independent?

   b.  Let $D = Span\{3x^2 - 2, -5x, 6x^2 - 10x - 4\}$. Find a basis for $D$.

**Solution: Part a:** We know that these polynomial live in $P_3$, and as a vector space $P_3$ is the same as $\mathbb{R}^3$. According to result 3.2.2, we need to make each polynomial equivalent to $\mathbb{R}^n$ as follows:

$3x^2 - 2 = -2 + 0x + 3x^2 \leftrightarrow (-2, 0, 3)$

$-5x = 0 - 5x + 0x^2 \leftrightarrow (0, -5, 0)$

$6x^2 - 10x - 4 = -4 - 10x + 6x^2 \leftrightarrow (-4, -10, 6)$

Now, we need to write these vectors as a matrix.

$\begin{bmatrix} -2 & 0 & 3 \\ 0 & -5 & 0 \\ -4 & -10 & 6 \end{bmatrix}$ Each point is a row-operation. We need to reduce this matrix to Semi-Reduced Matrix.

Then, we apply the Row-Reduction Method to get the Semi-Reduced Matrix as follows:

$$\begin{bmatrix} -2 & 0 & 3 \\ 0 & -5 & 0 \\ -4 & -10 & 6 \end{bmatrix} -2R_1 + R_3 \rightarrow R_3 \begin{bmatrix} -2 & 0 & 3 \\ 0 & -5 & 0 \\ 0 & -10 & 0 \end{bmatrix}$$

$$-2R_2 + R_3 \rightarrow R_3 \begin{bmatrix} -2 & 0 & 3 \\ 0 & -5 & 0 \\ 0 & 0 & 0 \end{bmatrix}$$ This is a Semi-Reduced Matrix.

Since there is a zero-row in the Semi-Reduced Matrix, then these elements are dependent. Thus, the answer to this question is NO.

**Part b:** Since there are only 2 vectors survived after checking for dependency in part a, then the basis for $(0, -5,0) \leftrightarrow -5x$.

**Result C.2.4** Given $v_1, v_2, \ldots, v_k$ points in $\mathbb{R}^n$ where $k < n$. Choose one particular point, say $Q$, such that $Q = c_1 v_1 + c_2 v_2 + \cdots + c_k v_k$ where $c_1, c_2, \ldots, c_k$ are constants. If $c_1, c_2, \ldots, c_k$ are unique, then $v_1, v_2, \ldots, v_k$ are independent.

Note: The word "unique" in result 3.2.4 means that there is only one value for each of $c_1, c_2, \ldots, c_k$.

**Proof of Result C.2.4** By using proof by contradiction, we assume that $v_1 = \alpha_2 v_2 + \alpha_3 v_3 + \cdots + \alpha_k v_k$ where $\alpha_2, \alpha_3, \ldots, \alpha_k$ are constants. Our assumption means that it is dependent. Using algebra, we obtain:

$Q = c_1\alpha_2 v_2 + c_1\alpha_3 v_3 + \cdots + c_1\alpha_k v_k + c_2 v_2 + \cdots + c_k v_k$.

$Q = (c_1\alpha_2 + c_2)v_2 + (c_1\alpha_3 + c_3)v_3 + \cdots + (c_1\alpha_k + c_k)v_k + 0v_1$. Thus, none of them is a linear combination of the others which means that they are linearly independent. This is a contradiction. Therefore, our assumption that $v_1, v_2, \ldots,$ and $v_k$ were linearly dependent is false. Hence, $v_1, v_2, \ldots,$ and $v_k$ are linearly independent.                                   ◻

**Result C.2.5** Assume $v_1, v_2, \ldots, v_k$ are independent and $Q \in Span\{v_1, v_2, \ldots, v_k\}$. Then, there exists unique number $c_1, c_2, \ldots, c_k$ such that $Q = c_1 v_1 + c_2 v_2 + \cdots + c_k v_k$.

**Linear Transformation:**

**Definition C.2.1** $T: V \to W$ where $V$ is a domain and $W$ is a co-domain. $T$ is a linear transformation if for every $v_1, v_2 \in V$ and $\alpha \in \mathbb{R}$, we have the following: $T(\alpha v_1 + v_2) = \alpha T(v_1) + T(v_2)$.

**Example C.2.2** Given $T: \mathbb{R}^2 \to \mathbb{R}^3$ where $\mathbb{R}^2$ is a domain and $\mathbb{R}^3$ is a co-domain. $T((a_1, a_2)) = (3a_1 + a_2, a_2, -a_1)$.

  a. Find $T((1,1))$.
  b. Find $T((1,0))$.
  c. Show that $T$ is a linear transformation.

**Solution: Part a:** Since $T((a_1, a_2)) = (3a_1 + a_2, a_2, -a_1)$, then $a_1 = a_2 = 1$. Thus, $T((1,1)) = (3(1) + 1, 1, -1) = (4, 1, -1)$.

**Part b:** Since $T((a_1, a_2)) = (3a_1 + a_2, a_2, -a_1)$, then $a_1 = 1$ and $a_2 = 0$. Thus, $T((1,0)) = (3(1) + 0, 0, -1) = (3, 0, -1)$.

**Part c:** Proof: We assume that $v_1 = (a_1, a_2)$,

$v_2 = (b_1, b_2)$, and $\alpha \in \mathbb{R}$. We will show that $T$ is a linear transformation. Using algebra, we start from the Left-Hand-Side (LHS):

$\alpha v_1 + v_2 = (\alpha a_1 + b_1, \alpha a_2 + b_2)$

$T(\alpha v_1 + v_2) = T((\alpha a_1 + b_1, \alpha a_2 + b_2))$
$T(\alpha v_1 + v_2) = (3\alpha a_1 + 3b_1 + \alpha a_2 + b_2, \alpha a_2 + b_2, -\alpha a_1 - b_1)$

Now, we start from the Right-Hand-Side (RHS):

$\alpha T(v_1) + T(v_2) = \alpha T(a_1, a_2) + T(b_1, b_2)$
$\alpha T(v_1) + T(v_2) = \alpha(3a_1 + a_2, a_2, -a_1) + (3b_1 + b_2, b_2, -b_1)$
$= (3\alpha a_1 + \alpha a_2, \alpha a_2, -\alpha a_1) + (3b_1 + b_2, b_2, -b_1)$
$= (3\alpha a_1 + \alpha a_2 + 3b_1 + b_2, \alpha a_2 + b_2, -\alpha a_1 - b_1)$

Thus, $T$ is a linear transformation. $\qquad\qquad\square$

**Result C.2.6** Given $T: \mathbb{R}^n \to \mathbb{R}^m$. Then, $T((a_1, a_2, a_3, \ldots, a_n)) =$ Each coordinate is a linear combination of the $a_i's$.

**Example C.2.3** Given $T: \mathbb{R}^3 \to \mathbb{R}^4$ where $\mathbb{R}^3$ is a domain and $\mathbb{R}^4$ is a co-domain.

 a. If $T\big((x_1, x_2, x_3)\big) = (-3x_3 + 6x_1, -10x_2, 13, -x_3)$, is $T$ a linear transformation?
 b. If $T\big((x_1, x_2, x_3)\big) = (-3x_3 + 6x_1, -10x_2, 0, -x_3)$, is $T$ a linear transformation?

**Solution: Part a:** Since 13 is not a linear combination of $x_1, x_2$ and $x_3$. Thus, $T$ is not a linear transformation.

**Part b:** Since 0 is a linear combination of $x_1, x_2$ and $x_3$. Thus, $T$ is a linear transformation.

**Example C.2.4** Given $T: \mathbb{R}^2 \to \mathbb{R}^3$ where $\mathbb{R}^2$ is a domain and $\mathbb{R}^3$ is a co-domain. If $T\big((a_1, a_2)\big) = (a_1{}^2 + a_2, -a_2)$, is $T$ a linear transformation?

**Solution:** Since $a_1{}^2 + a_2$ is not a linear combination of $a_1$ and $a_2$. Hence, $T$ is not a linear transformation.

**Example C.2.5** Given $T: \mathbb{R} \to \mathbb{R}$. If $T(x) = 10x$, is $T$ a linear transformation?

**Solution:** Since **it** is a linear combination of $a_1$ such that $\alpha a_1 = 10x$. Hence, $T$ is a linear transformation.

**Example C.2.6** Find the standard basis for $\mathbb{R}^2$.

**Solution:** The standard basis for $\mathbb{R}^2$ is the rows of $I_2$.

Since $I_2 = \begin{bmatrix} 1 & 0 \\ 0 & 1 \end{bmatrix}$, then the standard basis for $\mathbb{R}^2$ is $\{(1,0), (0,1)\}$.

**Example C.2.7** Find the standard basis for $\mathbb{R}^3$.

**Solution:** The standard basis for $\mathbb{R}^3$ is the rows of $I_3$.

Since $I_3 = \begin{bmatrix} 1 & 0 & 0 \\ 0 & 1 & 0 \\ 0 & 0 & 1 \end{bmatrix}$, then the standard basis for $\mathbb{R}^3$ is $\{(1,0,0), (0,1,0), (0,0,1)\}$.

**Example C.2.8** Find the standard basis for $P_3$.

**Solution:** The standard basis for $P_3$ is $\{1, x, x^2\}$.

**Example C.2.9** Find the standard basis for $P_4$.

**Solution:** The standard basis for $P_4$ is $\{1, x, x^2, x^3\}$.

**Example C.2.10** Find the standard basis for $\mathbb{R}_{2\times2} = M_{2\times2}(\mathbb{R})$.

**Solution:** The standard basis for $\mathbb{R}_{2\times2} = M_{2\times2}(\mathbb{R})$ is $\{\begin{bmatrix} 1 & 0 \\ 0 & 0 \end{bmatrix}, \begin{bmatrix} 0 & 1 \\ 0 & 0 \end{bmatrix}, \begin{bmatrix} 0 & 0 \\ 1 & 0 \end{bmatrix}, \begin{bmatrix} 0 & 0 \\ 0 & 1 \end{bmatrix}\}$ because $\mathbb{R}_{2\times2} = M_{2\times2}(\mathbb{R}) = \mathbb{R}^4$ as a vector space where standard basis for $\mathbb{R}_{2\times2} = M_{2\times2}(\mathbb{R})$ is the rows of $I_4 = \begin{bmatrix} 1 & 0 & 0 & 0 \\ 0 & 1 & 0 & 0 \\ 0 & 0 & 1 & 0 \\ 0 & 0 & 0 & 1 \end{bmatrix}$ that are represented by $2 \times 2$ matrices.

**Example C.2.11** Let $T \colon \mathbb{R}^2 \to \mathbb{R}^3$ be a linear transformation such that
$T(2,0) = (0,1,4)$

$T(-1,1) = (2,1,5)$

Find $T(3,5)$.

**Solution:** The given points are $(2,0)$ and $(-1,1)$. These two points are independent because of the following:

$\begin{bmatrix} 2 & 0 \\ -1 & 1 \end{bmatrix} \frac{1}{2}R_1 + R_2 \to R_2 \begin{bmatrix} 2 & 0 \\ 0 & 1 \end{bmatrix}$

Every point in $\mathbb{R}^2$ is a linear combination of $(2,0)$ and $(-1,1)$. There exists unique numbers $c_1$ and $c_2$ such that $(3,5) = c_1(2,0) + c_2(-1,1)$.
$3 = 2c_1 - c_2$
$5 = c_2$
Now, we substitute $c_2 = 5$ in $3 = 2c_1 - c_2$, we obtain:
$3 = 2c_1 - 5$
$c_1 = 4$
Hence, $(3,5) = 4(2,0) + 5(-1,1)$.
$T(3,5) = T(4(2,0) + 5(-1,1))$

$T(3,5) = 4T(2,0) + 5T(-1,1)$
$T(3,5) = 4(0,1,4) + 5(2,1,5) = (10,9,41)$
Thus, $T(3,5) = (10,9,41)$.

**Example C.2.12** Let $T: \mathbb{R} \to \mathbb{R}$ be a linear transformation such that $T(1) = 3$. Find $T(5)$.

**Solution:** Since it is a linear transformation, then $T(5) = T(5 \cdot 1) = 5T(1) = 5(3) = 15$. If it is not a linear transformation, then it is impossible to find $T(5)$.

# C.3 Kernel and Range

In this section, we discuss how to find the standard matrix representation, and we give examples of how to find kernel and range.

**Definition C.3.1** Given $T: \mathbb{R}^n \to \mathbb{R}^m$ where $\mathbb{R}^n$ is a domain and $\mathbb{R}^m$ is a co-domain. Then, Standard Matrix Representation is a $m \times n$ matrix. This means that it is $dim(Co - Domain) \times dim(Domain)$ matrix.

**Definition C.3.2** Given $T: \mathbb{R}^n \to \mathbb{R}^m$ where $\mathbb{R}^n$ is a domain and $\mathbb{R}^m$ is a co-domain. Kernel is a set of all points in the domain that have image which equals to the origin point, and it is denoted by $Ker(T)$. This means that $Ker(T) = Null\ Space\ of\ T$.

**Definition C.3.3** Range is the column space of standard matrix representation, and it is denoted by $Range(T)$.

**Example C.3.1** Given $T: \mathbb{R}^3 \to \mathbb{R}^4$ where $\mathbb{R}^3$ is a domain and $\mathbb{R}^4$ is a co-domain.

$T\big((x_1, x_2, x_3)\big) = (-5x_1, 2x_2 + x_3, -x_1, 0)$

    a.  Find the Standard Matrix Representation.
    b.  Find $T((3,2,1))$.
    c.  Find $Ker(T)$.
    d.  Find $Range(T)$.

**Solution: Part a:** According to definition 3.3.1, the Standard Matrix Representation, let's call it $M$, here is

$4 \times 3$. We know from section 3.2 that the standard basis for domain (here is $\mathbb{R}^3$) is $\{(1,0,0),(0,1,0),(0,0,1)\}$. We assume the following:

$$v_1 = (1,0,0)$$
$$v_2 = (0,1,0)$$
$$v_3 = (0,0,1)$$

Now, we substitute each point of the standard basis for domain in $T\big((x_1, x_2, x_3)\big) = (-5x_1, 2x_2 + x_3, -x_1, 0)$ as follows:

$T\big((1,0,0)\big) = (-5,0,-1,0)$
$T\big((0,1,0)\big) = (0,2,0,0)$
$T\big((0,0,1)\big) = (0,1,0,0)$

Our goal is to find $M$ so that $T\big((x_1, x_2, x_3)\big) = M \begin{bmatrix} x_1 \\ x_2 \\ x_3 \end{bmatrix}$.

$M = \begin{bmatrix} -5 & 0 & 0 \\ 0 & 2 & 1 \\ -1 & 0 & 0 \\ 0 & 0 & 0 \end{bmatrix}$ This is the Standard Matrix

Representation. The first, second and third columns represent $T(v_1), T(v_2)$ and $T(v_3)$.

**Part b:** Since $\big((x_1, x_2, x_3)\big) = M \begin{bmatrix} x_1 \\ x_2 \\ x_3 \end{bmatrix}$, then

$$T\big((3,2,1)\big) = \begin{bmatrix} -5 & 0 & 0 \\ 0 & 2 & 1 \\ -1 & 0 & 0 \\ 0 & 0 & 0 \end{bmatrix} \begin{bmatrix} 3 \\ 2 \\ 1 \end{bmatrix}$$

$$T\big((3,2,1)\big) = 3 \cdot \begin{bmatrix} -5 \\ 0 \\ -1 \\ 0 \end{bmatrix} + 2 \cdot \begin{bmatrix} 0 \\ 2 \\ 0 \\ 0 \end{bmatrix} + 1 \cdot \begin{bmatrix} 0 \\ 1 \\ 0 \\ 0 \end{bmatrix} = \begin{bmatrix} -15 \\ 5 \\ -3 \\ 0 \end{bmatrix}$$

$\begin{bmatrix} -15 \\ 5 \\ -3 \\ 0 \end{bmatrix}$ is equivalent to $(-15,5,-3,0)$. This lives in the

co-domain. Thus, $T\big((3,2,1)\big) = (-15,5,-3,0)$.

**Part c:** According to definition 3.3.2, $Ker(T)$ is a set of all points in the domain that have image= $(0,0,0,0)$. Hence, $T\big((x_1,x_2,x_3)\big) = (0,0,0,0)$. This means the following: $M \begin{bmatrix} x_1 \\ x_2 \\ x_3 \end{bmatrix} = \begin{bmatrix} 0 \\ 0 \\ 0 \\ 0 \end{bmatrix}$

$$\begin{bmatrix} -5 & 0 & 0 \\ 0 & 2 & 1 \\ -1 & 0 & 0 \\ 0 & 0 & 0 \end{bmatrix} \begin{bmatrix} x_1 \\ x_2 \\ x_3 \end{bmatrix} = \begin{bmatrix} 0 \\ 0 \\ 0 \\ 0 \end{bmatrix}$$

Since $Ker(T) = Null(M)$, then we need to find $N(M)$ as follows:

$$\left(\begin{array}{ccc|c} -5 & 0 & 0 & 0 \\ 0 & 2 & 1 & 0 \\ -1 & 0 & 0 & 0 \\ 0 & 0 & 0 & 0 \end{array}\right) -\frac{1}{5}R_1 \left(\begin{array}{ccc|c} 1 & 0 & 0 & 0 \\ 0 & 2 & 1 & 0 \\ -1 & 0 & 0 & 0 \\ 0 & 0 & 0 & 0 \end{array}\right) R_1 + R_3 \to R_3$$

$$\begin{pmatrix} 1 & 0 & 0 & | & 0 \\ 0 & 2 & 1 & | & 0 \\ 0 & 0 & 0 & | & 0 \\ 0 & 0 & 0 & | & 0 \end{pmatrix} \frac{1}{2} R_2 \begin{pmatrix} 1 & 0 & 0 & | & 0 \\ 0 & 1 & 0.5 & | & 0 \\ 0 & 0 & 0 & | & 0 \\ 0 & 0 & 0 & | & 0 \end{pmatrix}$$ This is a Completely-

Reduced Matrix. Now, we need to read the above matrix as follows:

$x_1 = 0$

$x_2 + \frac{1}{2} x_3 = 0$

$0 = 0$

$0 = 0$

To write the solution, we need to assume that

$x_3 \in \mathbb{R}$ $(Free\ Variable)$.

Hence, $x_1 = 0$ and $x_2 = -\frac{1}{2} x_3$.

$N(M) = \{(0, -\frac{1}{2} x_3, x_3) | x_3 \in \mathbb{R}\}$.

By letting $x_3 = 1$, we obtain:

$Nullity(M) = Number\ of\ Free\ Variables = 1$, and

$Basis = \{(0, -\frac{1}{2}, 1)\}$

Thus, $Ker(T) = N(M) = Span\{(0, -\frac{1}{2}, 1)\}$.

**Part d:** According to definition 3.3.3, $Range(T)$ is the column space of $M$. Now, we need to locate the columns in the Completely-Reduced Matrix in part c that contain the leaders, and then we should locate them to the original matrix as follows:

$$\begin{pmatrix} \boxed{1} & 0 & 0 & | & 0 \\ 0 & \boxed{1} & 0.5 & | & 0 \\ 0 & 0 & 0 & | & 0 \\ 0 & 0 & 0 & | & 0 \end{pmatrix}$$ Completely-Reduced Matrix

$$\begin{pmatrix} -5 & 0 & 0 & 0 \\ 0 & 2 & 1 & 0 \\ -1 & 0 & 0 & 0 \\ 0 & 0 & 0 & 0 \end{pmatrix} \text{Orignial Matrix}$$

Thus, $Range(T) = Span\{(-5,0,-1,0),(0,2,0,0)\}$.

**Result C.3.1** Given $T: \mathbb{R}^n \to \mathbb{R}^m$ where $\mathbb{R}^n$ is a domain and $\mathbb{R}^m$ is a co-domain. Let $M$ be a standard matrix representation. Then,
$Range(T) = Span\{Independent\ Columns of\ M\}$.

**Result C.3.2** Given $T: \mathbb{R}^n \to \mathbb{R}^m$ where $\mathbb{R}^n$ is a domain and $\mathbb{R}^m$ is a co-domain. Let $M$ be a standard matrix representation. Then, $dim\big(Range(T)\big) = Rank(M) =$ Number of Independent Columns.

**Result C.3.3** Given $T: \mathbb{R}^n \to \mathbb{R}^m$ where $\mathbb{R}^n$ is a domain and $\mathbb{R}^m$ is a co-domain. Let $M$ be a standard matrix representation. Then, $dim\big(Ker(T)\big) = Nullity(M)$.

**Result C.3.4** Given $T: \mathbb{R}^n \to \mathbb{R}^m$ where $\mathbb{R}^n$ is a domain and $\mathbb{R}^m$ is a co-domain. Let $M$ be a standard matrix representation. Then, $dim\big(Range(T)\big) + dim(Ker(T)) = dim(Domain)$.

**Example C.3.2** Given $T: P_2 \to \mathbb{R}$. $T(f(x)) = \int_0^1 f(x)dx$ is a linear transformation.

    a. Find $T(2x - 1)$.
    b. Find $Ker(T)$.
    c. Find $Range(T)$.

**Solution:**

**Part a:** $T(2x - 1) = \int_0^1 (2x - 1)dx = x^2 - x \Big|_{\substack{x = 1 \\ x = 0}} = 0.$

**Part b:** To find $Ker(T)$, we set equation of $T = 0$, and $f(x) = a_0 + a_1 x \in P_2$.

Thus, $T(f(x)) = \int_0^1 (a_0 + a_1 x)dx = a_0 x + \frac{a_1}{2} x^2 \Big|_{\substack{x = 1 \\ x = 0}} = 0$

$a_0 + \frac{a_1}{2} - 0 = 0$

$a_0 = -\frac{a_1}{2}$

Hence, $Ker(T) = \{-\frac{a_1}{2} + a_1 x | a_1 \in \mathbb{R}\}$. We also know that $dim(Ker(T) = 1$ because there is one free variable. In addition, we can also find basis by letting $a_1$ be any real number not equal to zero, say $a_1 = 1$, as follows:

$Basis = \{-\frac{1}{2} + x\}$

Thus, $Ker(T) = Span\{-\frac{1}{2} + x\}$.

**Part c:** It is very easy to find range here. $Range(T) = \mathbb{R}$ because we linearly transform from a second degree polynomial to a real number. For example, if we linearly transform from a third degree polynomial to a second degree polynomial, then the range will be $P_2$.

*This page intentionally left blank*

# Answers to Odd-Numbered Exercises

## 1.7 Exercises

**1.** $\mathcal{L}^{-1}\left\{\frac{10}{(s-4)^4}\right\} = \frac{10}{6}e^{4x}x^3$

**3.** $\mathcal{L}^{-1}\left\{\frac{s+5}{(s+3)^4}\right\} = \frac{1}{2}x^2e^{-3x} + 6x^3e^{-3x}$

**5.** $\mathcal{L}^{-1}\left\{\frac{2}{s^2-6s+13}\right\} = e^{3x}\sin(2x)$

**7.** $y(x) = -4e^{-3x}$

**9.** $\mathcal{L}^{-1}\left\{\frac{4}{(s-1)^2(s+3)}\right\} = \frac{1}{4}e^x + xe^x + \frac{1}{4}e^{-3x}$

**11.** $\mathcal{L}\{U(x-2)e^{3x}\} = \frac{e^{6-2s}}{s-3}$

**13.** $\mathcal{L}^{-1}\left\{\frac{se^{-4x}}{s^2+4}\right\} = U(x-4)\cos(2x-8)$

**15.** Assume $W(x) = -\frac{3}{8} + \frac{1}{3}e^x + \frac{1}{24}e^{-8x}$. Then, we obtain:

$y(x) = W(x) - 3U(x-5)W(x-5) - 2U(x-5)W(x-5)$

**17.** $y(x) = \mathcal{L}^{-1}\left\{\frac{e^{-3x}}{s^2(s^2+1)}\right\}$

**19.** $\mathcal{L}^{-1}\left\{\frac{2s}{(s^2+4)^2}\right\} = \frac{4s}{(s^2+4)^2}$

**21.** $w(t) = 1 + t^2$ and $h(t) = 2t$

## 2.3 Exercises

**1.** $y_{homogenous}(x) = c_1e^{-x} + c_2xe^{-x}$ for some $c_1, c_2 \in \Re$

**3.** $y_{general}(x) = (c_1 + c_2x + c_3x^2 + c_4e^x) + (a_0 + a_1x + a_2x^2)x^3$ for some $c_1, c_2, c_3, c_4, a_0, a_1, a_2 \in \Re$

**5.** $y_{general}(x) = (c_1 + c_2 x + c_3 e^x) + 0.05 \sin(2x) +$

$0.1\cos(2x)$ for some $c_1, c_2, c_3 \in \Re$

**7.** It is impossible to describe $y_{particular}(x)$

# 3.3 Exercises

**1.** $y_{homogenous}(x) = e^{-\frac{1}{2}x}\left[c_1 \cos\left(\frac{\sqrt{15}}{2}x\right) + c_2 \sin\left(\frac{\sqrt{15}}{2}x\right)\right]$for

some $c_1, c_2 \in \Re$

**3.** $y_{homogenous}(x) = c_1 + c_2 x^{\frac{3}{2}} \cos\left(\frac{\sqrt{3}}{2}\ln(x)\right) +$

$c_3 x^{\frac{3}{2}} \sin\left(\frac{\sqrt{3}}{2}\ln(x)\right)$for some $c_1, c_2, c_3 \in \Re$

**5.** It is impossible to use Cauchy-Euler Method
because the degrees of $y'$ and $y''$ are not equal to each
other when you substitute them in the given
differential equation.

# 4.6 Exercises

**1.** $y(x) = \sqrt[3]{(x + 1)^3 + (x + 1)ce^{-3x}}$

**3.** $tan^{-1}(y) - tan^{-1}(x) = c$ for some $c \in \Re$

**5.** $-\frac{1}{5}e^{-5y} - 3xe^x + 3e^x = c$ for some $c \in \Re$

**7.** $\ln|\sin(5x + y)| - 2(5x + y) - x = c$ for some $c \in \Re$

# Index

## A

Applications of Differential Equations, 96

## B

Basis, 123
Bernoulli Method, 78

## C

Cauchy-Euler Method, 74
Constant Coefficients, 51
Cramer's Rule, 46

## D

Water Tank Application, 104

*This page intentionally left blank*